**Spatial-temporal risk field-based coupled dynamic-static driving risk assessment and trajectory planning in weaving segments**

Guodong Ma[a], Baofeng Sun[a,*], Hongchao Liang[a], Wenyu Yang[a], Huxing Zhou[a]

a. School of Transportation, Jilin University, Changchun 130022, China

## Abstract

As connected and automated vehicles (CAVs) gradually penetrate the existing transportation system, the inherent turbulence within weaving segments is expected to be mitigated through CAV technologies. However, despite making progress in capturing static and dynamic risk factors, traditional CAV technologies still lack foreseeability for dynamic risks. This leads to suboptimal results in trajectory planning, thereby hindering the maximization of expected benefits. To fill these gaps, we first propose a spatial-temporal coupled risk assessment paradigm by constructing a three-dimensional spatial-temporal risk field (STRF). Specifically, we introduce spatial-temporal distances to quantify the impact of future trajectories of dynamic obstacles. We also incorporate a geometrically configured specialized field for weaving segments to constrain vehicle movement directionally. To enhance the STRF's accuracy, we further developed a parameter calibration method using real-world aerial video data, leveraging YOLO-based machine vision and dynamic risk balance theory. A comparative analysis with traditional risk field shows that the STRF possesses superior risk foreseeability. Building on these results, we final design a STRF-based CAV trajectory planning method in weaving segments. We integrate spatial-temporal risk occupancy maps, dynamic iterative sampling, and quadratic programming to enhance safety, comfort, and efficiency. By incorporating both dynamic and static risk factors during the sampling phase, our method ensures robust safety performance. Additionally, the proposed method simultaneously optimizes path and speed using a parallel computing approach, reducing computation time. Real-world cases show that, compared to the baseline schemes, and real human driving trajectories, our method significantly improves safety, reduces lane-change completion time, and minimizes speed fluctuations.

* Corresponding author at: No. 5988 Renmin Street, Changchun 130022, China.
E-mail address: magd22@mails.jlu.edu.cn (Guodong Ma), sunbf@jlu.edu.cn (Baofeng Sun), lianghc0312@163.com (Hongchao Liang), yangwy24@mails.jlu.edu.cn (Wenyu Yang), zhouhx@jlu.edu.cn (Huxing Zhou).

## Introduction

Expressways serve as the primary arteries of urban transportation networks, accommodating high-intensity traffic flows (Yuan et al., 2024). Weaving segments are a critical component of expressways, facilitating vehicle entry and exit via on-ramps and off-ramps (Ouyang et al., 2023). Traffic flow in these zones is characterized by three primary interactions: on-ramp merging, off-ramp diverging, and main-line through-traffic. The frequent lane-changing behaviors among these three traffic flows lead to complex interactions, making the weaving segments extremely prone to turbulence effects. (van Beinum et al., 2018). As a result, weaving segments emerge as critical bottlenecks with high accident risks and severe congestion (Chen and Ahn, 2018; Ouyang et al., 2023).

In recent years, the development of connected and automated vehicle (CAV) technology has provided a highly promising solution for reducing driving risks in weaving segments through enhanced perception, decision-making, and control capabilities (Chen et al., 2021a; Papadoulis et al., 2019). However, the prolonged coexistence of CAVs and human-driven vehicles (HDVs) in mixed traffic environments poses critical challenges, hindering the full utilization of CAV technological advantages. First, the inherent unpredictability of traditional traffic flow (pure HDVs) and the additional variance in driving behaviors within mixed traffic environments significantly diminish the effectiveness of traditional driving risk assessment methods in high-speed, dynamic, and multisource risk scenarios. Traditional methods relying on two-dimensional static fields fail to capture risk predictability, thereby ignoring the impact of the future trajectories of high-speed dynamic obstacles on current risk assessment (Han et al., 2023). These models typically determine the risk level based solely on the state of risk sources at the observation moment, leading to delayed responses to emergent events and inaccurate safety redundancy designs. Furthermore, the presence of multisource risks in weaving segments presents a major challenge to existing trajectory planning methods. Existing methods primarily consider static obstacles and overly simplified dynamic risk assessments (Yao and Sun, 2025), rendering them inadequate for mitigating dynamic obstacle risks in high-speed environments.

Second, with the advancement of autonomous driving technology, trajectory planning must comprehensively consider the coordinated optimization of multiple metrics, including safety, efficiency, comfort, and real-time performance. Among these, safety and real-time performance are particularly crucial, especially in high-interaction scenarios such as weaving segments that involve multiple intertwined traffic streams. While numerous previous studies have extensively explored the tradeoff among factors like safety and efficiency in trajectory planning, and despite achieving favorable progress in conventional scenarios, these conventional models face significant challenges when applied to urban expressway weaving segments where safety demands are highly pronounced. Real-time performance is equally critical in weaving segments; some studies establish 100 milliseconds as the time threshold for trajectory planning (Liu et al., 2025). The intense computational load caused by multivehicle interactions in weaving segments causes some models to experience computation timeouts under such environments. The unique geometric characteristics of weaving segments lead to a high degree of coupling between merging and diverging behaviors, making multivehicle interactions much more prominent and thus imposing far stricter requirements for safety and real-time performance than conventional road segments. Under these specialized scenarios, prioritizing the hard constraints of safety and real-time performance in a differentiated manner before considering other objectives has emerged as a critical challenge requiring immediate resolution. This necessitates that current research must achieve a multi-objective tradeoff specifically oriented toward the complex scenarios of merging and diverging interactions within weaving segments.

To address the aforementioned theoretical gaps and challenges specific to weaving segments, this study conducts an in-depth exploration from two dimensions: precise risk quantification and highly efficient multi-objective trajectory planning. First, at the safety assessment level, this study proposes an innovative spatiotemporal coupled risk assessment paradigm by constructing a three-dimensional spatiotemporal risk field (STRF) model. This model overcomes the limitations of traditional two-dimensional static fields by establishing a rigorous correlation between temporal and spatial dimensions, and it proactively quantifies the potential hazards of the future trajectories of dynamic obstacles using a spatiotemporal distance metric. Concurrently, the model deeply integrates real-world aerial trajectory data from weaving segments, extracted via YOLO machine vision, for parameter calibration, thereby achieving high fidelity micro level risk evolution modeling. Second, at the trajectory planning level, to address the specialized scenario of weaving segments, this paper introduces an STRF-based trajectory planning framework. This framework integrates a spatiotemporal risk occupancy map, dynamic iterative sampling, and quadratic programming. It is specifically tailored to resolve the multi-objective tradeoff challenge inherent in the complex merging and diverging interactions of weaving segments. By establishing safety and real-time execution as strict priorities, the framework successfully balances competing driving metrics. Notably, it enables the parallel optimization of path and speed. This architectural advantage significantly reduces computational load while guaranteeing physical safety via the risk field, fully satisfying the stringent time constraints of high-frequency multivehicle interactions. Ultimately, this research delivers a novel integrated solution for risk assessment and trajectory planning, ensuring the safe, efficient, and comfortable operation of CAVs within mixed traffic flows.

## 2 Literature review and main contributions

The literature review will be centered on two major challenges that contribute to congestion and safety issues in weaving segments: surrogate safety measures (SSMs) for driving risk assessment in weaving segments and the research trends of trajectory planning in weaving segments.

### 2.1 SSMs for driving risk assessment in weaving segments

Weaving segments are roadways with a high frequency of crashes, making risk research in these areas particularly valuable. Since crashes are relatively rare events, assessing road risk solely based on crash frequency is often infeasible, especially in the absence of crash data. Further, the CAV risk assessment needs to be integrated into the same framework as the subsequent trajectory planning exercise. Therefore, an alternative methodology for risk assessment is required (Zhang et al., 2023). Previous studies have addressed this issue by analyzing traffic conflicts using simulation models with SSMs. Current SSMs are based on three types of deterministic approaches: time-based, distance-based, and acceleration-based. Time-based SSM includes time to collision (TTC) (Zhang et al., 2022), time headway (THW), and post encroachment time (PET) (Howlader et al., 2024). Distance-based SSM primarily involves the minimum safe distance (MSD) (Winkler et al., 2016), while acceleration-based SSM includes the deceleration rate for avoiding collision (DRAC) (Zhang et al., 2022). However, these methods exhibit several limitations in practical applications:

(1) They assume vehicle states remain constant over a given time period. While these methods are computationally efficient, they fail to account for uncertainties in vehicle motion and dynamic external factors (Brechtel et al., 2014; de Gelder et al., 2023; Li et al., 2023b), which can lead to misjudgments.
(2) They focus on rear-end risks, neglecting lateral and omnidirectional risks (Joo et al., 2023; Ma et al., 2025a).

(3) They are suitable for simple scenarios with a single risk source but often underestimate risks in complex scenarios with multiple, diverse risk sources, limiting their applicability (Lu et al., 2021; Ma et al., 2025a; Xiong et al., 2023).

(4) They rely on predefined thresholds, making them sensitive to minor threshold variations (Chen et al., 2021b).

Field-based risk assessment methods can address these limitations by accounting for vehicle uncertainties, multi-source risks, and their superposition, providing a unified and effective metric for multi-directional risks (Jiang et al., 2025; Li et al., 2025; Zuo et al., 2025). Given the diverse nature and characteristics of real transportation problems in traffic flow scenarios, the construction of the risk field must adapt to varying scenarios and issues for extended studies. From **Table 1**, the modeling and improvement of risk fields are primarily addressed through four key aspects.

(1) **Expanding the perceptual factors of the risk field:** The traditional risk field is often symmetric. This symmetry arises from considering only the distance to obstacles, simplifying computations but compromising accuracy. Correlation improvements accurately characterize the superimposed effects of multiple factors on the distribution of risk fields.

(2) **Expanding the perceptual dimensions of the risk field:** Most existing risk fields focus on two-dimensional spatial effects, considering field strength at a single moment and neglecting future trajectories. By incorporating predicted trajectories into current risk field modeling, the enhanced model captures future vehicle movements and better predicts potential risks (Chen et al., 2024b; Han et al., 2023).

(3) **Describing differences in the subjective risk field:** The heterogeneity of vehicle driving behavior manifests in both risk generation and observation. In the risk generation phase, trajectory prediction uncertainty, particularly for drivers, increases risk. Studies on risk field modeling in mixed CAV-HDV traffic showed that behavioral differences between CAVs and HDVs result in distinct risk values (Sun et al., 2023). In the risk observation phase, subjective risk fields differ from objective ones: the objective field aggregates all risks, whereas the subjective field depends on the observer's attributes, potentially leading to missed or inaccurate risk assessments (Huang et al., 2019; Song et al., 2024).

(4) **Designing a special field for roads with different geometrical organizational structures:** On standard roads, environmental modeling typically considers road lines and boundaries. However, complex traffic scenarios, such as intersections, weaving areas, and tunnels, demand specialized fields to account for unique environmental factors like traffic lights, merging zones, and lighting variations. These specialized fields are essential for accurate risk quantification in such environments (Chen and Wen, 2022; Yan et al., 2022; Zong et al., 2022).

**Table 1** Some improvements of the risk field.

| Directions | Ref. | Thrust |
|---|---|---|
| Expand the perceptual factors. | (Li et al., 2020) | Introduce acceleration. |
| | (Yan et al., 2022) | Introduce road curvature |
| Expand the perceptual dimensions. | (Han et al., 2023) | Reconstruct a spatial-temporal 3D risk field |
| | (Wang et al., 2024b) | Decouple the spatial dimensionality into vertical and horizontal dimensions. |
| Describe differences in subject vehicles. | (Xia et al., 2024) | Introduce driver perception and steering characteristics. |
| | (Sun et al., 2023) | Introduce driver trust in CAVs. |
| | (Song et al., 2024) | Simulate subjective risk perceptions of drivers. |
| | (Sarvesh et al., 2021) | Verify the correlation between CAV and HDV risk perception. |
| Design special fields | (Zong et al., 2022) | Establish a risk field at signalized intersections |
| | (Wang et al., 2024a) | Establish a risk field at highway merging areas |
| | (Shen et al., 2024) | Establish a risk field in the tunnel |
| | (Sun et al., 2023) | Establish a risk field at weaving segments |

2.2 Trajectory planning for complex weaving segment scenarios

Weaving segments in expressway traffic are more complex than basic, merging, and diverging segments (Rim et al., 2023). Merging and diverging segments involve two interacting traffic flows, while weaving segments (as shown in **Fig. 1**) involve three. These three interacting traffic flows lead to frequent lane changes, which ultimately result in the emergence of the turbulence effect. A weaving segment forms when an on-ramp directly connects to an off-ramp, combining merging and diverging characteristics with coupling effects. Although studies on merging and diverging focus on order optimization and trajectory planning (Chen et al., 2024a; Gu et al., 2024; Zhu et al., 2024), weaving segments introduce added complexity due to inflow-outflow interactions. Simply superimposing findings from merging and diverging studies is insufficient, yet research on trajectory planning within weaving segments remains limited.

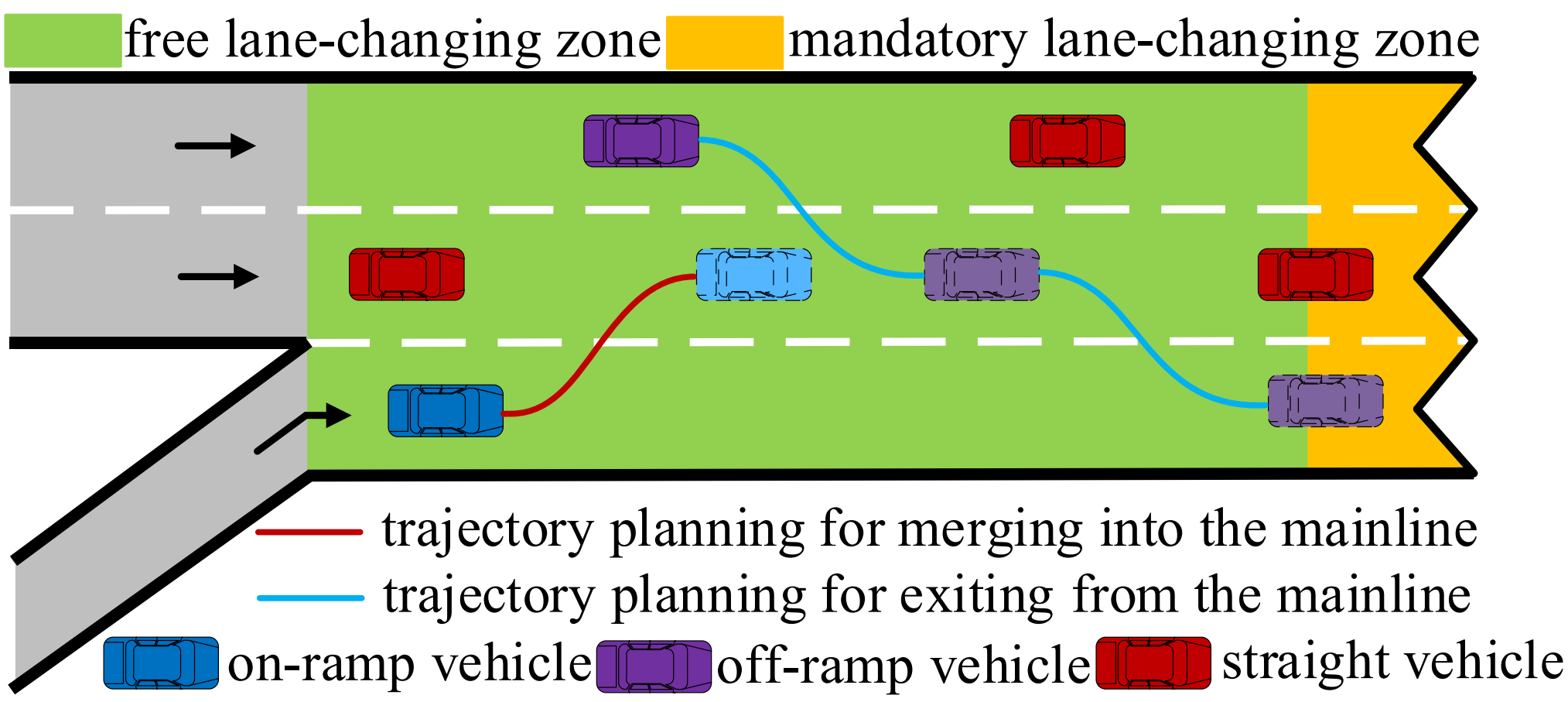

**Fig. 1** The traffic composition and the trajectory planning in weaving segments.

Existing trajectory planning methods in weaving segments can be broadly categorized into several technical approaches: graph search-based methods, sampling-based methods, dynamic programming (DP) combined with quadratic programming (QP), model predictive control (MPC)-based methods, interpolation curve-based methods, potential field-based methods, and learning-based methods. While these approaches have addressed certain challenges, they continue to face limitations that fail to capture the superimposed effects of multi-source risks and dynamic-static spatial-temporal coupling.

The deep integration of risk field with existing trajectory planning methods, alongside advancements in “vehicle-road-cloud integration (VRCI)”, offers new insights into overcoming these challenges. This integration facilitates multi-dimensional optimization and has been explored through several key approaches. Relevant studies are summarized in **Table 2**.

(1) **Direct trajectory generation:** This method generates low-risk discrete points or directional curves using risk field metrics as optimization objectives, often within frameworks like QP, MPC, or APF-based force guidance. While computationally efficient, it may produce trajectories that violate vehicle dynamics or smoothness constraints, requiring additional regularization.

(2) **Trajectory screening with risk field-based metrics:** A two-stage approach that first generates trajectory clusters via interpolation or curve fitting, then optimizes selection based on risk indicators while considering efficiency and comfort. This method allows for vehicle dynamics constraints but requires careful tuning of candidate trajectory numbers to balance computational cost and optimality.

(3) **Integration with reinforcement learning:** Risk field enhances deep reinforcement learning (DRL) by quantifying risk through physical properties. Used in reward, they improve vehicle interaction modeling and decision-making.

**Table 2** Risk field-based trajectory planning for CAV.

| Approach | Technical Features | Ref. |
|---|---|---|
| Using the risk field to directly generate trajectories. | The direction of the field force is initially employed as a motion guide to generate the trajectory, followed by secondary optimization to achieve trajectory smoothing. | (Liu et al., 2023; Xie et al., 2022) |
| | Risk-field related metrics are utilized as objective functions in numerical optimization, DP+QP, MPC, and other similar methods. | (Guo et al., 2025; Han et al., 2023; Hang et al., 2021; Xia et al., 2024; Yadollah et al., 2017) |
| Using risk field-based indicators as screening criteria. | A set of trajectories is first generated via sampling-based methods or fitting interpolation curves, followed by the selection of optimal trajectories using risk-field metrics, and finally, trajectory smoothing is performed. | (Fang et al., 2022; Tan et al., 2024; Wang et al., 2024c; Wu et al., 2023) |
| The risk field is combined with reinforcement learning. | Rewards are set up more scientifically to enhance interpretability. | (Li et al., 2023a; Wang et al., 2021; Wu et al., 2024) |

2.3 Research gaps and main contributions

Although existing research is comprehensive, the identified research gaps are as follows:

(1) Most risk field models are two-dimensional and static, considering only current risk sources while ignoring the future movement trends of dynamic obstacles and their impact on current risks, thus lacking risk foreseeability.

(2) The systematic calibration of risk field models remains insufficient due to the large number of involved parameters. The calibration approach relying on macro-statistical data, due to its inherent coarse-graining, tends to compromise model accuracy in specific scenarios. In contrast, the calibration approach based on micro-trajectory data exhibits limitations in its ability to generalize across different scenarios.

(3) Traditional trajectory planning methods struggle to account for risks from high-speed, coupled dynamic-static, and multi-source, often leading to imbalances in safety redundancy design or the neglect of dynamic risks.

(4) Balancing multi-dimensional optimization objectives while subject to strict real-time constraints and prioritizing safety remains challenging.

To fill these gaps, this study makes the following contributions:

(1) We innovatively introduce the concept of spatiotemporal distance, extending the conventional two-dimensional risk field into a three-dimensional STRF. Spatiotemporal distance effectively quantifies the impact of predicted trajectory and time on risk, addressing the gap between traditional driving risk assessment methods that fail to measure the impact of future movement trends of dynamic obstacles on current risk.

(2) We calibrate STRF parameters using real-world aerial data, leveraging machine vision techniques (YOLO) and dynamic risk balance theory. This provides a framework for high-precision model calibration work in specific interleaving zone scenarios.

(3) We propose a STRF-based CAV trajectory planning method in weaving segment. By deriving a spatiotemporal risk occupancy map (STROM) from the STRF, we constrain the sampling range within the dynamic time domain. Each sampling area corresponds to a unique STROM, ensuring that both dynamic and static risks are considered at every sampling point for enhanced safety. The proposed trajectory planning method addresses the shortcomings of traditional methods that are difficult to take into account risks from high speeds, dynamic and static coupling, and multiple sources.

(4) Our trajectory planning method also integrates STRF, dynamic iterative sampling, and QP to balance safety, efficiency, and comfort in a prioritized manner. Furthermore, to address the trade-off between multi-objective optimization and real-time performance, the dynamic iterative sampling and STROM we designed enable parallel path and speed optimization, and parallel computing can reduce the computation time. The integration of the new multi-objective optimization method with the parallel computing framework solves the challenging trade-off among safety, efficiency, comfort, and real-time performance.

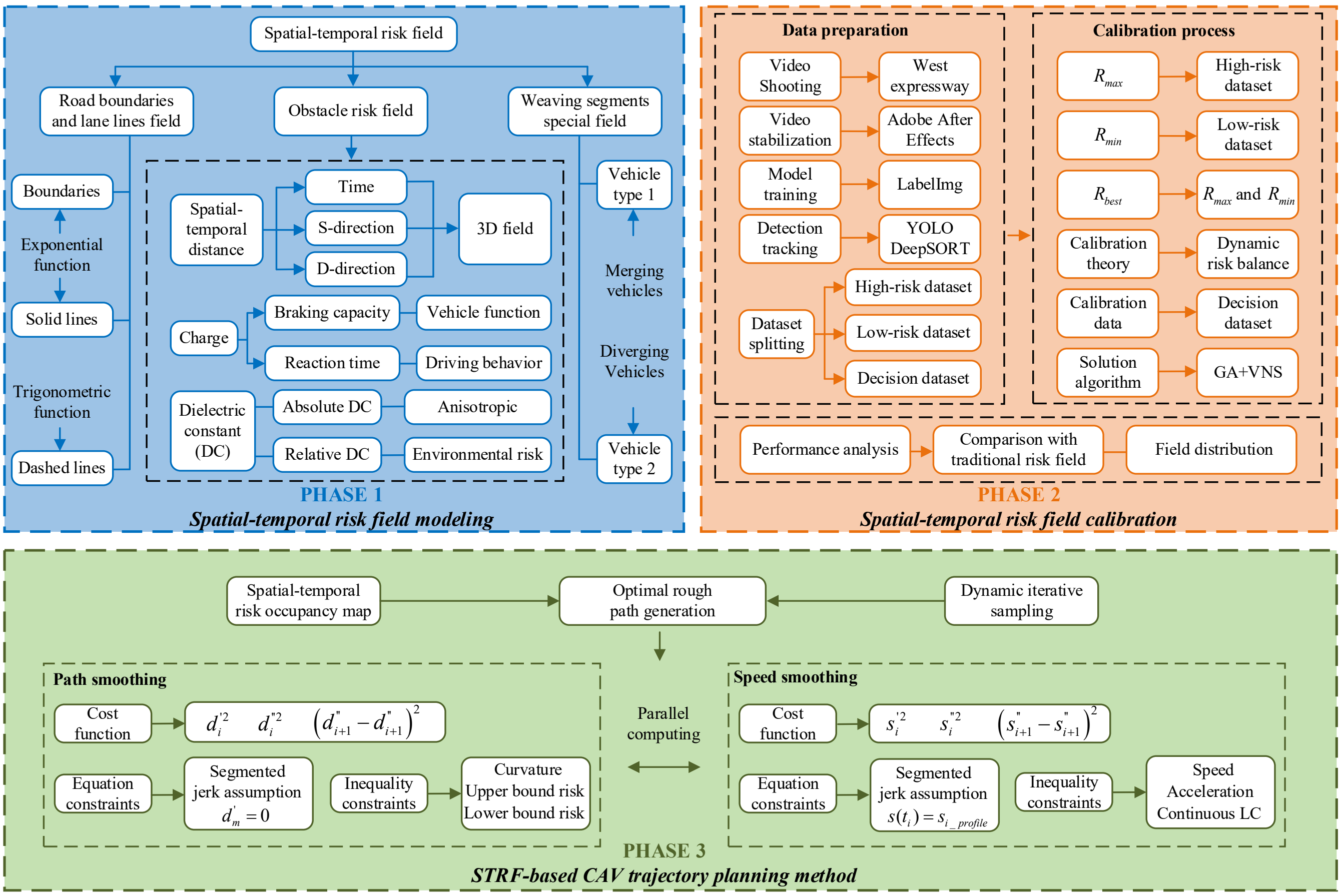

**Fig. 2** Overall framework diagram of this study.

The overall framework diagram of this study is shown in **Fig. 2**. The remainder of this paper is organized as follows: **Section 3** creates a 3D STRF model of the expressway weaving segment. **Section 4** calibrates the STRF's parameters using real-world aerial data based on the dynamic balance of risk theory and YOLO machine vision technology. **Section 5** designs a new STRF-based trajectory planning method in weaving segment. **Section 6** validates the safety, real-time performance, efficiency, and comfort of the proposed trajectory planning method using real-world cases. **Section 7** summarizes the study.

## 3 STRF: a new driving risk assessment paradigm

According to field theory, the entire traffic environment can be conceptualized as a risk field that influences vehicle movement according to specific rules. Unlike the basic segment, where vehicle interactions are relatively simple, the weaving segment involves multiple influencing factors. Consequently, the risk field in this segment is considered a superposition of fields generated by various traffic elements. Before analyzing the movement rules within this field, it is essential to first construct an appropriate function to generate the risk field. Directly formulating a risk field for the entire traffic environment is highly complex; however, this challenge can be addressed by decomposing the overall risk field into components corresponding to different traffic elements. These elements can be classified into three primary categories:

(1) Obstacles;
(2) Road boundaries and lane lines;
(3) Weaving segment-specific geometric characteristics.

The total risk field is obtained by summing the risk field of each element, as expressed in Eq. (1). To accurately model the risk field of three elements, we employ tailored construction methods, which are explained in the following subsections.

$$E_{total,A} = E_{obs,A} + E_{lane,A} + E_{geo,A} \tag{1}$$

where $E_{total,A}$ is the total field strength, $E_{obs,A}$ is the field strength imposed by all obstacles, $E_{lane,A}$ is the field strength imposed by road boundaries and lane lines, $E_{geo,A}$ is the field strength imposed by the on- and off-ramps in the weaving segment. The subscript $A$ indicates the observation position.

### 3.1 STRF triggered by obstacles

The form of the obstacle field can be analogized from the electrostatic field (Eq. (2)) in physics. In the strength formula of the electrostatic field, the charge of the particle, the distance from a point outside the particle to the particle itself, and the dielectric constant are the three main parameters. Analogously, we construct the overall form of the risk field as in Eq. (3), mapping the obstacle's intrinsic properties (mass, geometric characteristics) and state of motion (kinematic characteristics, speed, acceleration, yaw angle, etc.) to the particle's charge ($q$), the spatial-temporal distance between the observation point and the obstacle to the distance ($r$), and the characteristics of the roadway environment (traffic flow, density, velocity, roadway friction coefficient, ambient visibility, light intensity and its rate of change, road curvature, etc.) to the dielectric constant ($\varepsilon$). For ease of understanding, we provide a specific expression for the field strength as shown in Eq. (3). Next, we meticulously model each of the above three types of physical parameters based on the characterization of the obstacle risk field.

$$\boldsymbol{E}_{1\to 2} = \frac{1}{4\pi\varepsilon}\frac{q_1}{r^2}\frac{\boldsymbol{r}}{r} \tag{2}$$

$$
\begin{aligned}
E_{obs,A} &= \frac{1}{\varepsilon}\frac{q_{obs}}{r} \\
&= \vartheta \cdot m_{obs} \cdot \frac{exp\left[k \cdot cos(\psi)\cdot\left(\beta_1 \cdot v_{obs} + \beta_2 \cdot a_{obs}\right) + v_{obs}\cdot\gamma\cdot\frac{T_r}{a_{max}}\right]}{\sqrt{{T_{min}^{*}}^{2} + \alpha\left(t_{min} - t_1\right)^2}}
\end{aligned} \tag{3}
$$

where $\boldsymbol{E}_{1\to 2}$ is the charged particle 1 on the charged particle 2 generated field strength, $\varepsilon$ is the dielectric constant, $q_1$ is the charged particle 1 of the amount of charge, $\boldsymbol{r}$ is the vector distance from charged particle 1 to the charged particle 2, $r$ is the length of the $\boldsymbol{r}$. $\vartheta$ is the environmental complexity, $m_{obs}$ is the mass of obstacle, $v_{obs}$ and $a_{obs}$ are the velocity and acceleration of the obstacle, separately, $\boldsymbol{d}' = (s_A - s_{obs}, d_A - d_{obs})$ is the obstacle pointing to the observation point of the position vector, $\psi$ is the angle between the direction of the obstacle and $\boldsymbol{d}'$, $T_r$ is the vehicle reaction time, which the sum of the driver reaction time and the vehicle system lag, $a_{max}$ is the maximum deceleration of the vehicle, $k$, $\beta_1$, $\beta_2$, $\gamma$, $\alpha$ are the coefficient to be determined, $T_{min}^{*}$ and $t_{min}$ are described below.

*3.1.1 Assumptions, definitions, and characteristics*

To effectively model the STRF of obstacles, it is necessary to clarify their characteristics, which are summarized as follows:

(1) **Trigger mechanism:** Once an obstacle enters the perceptual range, the target vehicle is affected by the obstacle's field strength.
(2) **Anisotropy:** Unlike the isotropic nature of electrostatic field, the STRF distribution of the obstacle is non-uniform across different directions. Specifically, the smaller the angle between the observation point and the obstacle's driving direction, the greater the range and intensity of the obstacle's influence.
(3) **Velocity and acceleration bias:** The STRF amplifies the anisotropy described in (2) based on the obstacle's velocity and acceleration. For instance, if the obstacle has a forward velocity, its potential field distribution is biased forward, with acceleration having a similar effect.
(4) **Risk foreseeability:** In the intelligent and connected environment, CAVs predict obstacles' trajectories to detect anomalies and avoid risks, emphasizing the impact of future trajectories of obstacles on risk assessment. Unlike the traditional 2D risk field, this study integrates future trajectories into distance calculations, forming a 3D STRF. The 3D STRF demonstrates excellent risk foreseeability due to the integration of temporal and spatial factors.
(5) **Obstacle geometry:** Obstacles are physical entities with shape and mass. Representing obstacles as point masses ignores the risk impact of their geometry. Therefore, the influence of obstacles' geometric contours on risk must be considered.

*3.1.2 Spatial-temporal distance*

Most existing studies measure distance under the same temporal dimension, which, to some extent, leads to incorrect risk estimation. Taking **Fig. 3** as an example, it records the historical trajectory of the ego vehicle, and two future trajectories are selected for analysis through a multimodal trajectory predictor. At time $t$, when measuring the risk value of the ego vehicle to point A, the distance tendentiously used in existing studies is $d$. However, this may actually result in incorrect risk estimation. For instance, along predicted trajectory 1, the ego vehicle keeps moving toward point A in the subsequent period; thus, the actual risk should be higher than that calculated using $d$. Conversely, if the future trajectory shows the ego vehicle moving away from point A (e.g., predicted trajectory 2), the actual risk may be lower than the result calculated based on $d$. It can be seen that in the quantification of risk, the risk value is not only related to the spatial distance at the current moment but also to the future trajectory trend of the obstacle. Therefore, there exists a temporal effect that influences the risk value. Based on this perspective, a new distance characterization method is proposed, which we term spatial-temporal distance. This method not only considers the current spatial distance between the obstacle and the observation point but also takes the obstacle's future trajectory into account, featuring a spatial-temporal three-dimensional attribute.

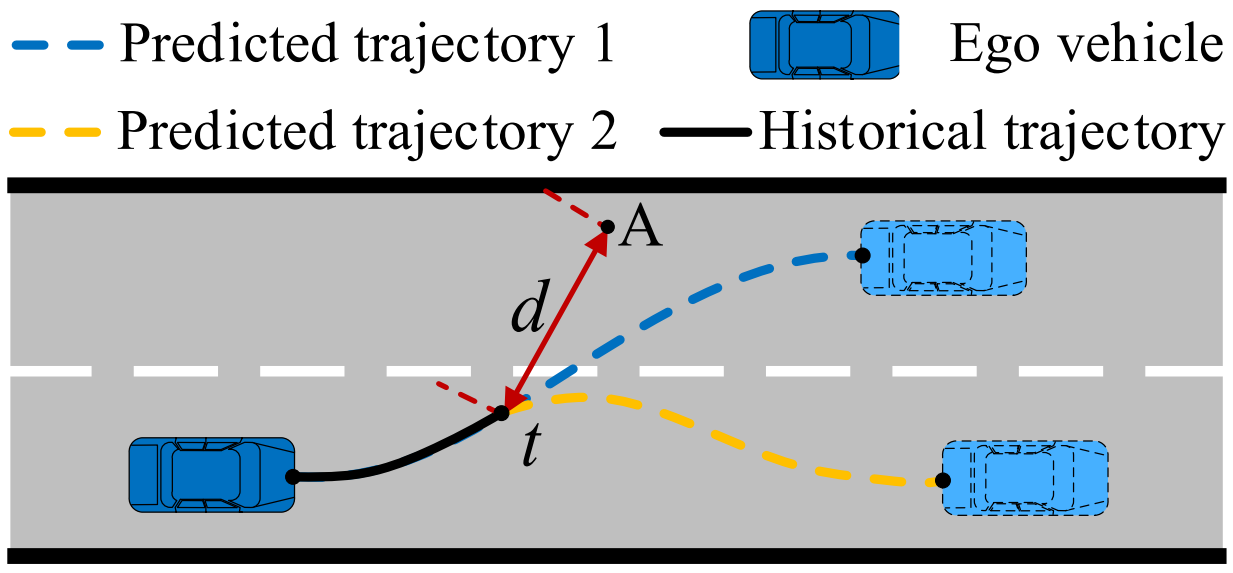

**Fig. 3** Case regarding spatial-temporal distance.

Based on the Frenet coordinate system, we establish a spatial-temporal, and three-dimensional coordinate system. The coordinate system comprise $S$-axes aligned with the road travel direction, $D$-axes perpendicular to the road direction, and the $T$-axis perpendicular to $S$-axes and $D$-axes. Using trajectory prediction, we plot the spatial-temporal trajectory of the obstacle. The spatial-temporal trajectories of the obstacle are projected onto a plane as spatial trajectories, with the observation point $A$ serving as the reference point and the $T_p$ as predictive time domain. The spatial-temporal trajectory of the obstacle in the range of $t_i$ to $t_j$ is $l_{i,j}$, and $t_m$ is the last position predicted by the trajectory. When observation time is $t_i$, the spatial-temporal distance between the obstacle and the observation point $A$ is the minimum spatial-temporal weighted distance from $A$ to the $l_{i,j}$. The visualization of the spatial-temporal distance is shown in **Fig. 4**, and the calculation formula is Eq. (4).

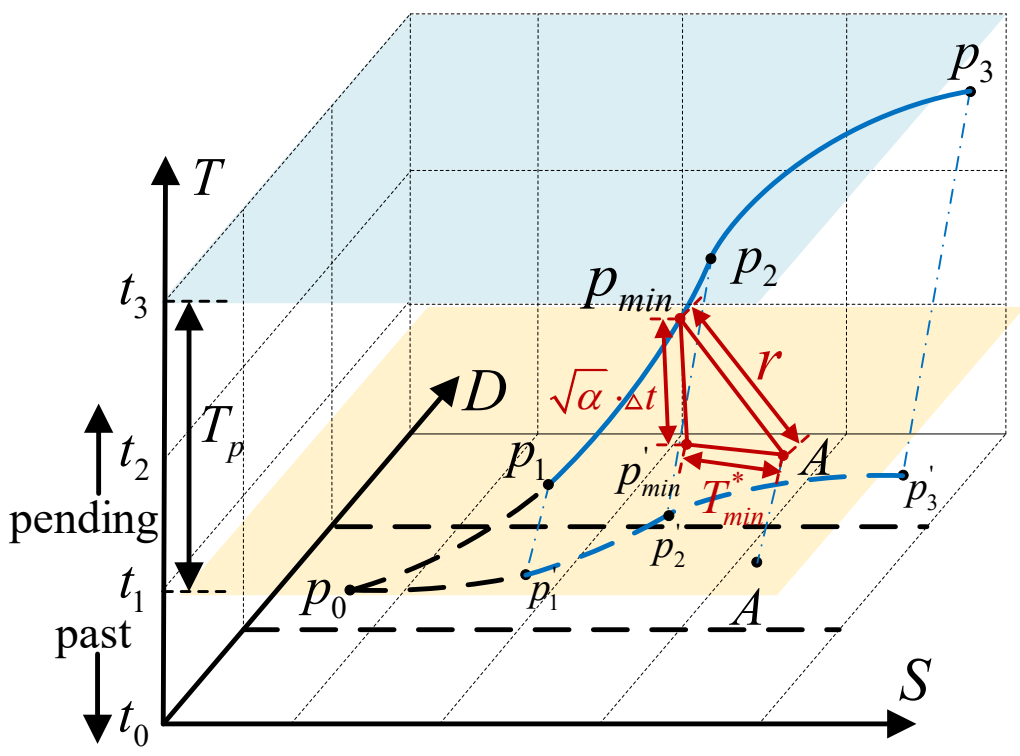

**Fig. 4** Spatial-temporal Frenet coordinate system and spatial-temporal distances.

$$r^i = \sqrt{{T^*_{min}}^2 + \alpha(t_{min} - t_i)^2} \tag{4}$$

where $r^i$ is the time distance between the observation point and the obstacle at the time $t_i$, $T^*_{min}$ is the equivalent time-based distance between the point (denoted as $p_{min}$) closest to observation point $A$ on the scattering point of the $l_{i,j}$ at the time $t_i$ and the observation point, the time required for the vehicle to arrive at the specified location is called the time-based distance, $t_{min}$ is the position of the $T$-axis of $p_{min}$ at the time $t_i$, and $\alpha$ is the weight of the time effect.

Since obstacle trajectories do not follow a specific function, we cannot calculate spatial-temporal distances directly. Calculating it requires a specific method. This method is divided into four steps:

**Step 1: Obtain the predicted trajectories of the obstacles.** In practical CAV system architectures, such future spatiotemporal state information is typically provided in real time by upstream perception and prediction modules. Because trajectory prediction is not the primary focus of this research, to ensure the input quality for the planning module and to validate the engineering feasibility of the proposed framework, this study adopts the advanced trajectory prediction algorithm introduced by Zhou et al. (Zhou et al., 2024) as the designated prediction module. This prediction model demonstrates exceptional accuracy and robust generalization capabilities within complex multivehicle interaction scenarios, and it can be seamlessly applied to urban expressway weaving segments following simple adaptations. The scatter set of obstacle trajectories is denoted as $P=\{p_i,p_{i+1},\cdots,p_m\}$, where $p_i$ represents the predicted spatial-temporal trajectory of the obstacle at the time $t_i$, $p_i=(s_i^{obs},d_i^{obs},t_i)$, $1\le i\le m$. $s_i^{obs},d_i^{obs},t_i$ are the position in the spatial-temporal coordinate system.

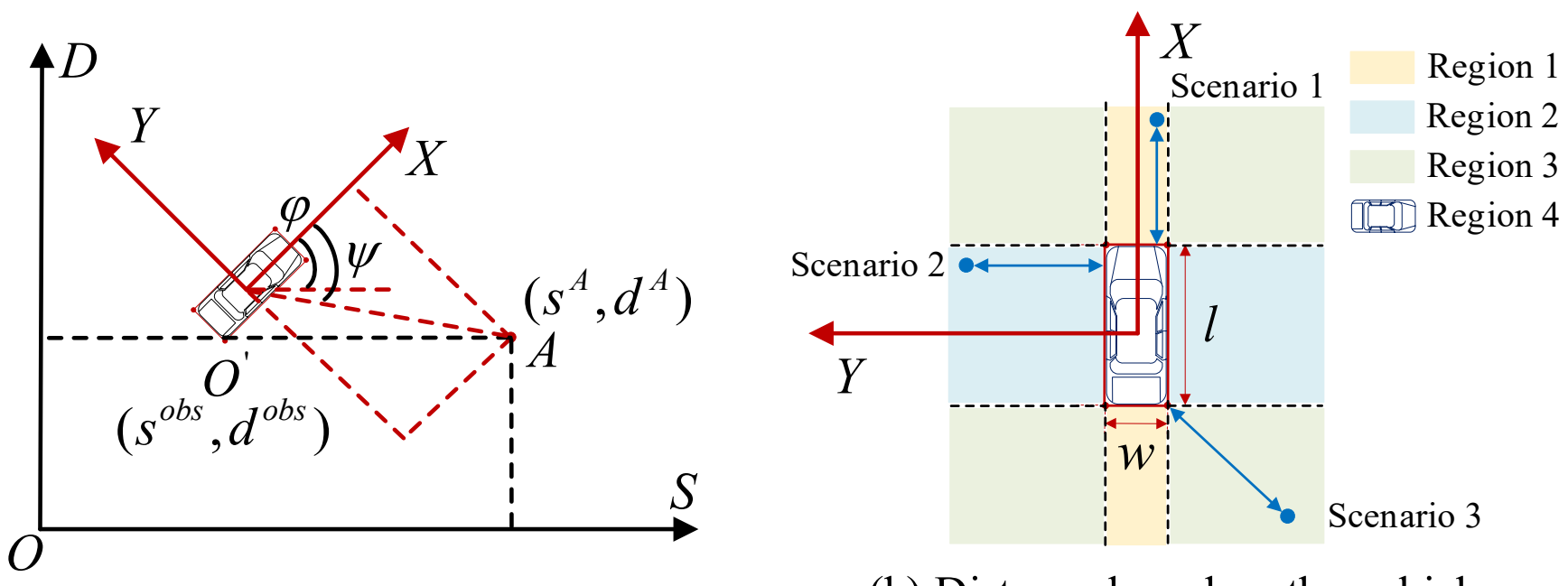

**Fig. 5** Schematic of coordinate transformation and distance calculation.

**Step 2: Calculate the spatial distance, considering the vehicle's geometric contour.** In this study, we consider the shape of the obstacle as a rectangle with definite values of length $l$ and width $w$, because most of the obstacles are interfering vehicles during driving.

$$\begin{bmatrix} x_A \\ y_A \end{bmatrix} = \begin{bmatrix} cos(\varphi) & sin(\varphi) \\ -sin(\varphi) & cos(\varphi) \end{bmatrix} \begin{bmatrix} s_A - s_{obs} \\ d_A - d_{obs} \end{bmatrix} \tag{5}$$

$$d = \begin{cases} |x_A| - \frac{l}{2}, & -\frac{w}{2} \le y_A \le \frac{w}{2} \wedge \frac{l}{2} \le |x_A| \\ |y_A| - \frac{w}{2}, & -\frac{l}{2} \le x_A \le \frac{l}{2} \wedge \frac{w}{2} \le |y_A| \\ min(d_1,d_2,d_3,d_4), & \frac{l}{2} \le |x_A| \wedge \frac{w}{2} \le |y_A| \\ 0, & -\frac{l}{2} \le x_A \le \frac{l}{2} \wedge -\frac{w}{2} \le y_A \le \frac{w}{2} \end{cases} \tag{6}$$

$$min(d_1,d_2,d_3,d_4) = min\left\{ \begin{aligned} &\sqrt{(x_A-\frac{l}{2})^2+(y_A-\frac{w}{2})^2}, \sqrt{(x_A-\frac{l}{2})^2+(y_A+\frac{w}{2})^2} \\ &,\sqrt{(x_A+\frac{l}{2})^2+(y_A-\frac{w}{2})^2}, \sqrt{(x_A+\frac{l}{2})^2+(y_A+\frac{w}{2})^2} \end{aligned} \right\} \tag{7}$$

where $(x_A, y_A)$ is the coordinate of the point $A$ in the vehicle's center coordinate system.

Since there is no change in the $T$-axes, we establish the spatial-temporal coordinate system and the vehicle's center coordinates system by ignoring the $T$-axes. As shown in **Fig. 5**, in the vehicle's center coordinates system, the center of mass of the vehicle is the origin, the direction along the vehicle is the $X$-axis, and the direction perpendicular to the vehicle is the $Y$-axis. Among them, under the spatial-temporal coordinate system, the coordinates of the point $A$ are $(s^A, d^A)$, and the coordinates of the center of mass of the obstacle are $(s^{obs}, d^{obs})$. In the vehicle's center coordinates system, the coordinates of the four vertices of the rectangle are $(l/2, w/2), (l/2, -w/2), (-l/2, w/2), (-l/2, -w/2)$. In addition, the yaw angle of the vehicle is $\varphi$, and the speed of the vehicle is $v^{obs}$. The coordinate transformation can be obtained by using the two steps of translation and rotation. The specific formula is shown in Eq. (5). We divided the spatial extent of the vehicle's center coordinates system into four regions, this is because each region has a unified formula to calculate the distance. The formulas for solving the spatial distances based on the body coordinate system are given as shown in Eqs. (6) and (7).

**Step 3: Calculate the equivalent time-based distance, accounting for the anisotropy of risk.** The Euclidean distance (spatial distances) in **Step 2** has two major drawbacks: **1)** it assumes that obstacles are equal in risk to objects in different locations but at equal distances, which is a clear violation of common sense. **2)** Euclidean distance is not the same physical property as time. There is no physical significance when weighting time with Euclidean distance. About these, we introduce the concept of equivalent time-based distance $T^*$. $T^*$ is a physical parameter that represents: when the interactive object $A$ is directly in front of the obstacle's driving path, and the interactive object $B$ is not directly in front of the vehicle's driving path, assuming that they have equal risk of collision, the temporal distance from the vehicle to the object $A$ is equal to the equivalent temporal distance from the vehicle to the object $B$. The problem of unequal risks in different directions can be handled by $T^*$. Furthermore, in terms of physical properties, it is time-based distance, which can be well weighted with the time dimension. We cite the formulas (Eqs. (8) and (9)) of Wang et al. (Wang et al., 2023) to transform Eqs. (6) and (7), the mathematical expressions can be seen in Eqs. (10) and (11).

$$T^* = \frac{\sqrt{y^2 + \mu^2 \cdot x^2}}{\mu \cdot v_{obs}} \tag{8}$$

$$\mu = 0.01476 + 0.8 v_{obs}^{-1} \tag{9}$$

$$T^* = \begin{cases} \dfrac{|x_A| - \frac{l}{2}}{v_{obs}}, & -\frac{w}{2} \le y_A \le \frac{w}{2} \wedge \frac{l}{2} \le |x_A| \\ \dfrac{\left(|y_A| - \frac{w}{2}\right)}{\mu \cdot v_{obs}}, & -\frac{l}{2} \le x_A \le \frac{l}{2} \wedge \frac{w}{2} \le |y_A| \\ min(T_1^*, T_2^*, T_3^*, T_4^*), & \frac{l}{2} \le |x_A| \wedge \frac{w}{2} \le |y_A| \\ 0, & -\frac{l}{2} \le x_A \le \frac{l}{2} \wedge -\frac{w}{2} \le y_A \le \frac{w}{2} \end{cases} \tag{10}$$

$$min(T_1^*, T_2^*, T_3^*, T_4^*) \tag{11}$$

$$= \frac{1}{\mu \cdot v_{obs}} \cdot min \left\{ \begin{array}{l} \sqrt{\mu^2 \cdot (x_A - \frac{l}{2})^2 + (y_A - \frac{w}{2})^2}, \sqrt{\mu^2 \cdot (x_A - \frac{l}{2})^2 + (y_A + \frac{w}{2})^2} \\ , \sqrt{\mu^2 \cdot (x_A + \frac{l}{2})^2 + (y_A - \frac{w}{2})^2}, \sqrt{\mu^2 \cdot (x_A + \frac{l}{2})^2 + (y_A + \frac{w}{2})^2} \end{array} \right\}$$

**Step 4: Combine the time-based distance and time to obtain the spatial-temporal distance.** The spatial-temporal distance $r_n^i$ at the time $t_i$ can be solved based on the formula $r_n^i = \sqrt{T_n^{*2} + \alpha(t_n - t_i)^2}$, exactly from the observation point to the scatter points of the vehicle trajectory. After the spatial-temporal distance of all scattering points is obtained, the minimum value of the spatial-temporal distance of each scattering point is selected as the spatial-temporal distance from $A$ to the obstacle. It is calculated as Eq. (12).

$$r^i = min(r_1^i, r_2^i, \cdots, r_n^i, \cdots, r_m^i) \tag{12}$$

where $r^i$ is the spatial-temporal distance between the obstacle and the observation point at the time $t_i$.

*3.1.3 Charge*

In the electrostatic field, charge is a fundamental physical quantity that characterizes an object's electrical properties (positive or negative) and its relationship to potential energy. Similarly, in the risk field, charge represents a combined property of vehicle function and driving behavior. We select two representative parameters from the combined attributes of vehicle function and driving behavior as the key influential factors of charge: braking performance and reaction time. These parameters are chosen because: (1) Braking performance reflects a vehicle's ability to avoid risk when faced with it; (2) Reaction time—which includes both driver response time and vehicle system latency—determines how effectively the driver and vehicle can perceive and respond to potential risks.

We measure the charge using a combination of the self-properties, kinematic properties, the vehicle mass, and acceleration. The formula for charge is defined as Eq. (13).

$$q_{obs} = m_{obs} \cdot exp(v_{obs} \cdot \gamma \cdot \frac{T_r}{a_{max}}) \tag{13}$$

*3.1.4 Dielectric constant*

In the electrostatic field, the dielectric constant is a physical quantity that describes a material's ability to enhance a capacitor's charge storage capacity. It is defined as the product of the absolute dielectric constant in a vacuum and the relative dielectric constant in a given environment. In the context of the risk field, we establish an analogy between the dielectric constant in the risk field and its counterpart in the electrostatic field. The similarities and characteristics of this analogy are analyzed as follows:

(1) **The absolute dielectric constant of STRF:** Unlike the isotropic dielectric constant in electrostatics, the risk field is inherently anisotropic. We introduce an anisotropic attenuation coefficient as its analogue, adopting a wedge-shaped potential field similar to prior research (Wang et al., 2023; Wolf and Burdick, 2008). We further refine the model by incorporating velocity and acceleration biases, enhancing its anisotropic characteristics.

(2) **The relative dielectric constant of STRF:** In the electrostatic field, the relative dielectric con-

stant adjusts the field strength based on environmental factors. Similarly, we define environmental complexity as an analogue, where roadway conditions influence driving risk. In safe environments (e.g., well-lit, dry pavement, free flow), the attenuation coefficient remains between 0 and 1, reducing risk. In hazardous conditions (e.g., dim lighting, wet pavement, congestion), it exceeds 1, amplifying risk.

We define the dielectric constant in the risk field as the product of the anisotropic attenuation coefficient and the ambient attenuation coefficient, as shown in Eq. (14).

$$\varepsilon = \frac{1}{\eta \cdot \vartheta} \tag{14}$$

where $\varepsilon$ is the dielectric constant, $\eta$ is the anisotropic attenuation coefficient, and $\vartheta$ is the environmental complexity. In this study, we let $\vartheta = 1$.

We define the anisotropic attenuation coefficient in Eq. (15), which accounts for variations in risk associated with different movement directions, speeds, and acceleration levels. Additionally, the environmental complexity represents the combined effect of all external environmental factors that impact vehicle safety.

$$\eta = \exp\left[k \cdot cos(\psi) \cdot \left(\beta_1 \cdot v_{obs} + \beta_2 \cdot a_{obs}\right)\right] \tag{15}$$

3.2 Risk field triggered by other elements

*3.2.1 Risk field triggered by road boundaries and lane lines*

Roadway boundaries and lane lines impose constraints on vehicle movement, influencing drivers to continuously adjust their lateral positioning to avoid penalties. As a result, the road boundary and lane line risk field vary based on lateral relative distance. The magnitude of this risk field is determined by the road boundary, the lane line configuration, and the vehicle's lateral position. The key characteristics of the risk field formed by road boundaries and lane lines are summarized as follows:

(1) **Triggering mechanism:** As long as a vehicle is on the road, it will continuously be affected by road boundaries and the road line field.
(2) **Nonlinear variation:** The risk field exerts minimal influence when a vehicle is centered within its lane. However, as the vehicle approaches a road boundary, the field strength increases nonlinearly, significantly amplifying its impact.
(3) **Risk-based classification:** Road boundaries and lane lines can be categorized into three distinct types based on their risk impact. **1)** Prohibited boundaries: Impassable road edges, such as guardrails and slopes, impose absolute restrictions on vehicle movement. **2)** Regulated boundaries: Areas that are physically passable but restricted by traffic regulations, such as solid lane lines that prohibit crossing. **3)** Permissible boundaries: Dashed lines that allow for lane changes, presenting the least restrictive conditions.

We selected two distinct functions to model the road boundary and road line risk field, based on the varying degree of change in field strength resulting from lateral position shifts. When a vehicle is near the centerline, small lateral position variations (between $-W/4$ and $W/4$) have minimal impact on the field strength, so we model it using a relatively stable trigonometric function. However, as the vehicle nears the road edge, even small lateral position changes induce significant variations in the field strength, we use the exponential function to capture these abrupt changes. The equations for the constructed road boundary and road line risk fields are presented in Eqs. (16) ~ (18).

$$E_{lane,\,A} = E_{lane1,\,A} + E_{lane2,\,A} + E_{lane3,\,A} \tag{16}$$

$$E_{lane\,i,\,A} = \begin{cases} \partial_i \cdot \left[ exp(\frac{W}{2} - | d_A - d_{lane\,i} |) - 1 \right], \; -\frac{W}{2} + d_{lane\,i} \le d_A \le \frac{W}{2} + d_{lane\,i}, i \in \{1, 2\} \\ 0, \, otherwise \end{cases} \tag{17}$$

$$E_{lane\,3,\,A} = \begin{cases} \partial_3 \cdot cos\left[ \frac{\pi}{W} \cdot \left( d_A - d_{lane\,3} \right) \right], \; -\frac{W}{2} + d_{lane\,3} \le d_A \le \frac{W}{2} + d_{lane\,3} \\ 0, \, otherwise \end{cases} \tag{18}$$

where $E_{lane,\,A}$ is the road boundary and the road line to the observation point $A$ of the field strength value, $E_{lane1,\,A}$ is the first type of road boundary to the observation point of the field strength, $E_{lane2,\,A}$ is second type of road solid line to the observation point of the field strength, $E_{lane3,\,A}$ is the third type of road dotted line to the observation point of the field strength, $W$ is the width of the road, $\partial_i$ and $\partial_3$ are coefficients to be determined, $d_{lane\,i}$ is the lateral position of the road line (road boundary) of type $i$.

*3.2.2 Risk field triggered by the specific geometry in weaving segments*

In the weaving segment, vehicles on the on-ramp, which need to merge into the mainline, must transition from the acceleration lane to the inside lane within specified spatial constraints. Similarly, vehicles that want to off-ramp must move from the mainline to the outside lane within the same spatial limitations. For a typical geometric configuration of the weaving segment, as illustrated in **Fig. 1**, the mandatory lane-changing behavior of vehicles on the on- and off-ramps is viewed as an adaptive decision influenced by the risk field associated with longitudinal acceleration and deceleration lanes. The characteristics of the risk field formed by the on- and off-ramps on the expressway are summarized as follows:

(1) **Trigger mechanism:** For vehicles on the on-ramp, the weaving segment risk field is generated when the following conditions are met simultaneously: 1) the vehicle intends to merge into the mainline, 2) the vehicle is in the acceleration lane, and 3) the vehicle is located within the mandatory lane change area. For vehicles on the mainline, the weaving segment off-ramp risk field is generated when the following conditions are met simultaneously: 1) the vehicle intends to exit onto the off-ramp, 2) the vehicle is in the mainline, and 3) the vehicle is within the mandatory lane-changing zone.

(2) **Non-linear characteristics:** The field strength does not change linearly with distance. It increases more rapidly as the lane-changing vehicle approaches the start of the off-ramp.

As an example, we calculate the field strength in the weaving segment as shown in Eq. (19). Similar forced lane changes occur in various traffic scenarios, including merging zones, work zones, and accident zones. While road and traffic conditions vary across these scenarios, the underlying principles remain consistent, leading to similar expressions.

$$E_{geo,A}=\begin{cases}\sigma_1\cdot\left[exp(x_{end}-x_A)^{\sigma_2}-exp(x_{end}-x_{start})^{\sigma_2}\right], p_A\in P_{mand}\wedge A\in\chi\vee A\in\delta\\ 0,\ A\in\gamma\end{cases} \tag{19}$$

where $E_{geo,A}$ is the weaving segment on and off the lane acceleration and deceleration lanes on the field strength of vehicle $A$, $s_{end}$ is the start of the off-ramp $S$-axis coordinates, $s_{start}$ is the beginning of the mandatory lane-changing zone $S$-axis coordinates, $p_A=(s_A,d_A)$ is the location of vehicle $A$, $P_{mand}$ is the mandatory lane-changing area range, $\partial_4$ is the parameter to be determined. $\chi$ is the set of vehicles on the on-ramp that intend to join the mainline, $\delta$ is the set of vehicles on the mainline that intend to go off-ramp, $\gamma$ is the set of remaining vehicles.

## 4 Parameter calibration of STRF with YOLO machine vision technology

The risk field model involves many undetermined parameters. Without effective parameter calibration, the performance and adaptability of the model will be greatly affected. There are three main existing calibration methods, which are introduced in detail as follows:

(1) **Macroscopic statistical data:** To enhance the model's applicability across diverse scenarios, macroscopic statistical data is frequently used as a calibration criterion (Ashutosh et al., 2023; Hua et al., 2022). This approach improves the model's generalizability by aligning its performance with aggregated real-world metrics. However, its coarse granularity often limits precision in specific scenarios, as it may obscure nuanced dynamic interactions.

(2) **Microscopic trajectory data:** To ensure accuracy in specific scenarios, microscopic trajectory data is employed as the calibration criterion (Li et al., 2020). This method captures detailed interactions with high fidelity in controlled environments. However, its effectiveness depends on data precision, which poses challenges when adapting to new environments. In such cases, recalibration with scenario-specific samples is often necessary, introducing inefficiencies and delays in model deployment.

(3) **SSMs:** To maintain consistency with traditional risk assessment methodologies, SSMs are used as the calibration criterion (Joo et al., 2023; Sun et al., 2023). This approach aligns the risk quantification outcomes of the model with well-established SSM frameworks, ensuring interpretability and comparability. Nonetheless, the calibrated model inherits both the strengths and limitations of the selected SSM, such as potential oversimplifications in representing complex risk dynamics.

Given we focus on on- and off-ramps in weaving segments, the trajectory data approach is the most suitable, maximizing accuracy while mitigating its limitations.

### 4.1 Data sources

To obtain high-precision micro-level vehicle motion states and ensure the broad representativeness of the research scenarios, this study employed UAV aerial photography to collect trajectory data in urban expressway weaving segments. Data collection was executed using DJI Mavic 2 Zoom UAVs. To ensure the complete capture of the full interaction process from the on-ramp to the off-ramp, the flight altitude of the UAVs was strictly maintained at the minimum visual height required to cover the entire weaving segment. Furthermore, to minimize the interference of high-altitude vibration and vehicle shadows on the machine vision recognition algorithms, all video data were recorded during the midday hours of clear, windless weekdays. The data collected during these periods encompass highly diverse traffic flow conditions, objectively and accurately reflecting microscopic vehicle interaction characteristics under varying traffic densities.

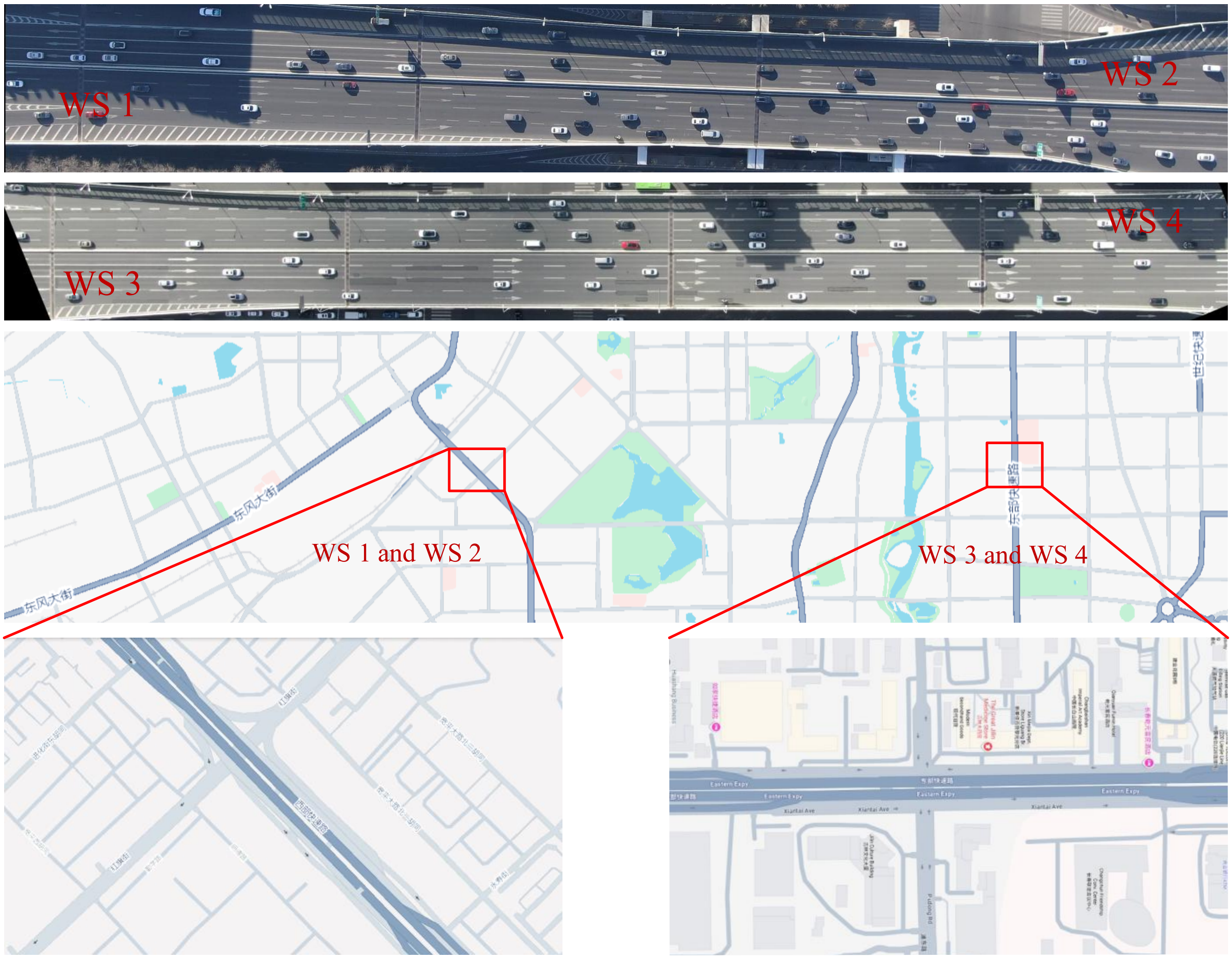

**Fig. 6** Study case.

**Table 3** Descriptive indicators of the studied weaving segments.

| | WS 1 | WS 2 | WS 3 | WS 4 |
|---|---|---|---|---|
| Time period | 10:00-12:00 | 10:00-12:00 | 12:00-14:00 | 12:00-14:00 |
| Video length (s) | 1547 | 1547 | 1437 | 1437 |
| WS length (m) | 290 | 290 | 375 | 375 |
| Resolution/frame | 2688×1512 60 FPS | 2688×1512 60 FPS | 3840×2160 30 FPS | 3840×2160 30 FPS |
| Location | WE | WE | EE | EE |
| Number of Lanes | 4 | 4 | 4 | 5 |
| Speed limit (km/h) | 80 | 80 | 80 | 80 |
| Number of total Vs | 3163 | 2873 | 2047 | 2715 |
| Usage | Calibration | Calibration | Validation | Validation |

Note: WE, EE represent Western and Eastern Expressway, respectively, WS represents weaving segment.

To guarantee the validity and applicability of the model parameter calibration, we constructed a dataset architecture based on a cross-validation mechanism. Specifically, two typical weaving segments were selected on both the Western Expressway and the Eastern Expressway, yielding empirical data from a total of four weaving segments. The data from Weaving Segment 1 and Weaving Segment 2 were exclusively designated for model parameter calibration, whereas the data from Weaving Segment 3 and Weaving Segment 4 were independently utilized for cross-validation. The geographical locations, geometric designs, and road channelization configurations of these four weaving segments are illustrated in **Fig. 6**. Detailed descriptive indicators for these weaving segments are provided in **Table 3**.

### 4.2 Extracting trajectories using YOLOv8

#### *4.2.1 Trajectory extraction procedure*

Based on YOLOv8, this study developed a specialized framework utilizing machine vision technology to extract microscopic vehicle trajectories from aerial videos. This framework is primarily composed of the following six core components:

**Step 1:** Video stabilization and cropping. When unmanned aerial vehicles capture footage from high altitudes, they are inevitably affected by wind and airflow, causing minor rotational and translational jitter in the video. This leads to deviations in the pixel coordinates of vehicles between frames. To address this, we utilized Adobe After Effects software to preprocess the video data. Researchers manually selected stationary landmark buildings in the footage as reference baselines to extract the overall rotational and translational trajectories. Subsequently, the system automatically applied reverse motion compensation technology to stabilize the video, thereby ensuring the accuracy of trajectory extraction. Furthermore, by manually cropping out irrelevant adjacent road sections, the processing efficiency of the subsequent automatic detection algorithms was significantly improved.

**Step 2:** Model training and testing. Most existing pretrained object detection models are trained on datasets acquired from conventional tilted monitoring perspectives, resulting in poor performance when processing vertical overhead data, as they struggle to adapt to the unique physical characteristics of a vertical view. However, a vertical overhead perspective provides more comprehensive and distinct vehicle contour information. Accordingly, this study utilized the LabelImg image annotation tool to manually construct an exclusive vertical overhead vehicle dataset. Subsequently, the YOLOv8 model was retrained using this dataset to ensure outstanding detection robustness and generalization capability across various flight altitudes.

**Step 3:** Object detection and multi-target tracking. Compared to two-stage detection algorithms like Faster RCNN, new single-stage detection algorithms like the YOLO series significantly enhance detection speed and effectively support real-time detection by directly predicting the position, size, and category of the target within a single network. The YOLOv8 algorithm employed in this study, which further optimizes detection accuracy and efficiency based on earlier versions, was responsible for precisely identifying vehicle targets in the video frame by frame. After completing object detection, the recognition results were deeply integrated with the DeepSORT multi-target tracking algorithm to achieve continuous tracking of vehicle trajectories. As an upgraded version of the SORT algorithm, DeepSORT incorporates a deep learning based appearance feature extraction module. Benefiting from the occlusion-free advantage provided by the vertical overhead perspective, this automatic tracking framework can accurately distinguish target vehicles based on unique appearance features, even when facing temporary vehicle occlusion or sudden changes in environmental lighting. This greatly reduces the expected ID switch rate and target loss rate, demonstrating exceptional algorithmic robustness in complex traffic scenarios.

**Step 4:** Lane line recognition and assignment. To accurately determine the physical lane positions of the vehicles, we extracted representative stable video frames from the processed footage and manu-ally calibrated the physical boundaries of all lane lines using the LabelMe tool. The resulting lane coordinate data matrix serves as the core reference baseline for the subsequent fully automatic determination of vehicle lane changing behaviors and lane assignment results.

**Step 5**: Coordinate system transformation. Constrained by the battery capacity of the recording equipment, slight variations in the actual flight altitude and hovering position of each video group are inevitable, leading to spatial positional deviations between different videos captured at the same location. Furthermore, original pixel coordinates cannot be directly used to accurately describe microscopic vehicle kinematics. To overcome this technical challenge, this study referred to the methodology of Mardiati et al. (Mardiati et al., 2019) to obtain authentic vehicle trajectory coordinates within a predefined reference coordinate system. Building upon this, the research team extracted the spatial position information of road boundaries and physical lane lines. Using this as a transformation reference, all two-dimensional coordinate information was further mapped into the Frenet coordinate system.

**Step 6:** Trajectory smoothing. Inevitable minor positional errors from pixel mapping exist in the raw positioning data. These errors are significantly amplified when calculating the derivatives for speed and acceleration, causing severe short-term fluctuations in the trajectory data. In light of this, the workflow introduced the Kalman filtering algorithm to smooth all extracted raw trajectories. After rigorous filtering optimization, sudden speed variations between adjacent video frames were effectively suppressed, and the overall position error was strictly confined to a reasonable and minimal range.

*4.2.2 Verification of trajectory extraction accuracy*

To ensure the data reliability for the subsequent spatial temporal risk assessment and trajectory planning models, this section conducts an exhaustive quantitative verification and impact analysis of the expected errors in the trajectory extraction work.

At the object detection level, the trained YOLOv8 model demonstrates favorable recognition precision and convergence stability. As illustrated by the model training process curves in **Fig. 7**, throughout the entire training cycle, the bounding box loss, classification loss, and distribution focal loss of the model exhibit a smooth downward trend on both the training and validation sets, ultimately converging to low levels. This indicates that the model learns the spatial and categorical features of the targets effectively, avoiding the issue of overfitting. Simultaneously, various core performance evaluation metrics experience a rapid climb in the early training stages and then maintain a stable plateau. Specifically, as shown in **Table 4**, the finally converged model achieves a precision of 0.9980 and a recall of 0.9948 on the validation set. The mAP50 metric reaches 0.9948, and the more stringent mAP50 to 95 metric also reaches 0.8205. Such stable training curves and evaluation metrics objectively demonstrate that the model can accurately bound the vehicle targets within the weaving segment (relevant examples can be seen in **Fig. 8**).

At the multiple object tracking level, detection is merely the first step. The core task of the tracking algorithm is to associate isolated detection boxes into continuous trajectories over the time series. In weaving segment scenarios, vehicles frequently execute transverse cut in and cut out maneuvers, resulting in complex trajectory intersections. This poses a significant challenge to the identity maintenance capability of the tracker. Therefore, evaluating tracking completeness and ID stability is necessary. MOTA and IDF1 are two important metrics for evaluating ID stability. MOTA integrates three types of errors (false negatives, false positives, and ID switches) and serves as the core metric for measuring comprehensive tracking performance. As shown in **Table 4**, IDSW is a relatively sensitive error type in weaving segment research. The MOTA of 0.9985 in this study verifies reliable tracking performance. IDF1 focuses on the time proportion during which the tracker identifies the same target, which better reflects the completeness of the trajectory. Herein, IDTP represents the count of bounding boxes correctly associated with the true ID, IDFP represents the count of bounding boxes incorrectly associated with an ID, and IDFN is the count of boxes failing to associate with an ID. In weaving segment scenarios characterized by frequent lane changing behaviors and large vehicle density, the IDF1 of 0.9992 indicates that the algorithm can stably maintain vehicle identities and ensure the completeness of vehicle trajectories.

In terms of trajectory position error and its impact assessment, although the model has high precision, certain errors are still inevitable in the initially extracted trajectory. However, the impact on the final risk assessment and trajectory planning is kept within a fully acceptable range. Because a vertical overhead perspective is adopted, trajectory breakpoints caused by physical occlusion are effectively avoided. To address the inherent pixel level mapping deviations during the coordinate system transformation, as well as the short term fluctuation noise amplified when calculating velocity and acceleration through position differentiation, this study introduces the Kalman filtering algorithm for post processing. The filtered trajectories smoothly eliminate unrealistic velocity mutations between frames, effectively reducing absolute position errors. Furthermore, the dynamic spatial temporal risk field model constructed in this paper inherently possesses robustness against slight perturbations in microscopic trajectory data at its foundational architecture level. Consequently, these residual minor errors, after appropriate control and algorithmic smoothing, do not exert a substantial impact on the precision of the final risk identification.

**Table 4** Indicators on the test set

| Precision | Recall | mAP50 | mAP50-95 | MOTA | IDF1 |
|---|---|---|---|---|---|
| 0.9980 | 0.9948 | 0.9948 | 0.8205 | 0.9985 | 0.9992 |

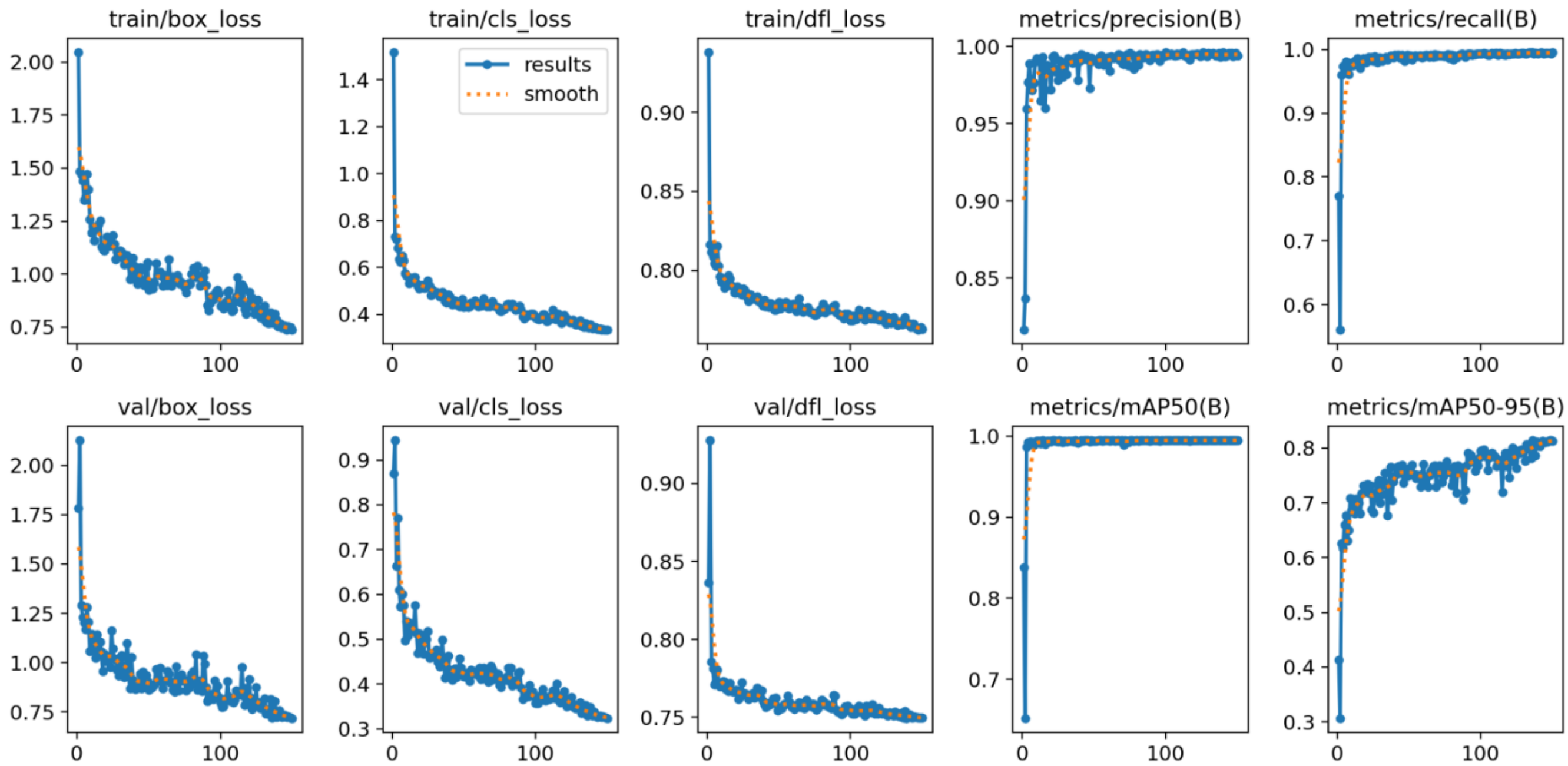

**Fig. 7** Convergence curves of the YOLOv8 model.

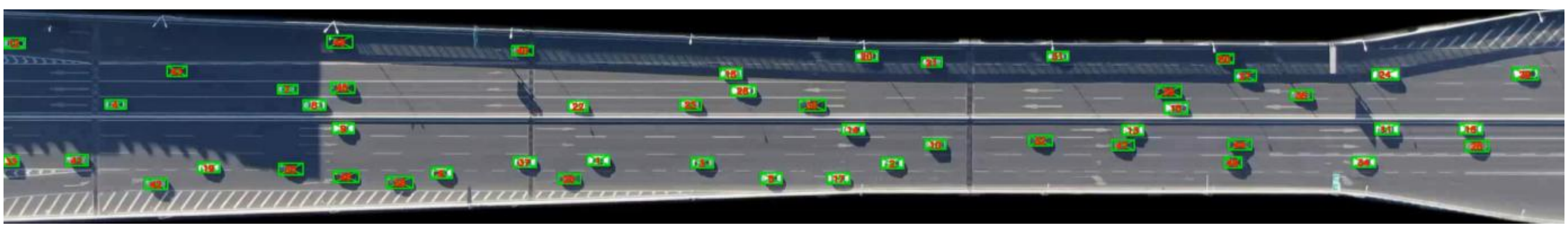

**Fig. 8** Vehicle identification and tracking effect from aerial photography.

### 4.3 Parameter calibration based on the dynamic risk balance theory

Inspired by Tan et al. (Tan et al., 2022), we observe that risk levels exhibit regular patterns before and after vehicle decision-making. Specifically, the risk value experienced by the vehicle after a decision fluctuates within a certain threshold. This phenomenon aligns with dynamic risk balance theory, which posits that drivers regulate risk within a subjectively acceptable range. This regulation results from a comprehensive assessment of expected behavioral benefits and potential health and safety impacts. During driving, drivers continuously perceive and evaluate risk levels, comparing them with their expected threshold and adjusting their behavior accordingly. If the perceived risk falls below the acceptable level, they tend to take greater risks; conversely, if it exceeds the acceptable level, they adopt a more cautious approach. Consequently, driver decision-making serves as an adaptive mechanism that aligns with the dynamic risk balance theory, ensuring that the vehicle and driver maintain risk exposure within predefined upper and lower limits. Based on this theory, we calibrate the risk field model using trajectory data.

#### *4.3.1 Calibration procedure*

Based on dynamic risk balance theory, we take the mean squared deviation of the distance between the sample risk and the acceptable risk space (i.e., the space comprising the upper and lower limits of risk) after decision-making as the objective function of the calibration model, which is expressed in Eq. (20) and Eq. (21). Calibration steps and processes can be seen in **Fig. 9**.

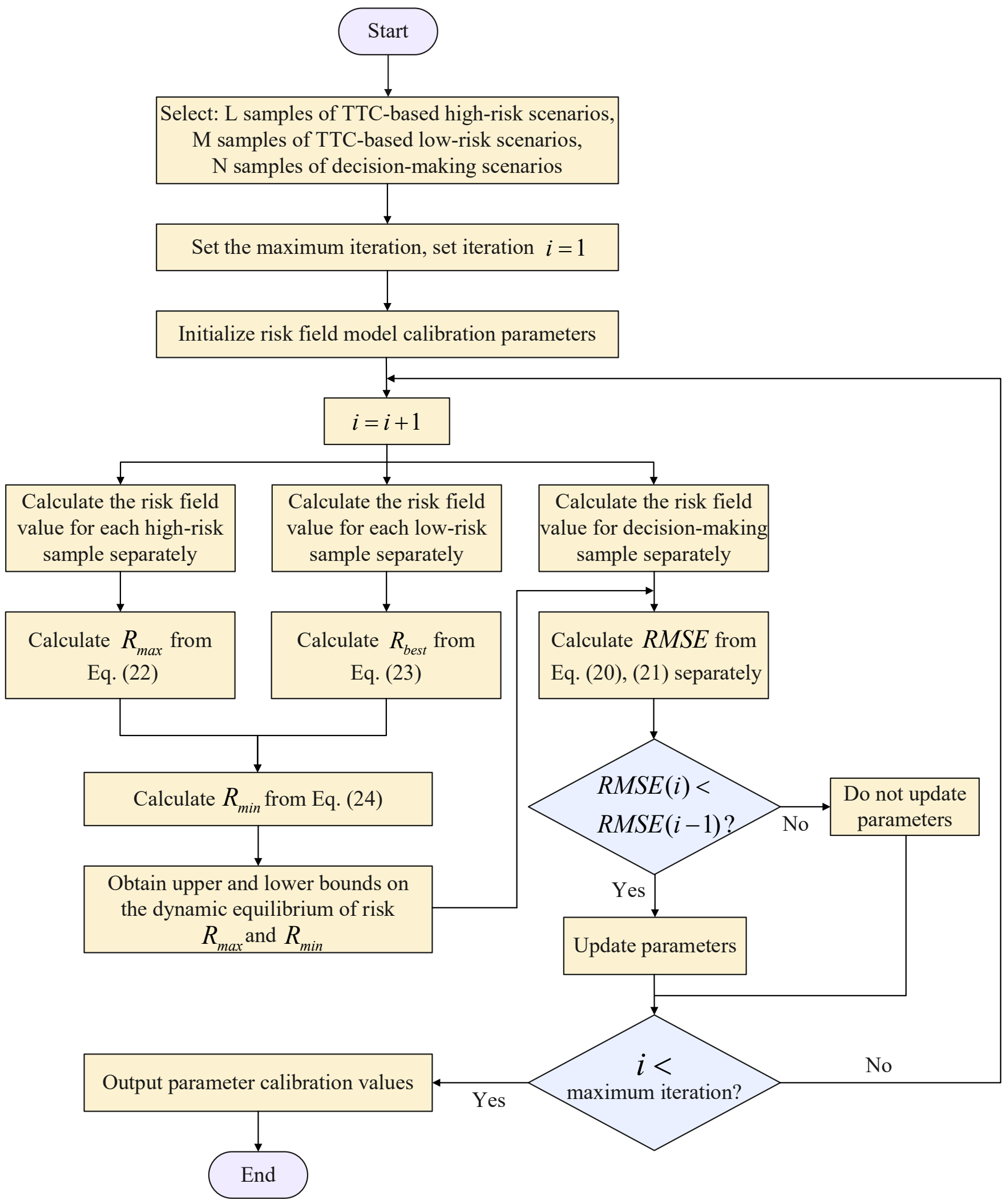

**Fig. 9** STRF calibration algorithm flow.

$$RMSE = \sqrt{\frac{\sum_{n=1}^{N} |\triangle R_n|^2}{N}} \tag{20}$$

$$\triangle R_n = \begin{cases} R_n - R_{max}, R_n > R_{max} \\ R_{min} - R_n, R_n < R_{min} \\ 0, R_{min} \le R_n \le R_{max} \end{cases} \tag{21}$$

where $RMSE$ is the objective function of the calibration model, $\triangle R_n$ distance between the post-decision risk of the $n^{th}$ sample and the acceptable risk space, $R_n$ is the value of risk at a specific moment after the decision of the $n^{th}$ sample, $R_{max}$ and $R_{min}$ are the upper limit risk and lower limit risk of the acceptable risk space, respectively, and $N$ is the number of samples of decision-making scenarios.

*4.3.2 Determination of $R_{max}$ and $R_{min}$*

We calculate and using the three types of datasets in **Fig. 9**: high-risk, low-risk, and decision datasets. Each dataset contains trajectories of an observation vehicle (ego vehicle) along with multiple surrounding vehicles. The high-risk dataset includes scenarios in which the observation vehicle is likely to encounter elevated risk, as assessed by TTC. We use a TTC threshold of 3 seconds. This dataset helps determine the upper limit, denoted as $R_{max}$. Conversely, the low-risk dataset includes scenarios where the observation vehicle remains consistently safe, as determined by the same TTC threshold. We use this dataset to establish the optimal risk level, denoted as $R_{best}$. Based on $R_{max}$ and $R_{best}$, we determine the $R_{min}$. The decision-making dataset includes scenarios where the observation vehicle engages in decision-making, such as lane-changing, sharp deceleration, or sharp acceleration. We use this dataset to determine $R_m$, which represents the level of risk the observation vehicle experiences after making a decision. By integrating three datasets, we fully define the objective function for the calibration model.

We randomly selected 30 sets of high-risk datasets, 100 sets of low-risk datasets, and 200 sets of decision-making datasets. Consequently, we can calculate the strength values of the surrounding vehicles to the observation vehicle under a specific parameter (which is a calibrated quantity). As a result, we obtained 30 strength values from the 30 sets of high-risk datasets, and took the smallest field strength value (according to Eq. (22)) as $R_{max}$. We also obtained 100 field strength values from the 100 sets of low-risk datasets, and took the average field strength value of them (according to Eq. (23)) as $R_{best}$. According to Eq. (24), the value of $R_{best}$ was determined.

$$R_{max} = min\{R_1, R_2 \cdots R_L\} \tag{22}$$

$$R_{best} = \frac{1}{M} \cdot \sum_{m=1}^{M} R_m \tag{23}$$

$$R_{min} = 2R_{best} - R_{\max} \tag{24}$$

*4.3.3 Calibration results*

To achieve an efficient solution and precise calibration of the complex multidimensional parameters in the model, this study introduces a hybrid variable neighborhood search and genetic algorithm (Ma et al., 2025b). Given that the parameter space of the model constructed in this paper possesses highly nonlinear and coupled characteristics, conventional single optimization algorithms often face the risk of search stagnation or premature convergence. This hybrid optimization framework deeply integrates the global population evolution advantages of the genetic algorithm with the local refined exploitation characteristics of the variable neighborhood search. While ensuring population diversity, it greatly enhances the ability of the algorithm to escape local optimal traps, thereby ensuring the stable acquisition of the global optimal solution.

In the specific implementation process of parameter optimization, we set the global maximum number of iteration cycles to 100 to ensure that the algorithm has ample evolutionary space for optimization. The fitness convergence curve for the entire parameter calibration process is detailed in **Fig. 10**. From this figure, the rapid downward trend of the algorithm in the early stage of iteration and the stable convergence state in the later stage can be clearly observed, verifying the superior performance of the hybrid algorithm. After rigorous optimization calculations, the final calibrated optimal values of all core parameters have been comprehensively summarized in **Table 5**. Furthermore, to ensure the scientific rigor and computational efficiency of the optimization process, we determined these initial values by referring to prior literature (Li et al., 2020; Wang et al., 2023) and theoretical inferences, which aligns with the physical boundary constraints and traffic flow dynamics specific to the urban expressway weaving segment scenarios in this study.

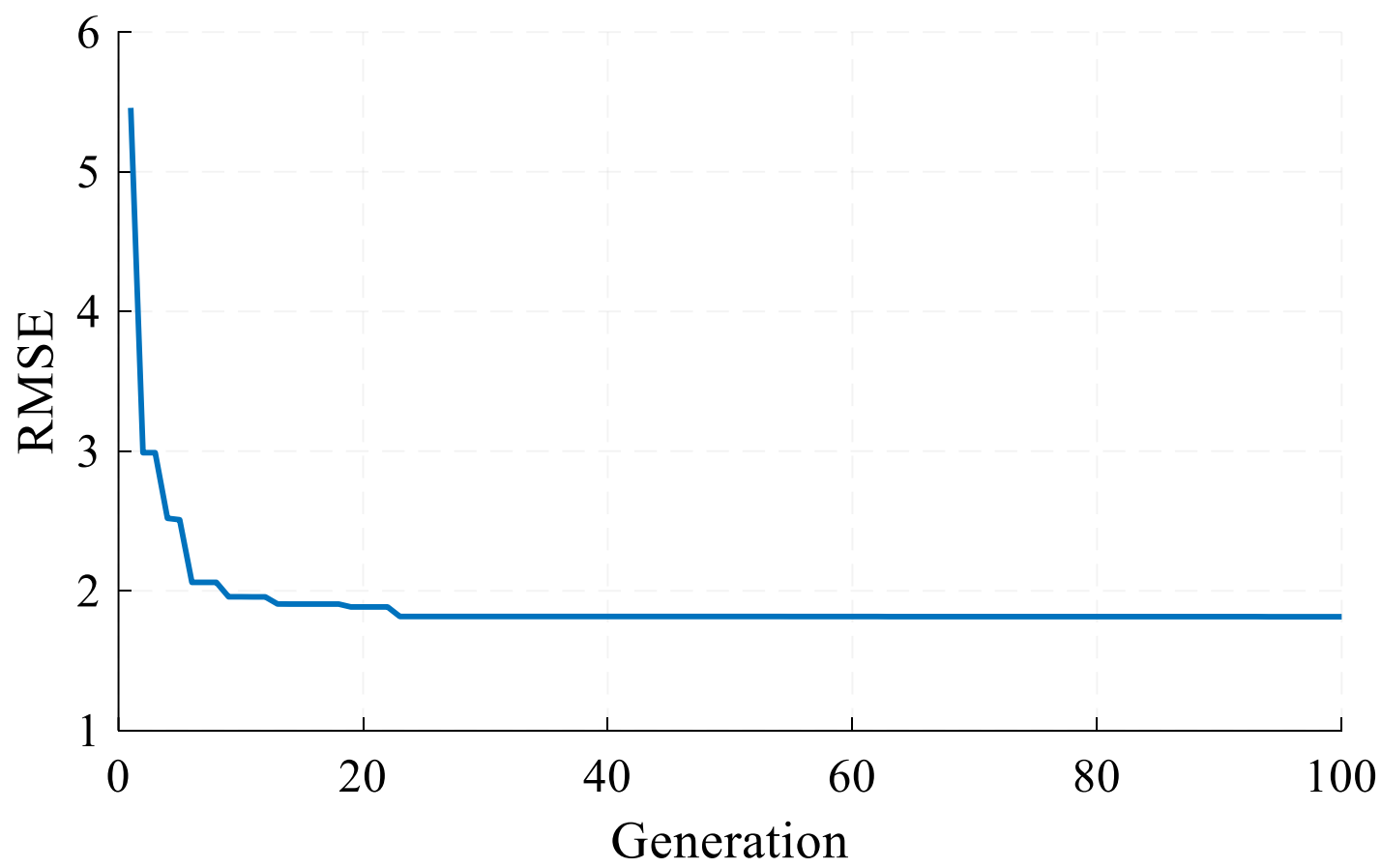

**Fig. 10** Calibration convergence curve of STRF.

**Table 5** Parameter calibration result.

| Parameters | $\alpha$ | $\beta_1$ | $\beta_2$ | $k$ | $\gamma$ | $\partial_1$ | $\partial_2$ | $\partial_3$ | $\sigma_1$ | $\sigma_2$ | $R_{max}$ | $R_{min}$ |
|---|---|---|---|---|---|---|---|---|---|---|---|---|
| Value | 1.72 | 0.07 | 0.25 | 0.56 | 0.09 | 2.02 | 1.06 | 2.05 | 9.55 | -0.45 | 4 | 1.2 |

## 4.4 Verification of applicability and stability of STRF

### *4.4.1 Applicability analysis of STRF*

To further verify whether the proposed model possesses strong generalization and applicability in other weaving segment scenarios, this section conducts an in-depth cross-validation analysis. In terms of experimental design, we precisely partitioned the data from WS 3 and WS 4 to construct four validation datasets (denoted as V1 to V4) that share the exact same specifications as the decision dataset (denoted as T1) utilized in the previous calibration phase. By directly applying the calibrated model to these four validation datasets, we calculated their respective root mean square errors according to Eq. (20), and conducted a comprehensive comparison with the performance indicators of the original decision dataset.

As illustrated in **Fig. 11**, when the model is transferred to the four validation datasets, its RMSE exhibits only an extremely minor fluctuation range. The various validation errors align highly with the baseline error, avoiding any overfitting phenomenon characterized by a sudden surge in error metrics. The comparison results powerfully demonstrate that the constructed model can still maintain high accuracy when processing data from other weaving segments with differentiated characteristics. This cross-scenario generalization capability lays a solid theoretical and empirical foundation for the future large-scale engineering deployment of the model in complex road networks.

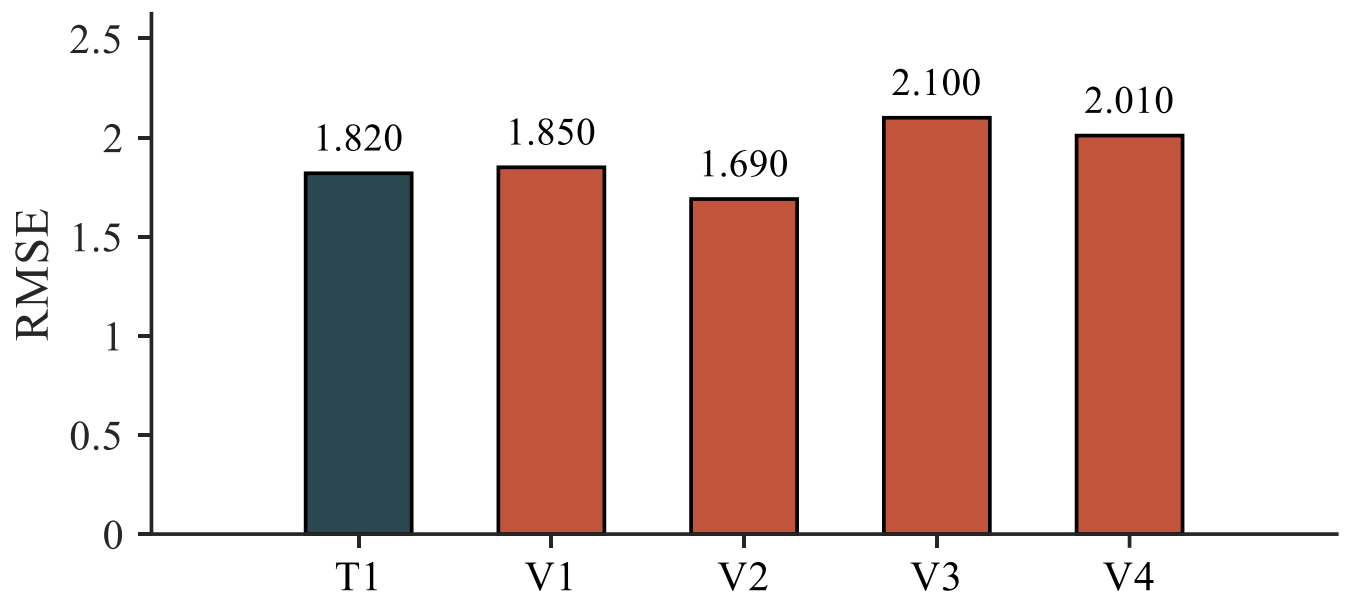

**Fig. 11** Comparison of RMSE across different data samples.

*4.4.2 Analysis of parameter sensitivity of STRF*

Given that the STRF model incorporates multiple calibrated parameters, evaluating their sensitivity in risk identification is crucial for validating the overall robustness of the model. To this end, this section selects three core parameters ($\beta_1$, $\beta_2$, and $\gamma$) that significantly dictate the shape and intensity distribution of the risk field to conduct an in-depth sensitivity analysis. In the experiment, we randomly selected nearly 8,000 vehicle groups from the two scenarios WS 3 and WS 4 as the research objects. To closely align with the core theme of this research, each group consists of one ego vehicle executing either an on-ramp or off-ramp maneuver, alongside multiple surrounding vehicles. Each data group was uniformly cropped to a 3-second duration. Within each group, we calculated the dynamic risk field strength imposed on the ego vehicle by surrounding vehicles and the traffic environment at every moment. If the field strength endured by the ego vehicle exceeded the maximum allowable risk threshold $R_{max}$ at any instant within the 3-second window, the sample was flagged as a high-risk interaction event.

By applying a 10% upward and downward perturbation to the calibrated values of these three parameters, we mapped the variations in the number of high-risk events and their proportions within the total sample across three pairwise parameter combinations. As illustrated in **Fig. 12**, when all parameters are set to their baseline calibrated values, the model accurately identifies 1186 high-risk events, accounting for 14.90% of the total samples. The color gradients in the heatmaps intuitively reveal the varying sensitivities of the model outputs to different parameters. First, the model exhibits exceptional stability regarding $\beta_2$ and $\gamma$. In the perturbation combinations involving $\beta_2$ and $\gamma$ (as shown in **Fig. 12** (c)), whether $\beta_2$ and $\gamma$ increases or decreases, the variations in the number and proportion of high-risk events remain extremely minimal. This indicates that even if there is a slight estimation deviation when evaluating the risk contribution of a vehicle's instantaneous acceleration, it will not compromise the overall risk determination. Second, although the model is more sensitive to $\beta_1$ compared to $\beta_2$ and $\gamma$, its overall output maintains excellent stability. This phenomenon perfectly aligns with physical expectations because $\beta_1$ determines the severity of the velocity bias of the risk field. In high-speed weaving segments, velocity is inherently the most critical kinetic factor dominating collision risks. Therefore, it is entirely normal and physically sound for the model to exhibit moderate sensitivity to velocity-related parameters.

In conclusion, these parameter sensitivity analyses clearly demonstrate that minor fluctuations in the core parameters have a limited impact on the overall accuracy of the model, and all variation trends remain within a reasonable and physically expected range. This thoroughly proves the outstanding robustness and engineering reliability of the proposed STRF model in parameter calibration.

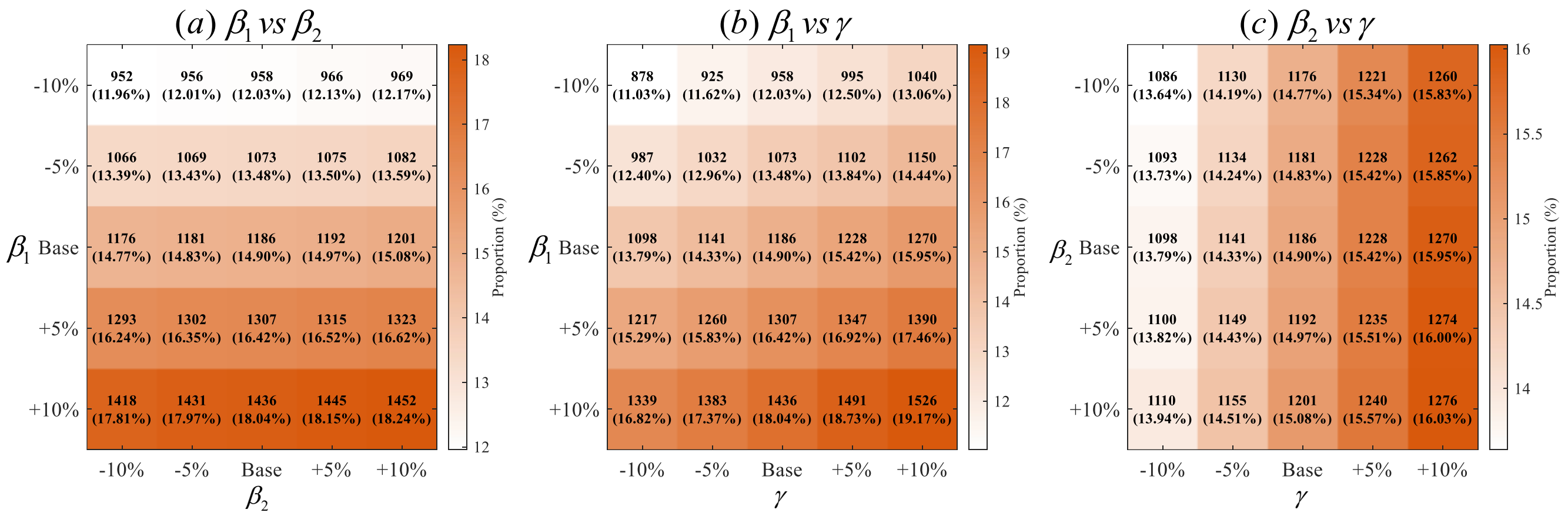

**Fig. 12** Impact of STRF parameter fluctuations on high-risk group identification.

## 4.5 STRF performance analysis

### *4.5.1 Risk predictability capability analysis*

To address the limitations of the traditional risk field, we introduce several key improvements in the STRF. These include: 1) incorporating the concept of spatial-temporal distance to construct a three-dimensional field, 2) accounting for the geometry of obstacles, 3) considering the coupling effects of kinematics, and 4) modeling the unique geometric characteristics of on- and off-ramps in weaving segments. However, among these improvements, the impact of spatiotemporal distance on risk assessment is the most significant. The introduction of spatiotemporal distance elevates the risk field from the traditional two-dimensional spatial field to a three-dimensional spatiotemporal field. It quantifies the influence of the future trajectory trends of dynamic obstacles on the current risk assessment, demonstrating the predictive ability of risk. To visualize its impact, we compared the STRF with the traditional two-dimensional risk field. To visualize the impact of spatiotemporal distance, the traditional risk field was modified; except for the spatiotemporal distance, all other model constructions are consistent with the STRF.

**Fig. 13** illustrates the schematic field strength of a traffic vehicle during an overtaking maneuver. **Fig. 13** (a), **Fig. 13** (c), **Fig. 13** (e), and **Fig. 13** (g) depict the field strengths formed by the STRF, while the remaining subfigures represent those formed by the traditional risk field. The following patterns can be observed:

(1) In the traditional risk field, the field strength exhibits a symmetrical shape along the vehicle's front-end direction. This is because the traditional field evaluates risk solely based on factors present at the moment of observation, without anticipatory perception based on predicted trajectories. In contrast, the symmetry of the field strength in the STRF depends on trajectory predictions.

(2) The field strength of the STRF is asymmetric in most cases, as seen in **Fig. 13** (a), **Fig. 13** (c), and **Fig. 13** (e). Taking **Fig. 13** (a) as an example, at T = 1s, the vehicle is moving longitudinally along the lane, but its trajectory shifts over the following period, resulting in a rightward lane-changing. Due to the trajectory prediction module's anticipatory capability, this information is incorporated into the spatial-temporal distance computation, causing the field strength to shift toward the predicted trajectory. Although the vehicle eventually returns to its original lane, the impact of returning to the original lane is significantly smaller than that of the rightward lane change due to the longer time interval. The combined effect of future trajectories results in the observed field strength distribution.

(3) In a few cases, the STRF remains symmetric, as shown in **Fig. 13** (g). At T = 7s, none of the predicted trajectories indicate a directional change, and consequently, the field strength remains symmetrical.

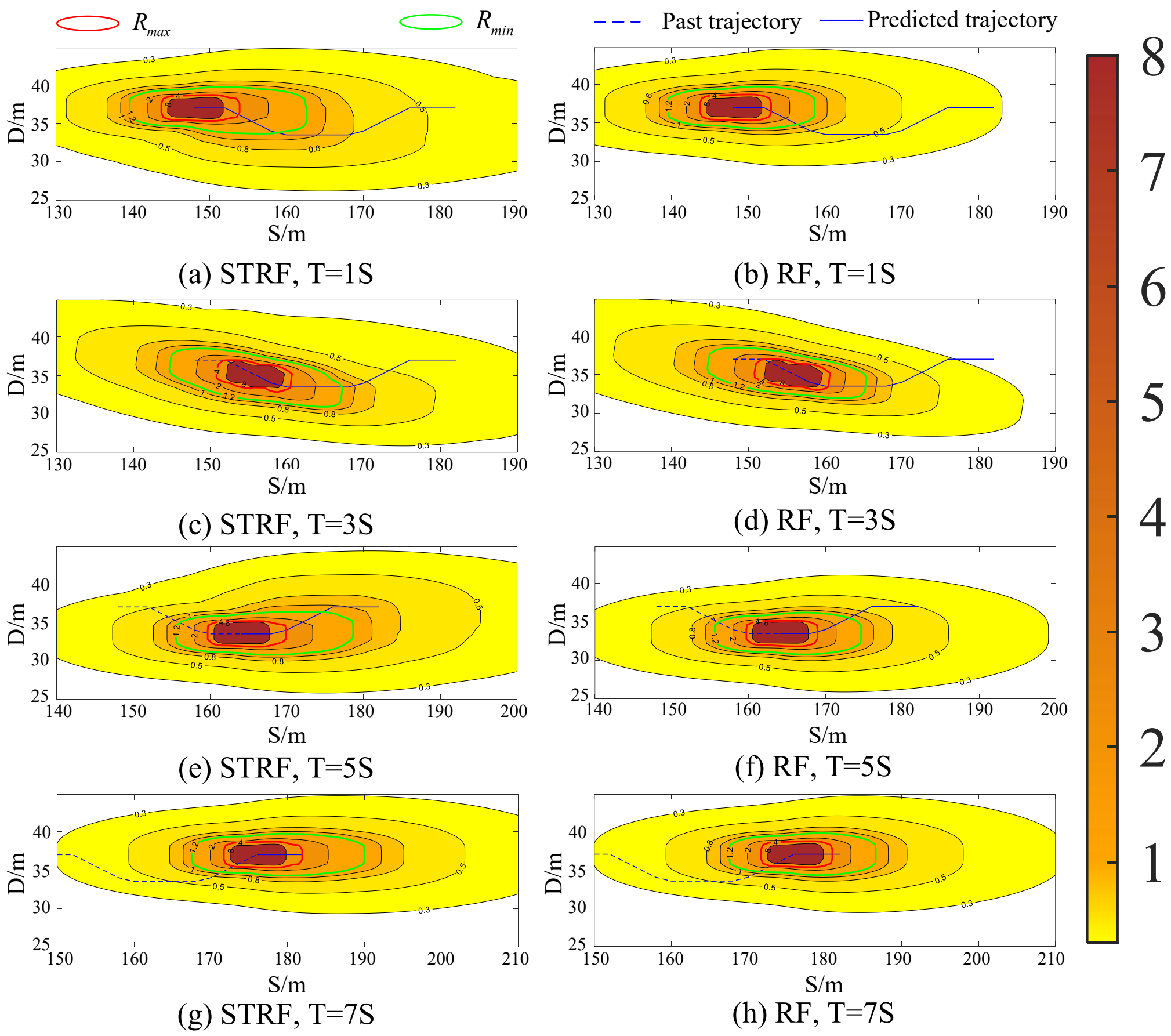

**Fig. 13** STRF and traditional risk field.

*4.5.2 STRF distribution: example of a real-world group of vehicles*

Based on aerial video trajectory data, we extracted a vehicle group, consisting of one off-ramp vehicle and five surrounding vehicles. To illustrate the relationship between physical quantities and the distribution of field strength: **Fig. 14** presents the motion posture of the surrounding vehicles. **Fig. 15** presents the potential field distribution of each surrounding vehicle. **Fig. 16** shows the total potential field distribution formed by combining the three types of elements. Several patterns emerge:

(1) **Anisotropy:** When driving at approximately 15 m/s, the potential field takes on a wedge-shaped form with an elongated distribution. In the absence of future trajectory changes, the affected region follows the pattern "front ≫ rear ≫ side," aligning with actual traffic phenomena.

(2) **Velocity and acceleration bias:** Comparing Car 2 and Car 3, their speed curves nearly overlap, but their acceleration curves differ. Car 3 decelerates by approximately 1 m/s initially before stabilizing, whereas Car 2 experiences a slight, steady acceleration (significantly less than 1 $m/s^2$). This acceleration bias influences field strength distribution: Car 3's field strength exhibits near-symmetry due to the counteracting effects of velocity and acceleration biases, while Car 2's distribution skews forward due to its strong velocity bias and weak acceleration bias.

(3) **Velocity and acceleration can drastically change the field strength distribution:** Car 4 exhibits the most heterogeneous field strength distribution, with the strongest forward bias among the five vehicles due to its high velocity and maximum acceleration.

(4) **The values of $R_{max}$ and $R_{min}$ are consistent with traffic phenomena:** The superposition distribution of the total obstacle potential field can be seen at the bottom of **Fig. 15**. It can be seen that the ego vehicle is currently located exactly between $R_{max}$ and $R_{min}$, which confirms the self-consistency of the dynamic risk balance theory.

(5) **$R_{max}$ and $R_{min}$ can assist with driving decision-making:** Car 2, which exhibits minimal acceleration fluctuations, maintains a stable speed of 15 m/s. The upper and lower risk threshold can be obtained from the red and green lines. Within approximately 20 m in front and 6 m behind, the risk exceeds the upper threshold, classifying this area as hazardous. Furthermore, if the following vehicle remains more than 18 m behind, the risk falls below the lower threshold, indicating overly conservative driving. In such cases, reducing the following distance slightly could increase risk in exchange for improved efficiency—consistent with real-world traffic behavior and safety principles.

(6) **The strength distribution of road boundaries and lane lines: Fig. 16** (d) shows road boundaries and lane lines forming a wavy pattern, with low points at the road's center and high points at the dotted lines. However, these markings do not impose strict restrictions on vehicle movement. The road boundary, resembling a wall, exhibits a rapidly increasing field strength as vehicles approach it. At the boundary itself, the field strength surpasses a critical threshold, effectively preventing vehicles from crossing.

(7) **The strength distribution of the off-ramp in weaving segments:** In **Fig. 16** (c), the field strength within the off-ramp weaving segment applies only to vehicles driving on the mainline, in the mandatory lane-changing zone, and intending to exit. The inner three lanes exhibit a consistent pattern of decreasing field strength from front to back, encouraging vehicles to merge into the outermost lane as soon as possible, where no field strength restrictions exist.

(8) **3D potential field distribution of the total risk field: Fig. 16** (a) integrates the other three subgraphs. Based on this figure, the risk levels at different locations can be visualized intuitively, which supports the trajectory planning task in the following text.

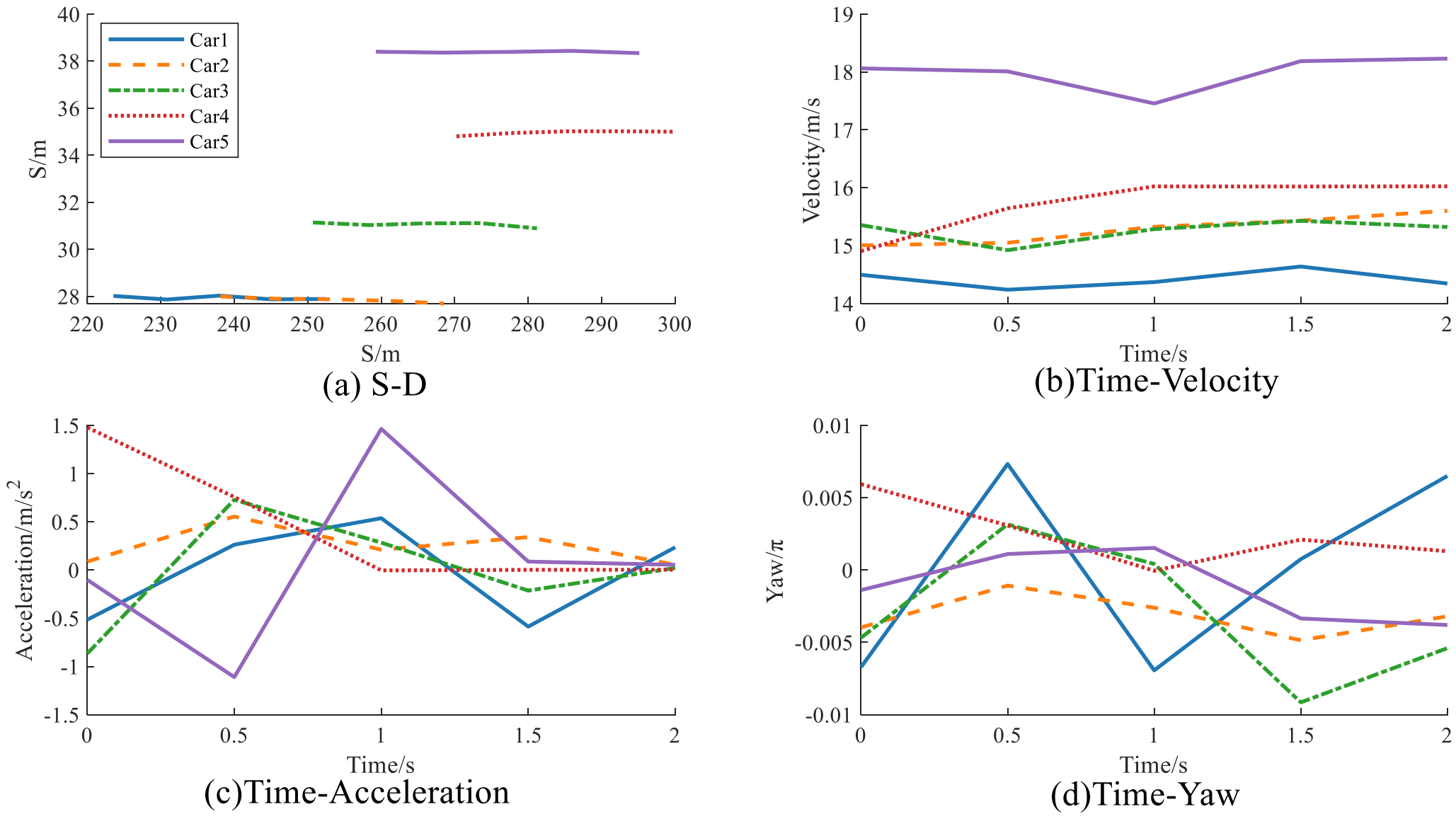

**Fig. 14** Surrounding vehicle motion posture.

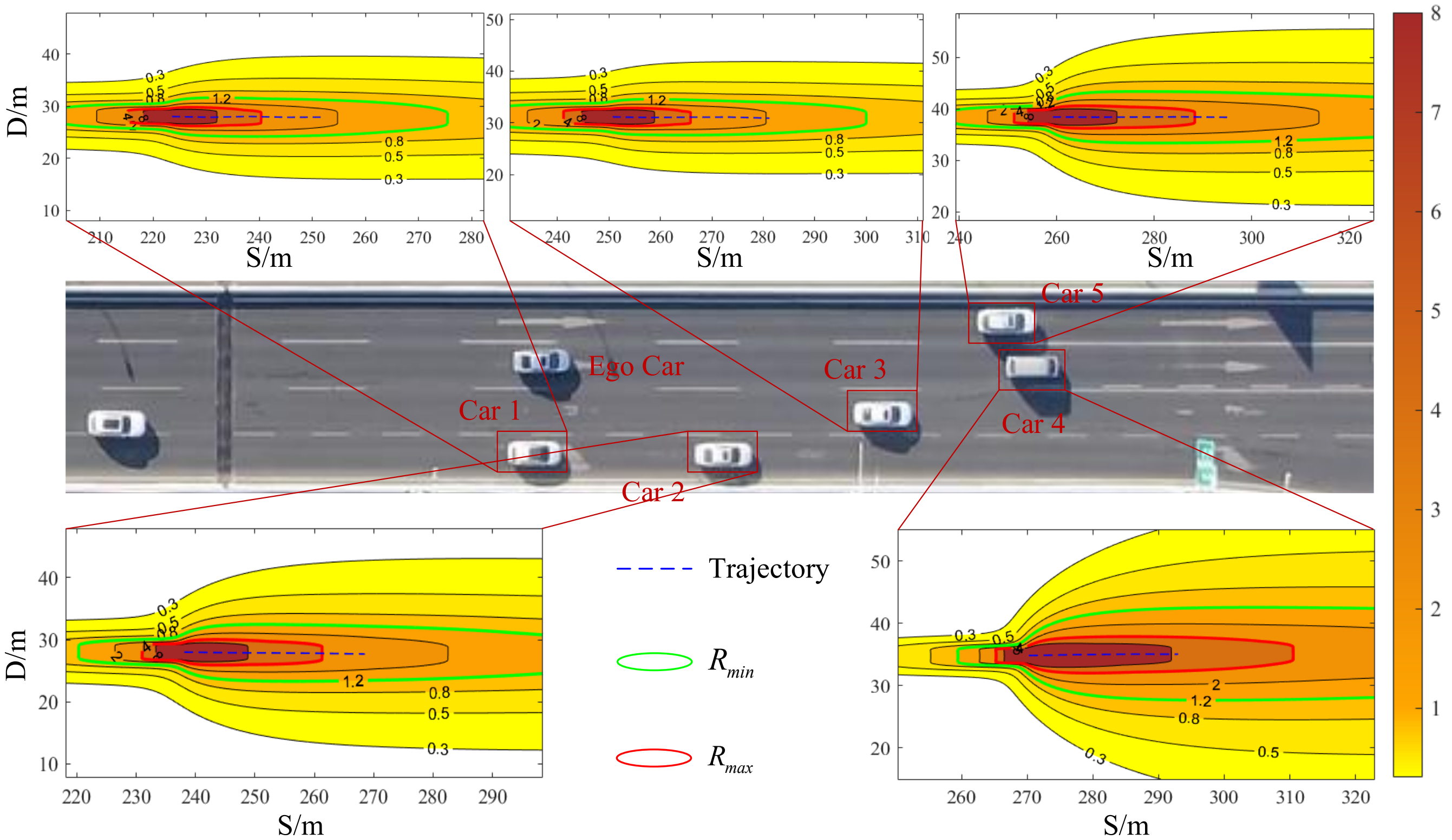

**Fig. 15** Typical vehicle groups with resulting obstacle risk field.

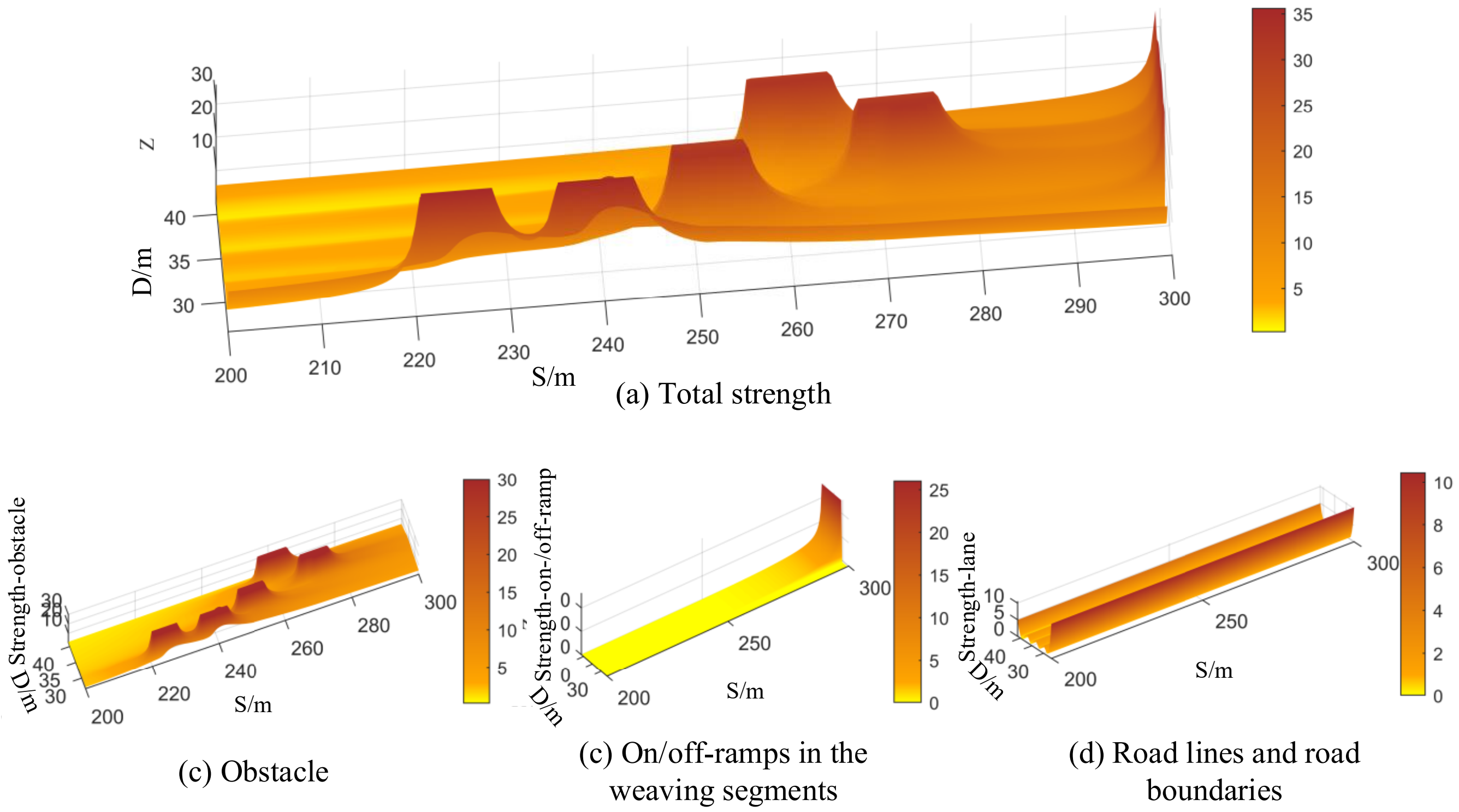

**Fig. 16** Three-dimensional potential field distribution of the total risk field.

## 5 A STRF-based CAV trajectory planning method in weaving segments

To establish a solid methodological bridge between the STRF safety assessment model constructed earlier and the specific trajectory planning framework in this section, it is first necessary to elucidate in detail how the continuous 3D risk field is effectively translated into actual constraints for underlying vehicle control. The STRF provides a macroscopic and precise quantification of global dynamic risks. It not only includes the geometric boundary constraints of static road structures but also proactively integrates the future movement trends of dynamic obstacles. However, applying this continuous and highly nonlinear 3D field directly as a cost function in numerical optimization imposes an immense computational burden. Furthermore, it easily causes the solver to fall into local optima, making it difficult to satisfy the strict real-time requirements of autonomous driving.

Therefore, this study introduces a transformation mechanism to refine the safety assessment results into a spatial-temporal risk occupancy map (STROM). Specifically, this transformation process involves three key steps. The first step is time domain slicing, where we discretize the continuous prediction horizon into multiple continuous time steps and extract the corresponding 2D planar risk field at each specific moment. The second step is spatial domain rasterization, which precisely maps the aforementioned 2D continuous planes into high-resolution discrete grid matrices, allowing each grid cell to carry the specific risk quantification value at that exact time and location. The third step is risk threshold binarization. By introducing a rigorously calibrated maximum allowable risk threshold from the previous stage, a strict logical judgment is applied to the grid matrices. Grids exceeding the threshold are directly marked as impassable, absolute high-risk areas, while grids below the threshold are designated as safe drivable areas.

Through the aforementioned slicing, rasterization, and binarization processes, the macroscopic safety assessment model is successfully reduced in dimensionality and transformed into a series of 2D safe feasible domains expanded chronologically. In subsequent trajectory planning tasks, this occupancy map acts as extremely stringent physical boundaries. It not only directly determines the spatial distribution range of feasible sampling points during the front-end dynamic iterative sampling phase but also explicitly provides the back-end quadratic programming solver with the specific mathematical boundaries for underlying obstacle avoidance hard constraints. This mechanism ensures that all subsequent trajectory sampling and optimization processes are fundamentally restricted to a safe domain free of collisions in both time and space, achieving a methodological alignment from dynamic risk situational awareness to microscopic vehicle safe trajectory generation.

In the traditional DP + QP hierarchical method, trajectory sampling points are selected based solely on the initial positions of obstacles. This static sampling method fails to account for dynamic obstacles over time, necessitating the use of an S-T diagram in high-speed scenarios to mitigate dynamic conflicts. However, incorporating S-T diagrams and differentiating between static and dynamic obstacles introduce the following two limitations:

(1) The trajectory planning completion time depends on the S-T diagram and speed planning, while the S-T diagram and speed planning are, in turn, dependent on path planning. If the initial path planning is not optimal, there may be cases where no feasible solution can be found on the S-T diagram. In addition, even if a solution exists, inefficiency issues caused by avoiding dynamic obstacles may arise, such as the extension of trajectory planning time.
(2) Since S-T diagram plotting and speed planning are contingent on path planning, speed planning can only commence once path planning is completed. This sequential dependency increases computational time and complexity.

To overcome the aforementioned limitations and deeply integrate the underlying safety assessment model, the proposed STROM comprehensively and deeply considers the risk impacts of static and dynamic obstacles directly during the initial sampling phase. This fundamental shift in underlying logic not only eliminates the reliance of the algorithm on serial ST diagrams but also significantly reduces the probability of generating infeasible trajectories. More crucially, because the dynamic sampling method explicitly delineates the absolutely safe feasible domain in advance through STROM, it enables the QP of paths and the speeds to achieve true independent parallel evolution. This parallel computing architecture reduces underlying computational demands and guarantees the absolute physical safety of intelligent vehicles in complex and high frequency weaving segment vehicle interaction environments.

5.1 Spatial-temporal risk occupancy map

Based on the predicted trajectory, we formulate a three-dimensional dynamic risk field using the Frenet coordinate system, as illustrated in **Fig. 17**. The resulting 3D dynamic risk field exhibits anisotropic density distributions. To optimize computational efficiency in trajectory planning, we integrate the concept of an occupancy grid map (OGM). By combining the STRF with OGM, we construct a STROM. The process involves the following steps:

**Step 1: Slicerization:** We segment the $T$ -axes into slices, with each slice's time domain matching that of the time domain of trajectory planning, yielding a slice-based representation of the 3D dynamic risk field.

**Step 2: Rasterization:** A uniform distribution strategy places sampling points on the road at a grid resolution of 0.5 m. Each sampling point assesses both dynamic and static risks within a 0.5 m × 0.5 m region, enabling efficient spatial risk evaluation.

**Step 3: Risk value assignment:** We assume a uniform risk value within each grid and replace the entire grid's risk value with that of its center point. We map this quantified risk value ped onto the rasterized grid.

**Step 4: Risk correction and constraint application:** Risks below a designated threshold are assigned a value of 0 (marked in orange in **Fig. 18** (b)), while risks above the threshold are assigned a value of 1 (marked in green in **Fig. 18** (b)). This results in a STROM with slices (**Fig. 18** (b)), serving as a constraint for defining dynamic sampling regions during trajectory planning.

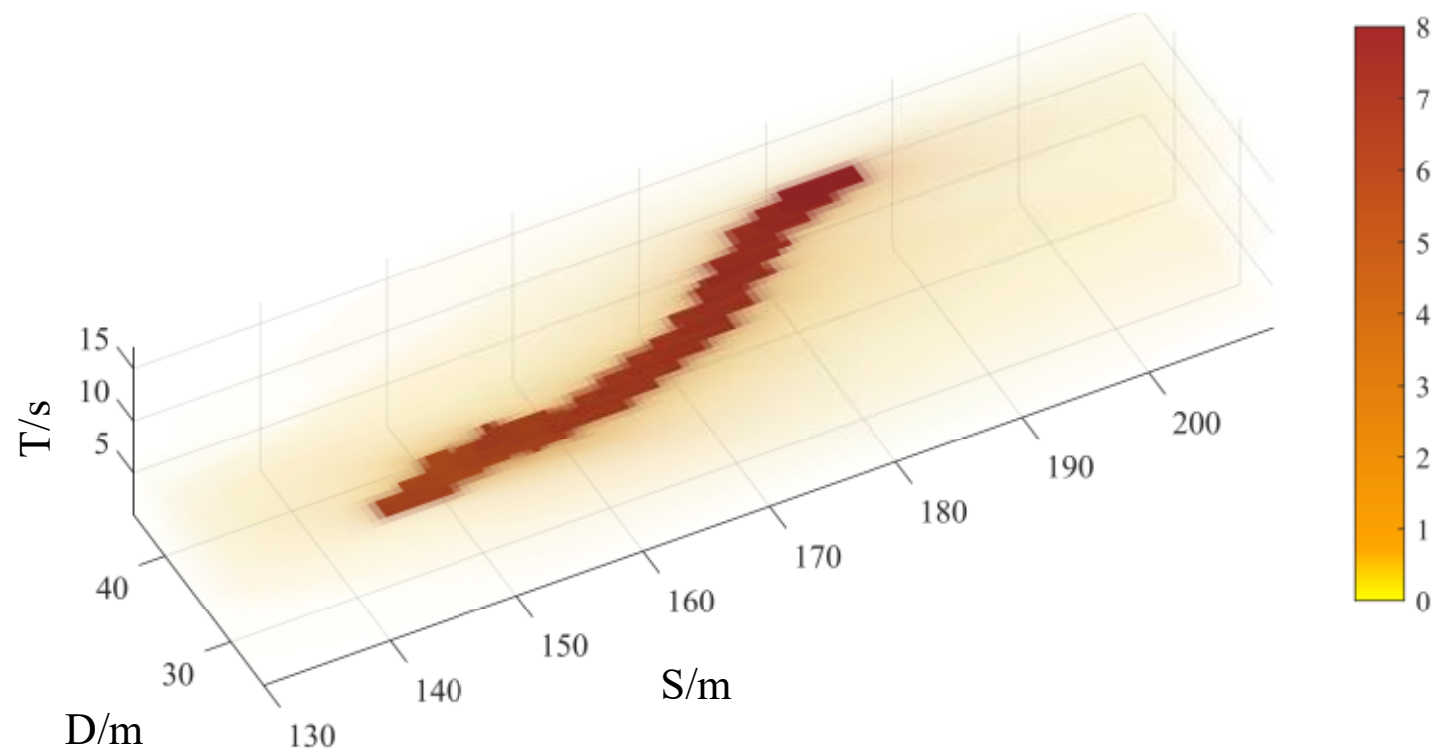

**Fig. 17** 3D dynamic risk field.

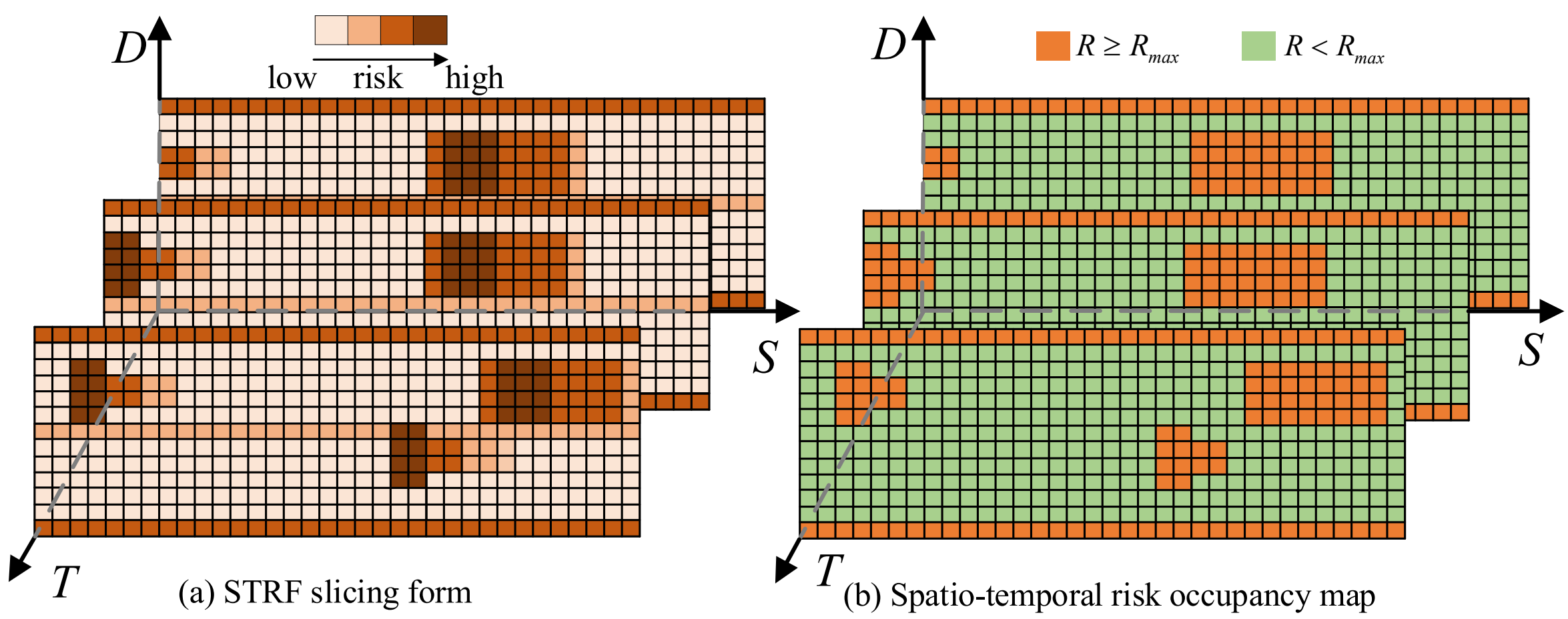

**Fig. 18** STRF slicing form and spatial-temporal risk occupancy map.

During trajectory planning, the system extracts STROM for specific time domains and sampling ranges. The planned trajectory must then avoid high-risk regions both spatially and temporally to ensure safety.

5.2 Dynamic iterative sampling method based on STROM

In conventional sampling, the location of each sampling point is determined based solely on the initial lane-change information. As mentioned at the beginning of **Section 5**, this sampling method can only avoid static obstacles, not high-speed dynamic moving obstacles. The implementation of S-T raises another set of issues. To accommodate both dynamic and static obstacles during the sampling phase, in our study, sampling points are dynamic, influenced by the STROM to incorporate more precise risk assessments during the sampling stage. This method ensures that each sampling point is not only linked to the previous sampling point but also constrained by the maximum feasible range of vehicle movement from the previous step. However, new sampling points in each step cannot be fully shared with other points in the same stratum because some neighboring steps may be inaccessible. To address this, we propose a global sampling approach and hierarchical attribution, which simplifies computations while fully accounting for STROM and adaptive sampling regions. The sampling framework, illustrated in **Fig. 19**, operates as follows:

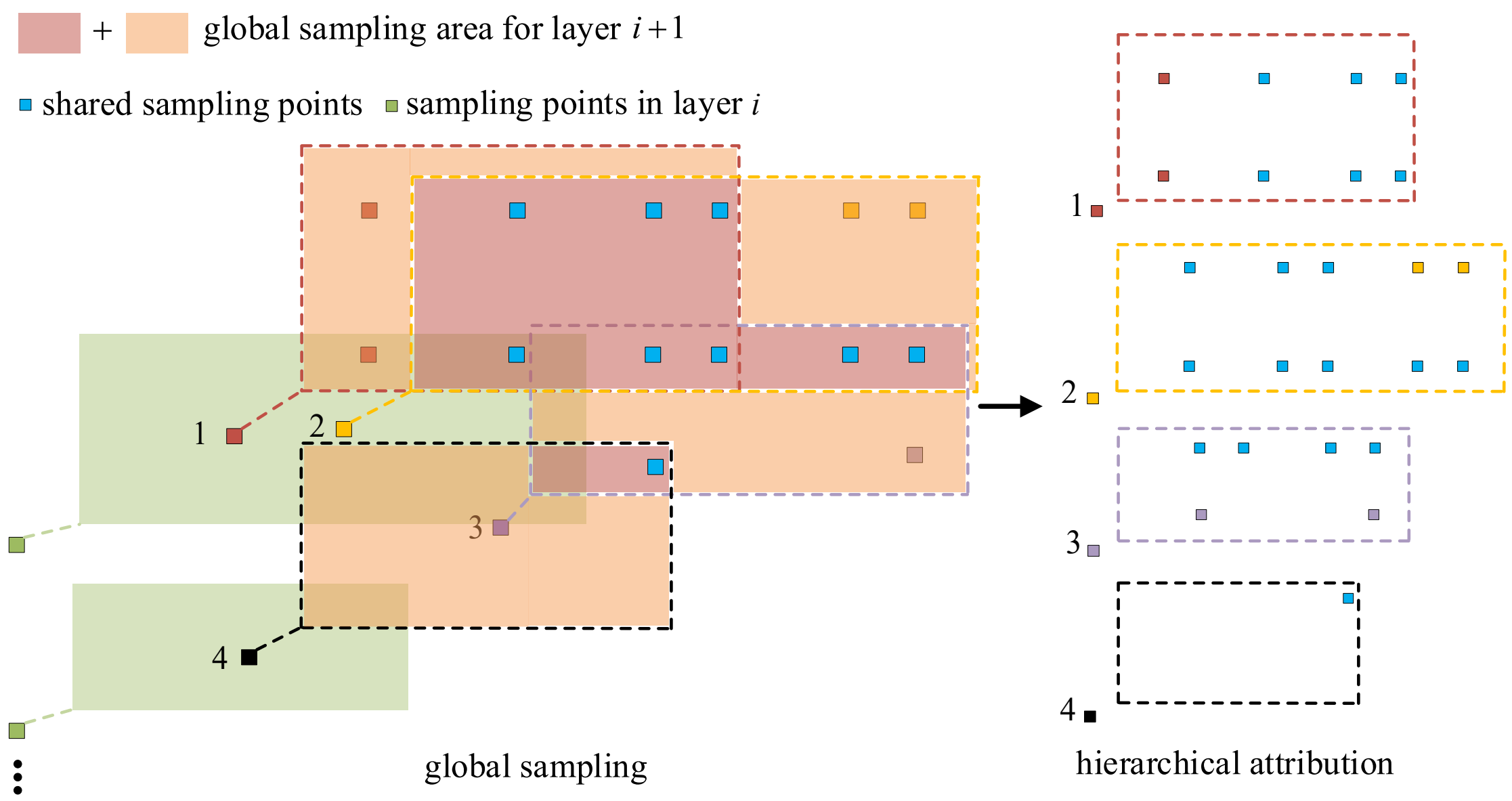

**Fig. 19** Global sampling and hierarchical attribution.

When we know the sampling points from the previous time domain (e.g., points labeled “1”, “2”, “3”, and “4” in **Fig. 19**), we first determine the maximum movement range for each point in the next time domain. We construct a global sampling range by union all individual sampling ranges. Once sampling is complete, we assign points in the next time domain based on their respective sampling ranges. The formula for calculating the maximum movement range will be detailed below (**Section 5.3.1**).

5.3 Optimal rough path generation

*5.3.1 Dynamic iterative sampling-based candidate trajectory generation strategy*

The sampling steps are schematically shown in **Fig. 20**. The candidate trajectory generation strategy, including five steps, is detailed as follows:

**Step 1: Determine the basic parameters.** At the initial moment, the vehicle decision module outputs the starting position of the lane change noted as $(s_0, d_0)$, and the longitudinal and transverse velocities are $v_{s,0}$, $v_{d,0}$, respectively.

**Step 2: Determining the range of STROM for the planning time domain.** Once the sampling points in the previous time domain are determined, each point must be assigned a range within the STROM for the next planning time domain. This range defines: **1)** The number of slices of the STROM selected. **2)** The area in which vehicles can move within a defined planning time domain on a given STROM. The $t^{th}$ planning time domain corresponds to the $t^{th}$ slice of the STROM. For lane-change trajectory planning, the area is modeled as a rectangle, where: the vehicle's motion is constrained to avoid opposite-direction movement. Each sampling point from the previous layer generates a sampling area, defined by the maximum and minimum distances the vehicle can move longitudinally and transversely in the given time domain. The mathematical formulation for these constraints is provided in Eqs. (24) ~ (29). At the initial moment, the STROM is defined by four endpoints: $(s_0 + s_{1,min}, d_0 + d_{1,min})$ , $(s_0 + s_{1,min}, d_0 + d_{1,max})$, $(s_0 + s_{1,max}, d_0 + s_{1,min})$, and $(s_0 + s_{1,max}, d_0 + d_{1,max})$. which correspond to the first slice. When the $t \geq 2$, the $n_{t-1}$ sampling points generated in the $t-1^{th}$ step, results in $n_{t-1}$ side-by-side and distinct STROM usage ranges. Consequently, the four endpoints of the $t^{th}$ slice are $(s_{t-1,k} + s_{t,k,min}, d_{t-1,k} + d_{t,k,min}), (s_{t-1,k} + s_{t,k,min}, d_{t-1,k} + d_{t,k,max}), (s_{t-1,k} + s_{t,k,max}, d_{t-1,k} + d_{t,k,min}), (s_{t-1,k} + s_{t,k,max}, d_{t-1,k} + d_{t,k,max})$.

$$s_{t,k,max} = \bar{v}_{s,t-1,k} \cdot t_D + \frac{1}{2} a_{s,t,k} \cdot t_D^2 \tag{25}$$

$$s_{t,k,min} = \bar{v}_{s,t-1,k} \cdot t_D - \frac{1}{2} a_{s,t,k} \cdot t_D^2 \tag{26}$$

$$d_{t,k,max} = \bar{v}_{d,t-1,k} \cdot t_D + \frac{1}{2} a_{d,t,k} \cdot t_D^2 \tag{27}$$

$$d_{t,k,min} = \bar{v}_{d,t-1,k} \cdot t_D - \frac{1}{2} a_{d,t,k} \cdot t_D^2 \tag{28}$$

$$\bar{v}_{s,t-1,k} = \begin{cases} v_{s,0}, t = 0 \\ \dfrac{s_{t-1,k} - s_{t-2,k}}{t_D}, t \geq 1 \end{cases} \tag{29}$$

$$\bar{v}_{d,t-1,k} = \begin{cases} v_{d,0}, t = 0 \\ \dfrac{d_{t-2,k} - d_{t-2,k}}{t_D}, t \geq 1 \end{cases} \tag{30}$$

where $s_{t,k,max}$ and $d_{t,k,max}$ are the maximum distances that can be reached by the vehicle longitudinally and horizontally a planning time domain at step $t$ in the Frenet coordinate system when the $k^{th}$ point sampled at step $t-1$ (denoted as $p_{t-1,k}$) is the sampling point, respectively. $s_{t,k,min}$ and $d_{t,k,min}$ are the corresponding minimum distances, respectively. $\bar{v}_{s,t,k}$ and $\bar{v}_{d,t,k}$ are the vehicle's longitudinal and transverse speeds, respectively, when the vehicle is at the step $t$ with the sampling start point $p_{t-1,k}$. In the first step, the actual state of the vehicle is known and is replaced by $v_{s,0}$ and $v_{d,0}$. When at step $t$ $(t \geq 2)$, the actual speed of the vehicle has not yet been planned and is replaced by the average speed planned at the previous step, Eqs. (29) and (30) explain this process. is the position of the point $p_{t-1,k}$. $(s_{t-1,k}, d_{t-1,k})$ is the location of the point where the previous step connects to $p_{t-1,k}$. When more than one point is connected, choose the nearest point. $t_D$ is the planning time domain. $a_{s,t,k}$ and $a_{d,t,k}$ are the maximum acceleration of the vehicle in the longitudinal and lateral directions, respectively, at $t$ step with the $k^{th}$ point as the sampling start point.

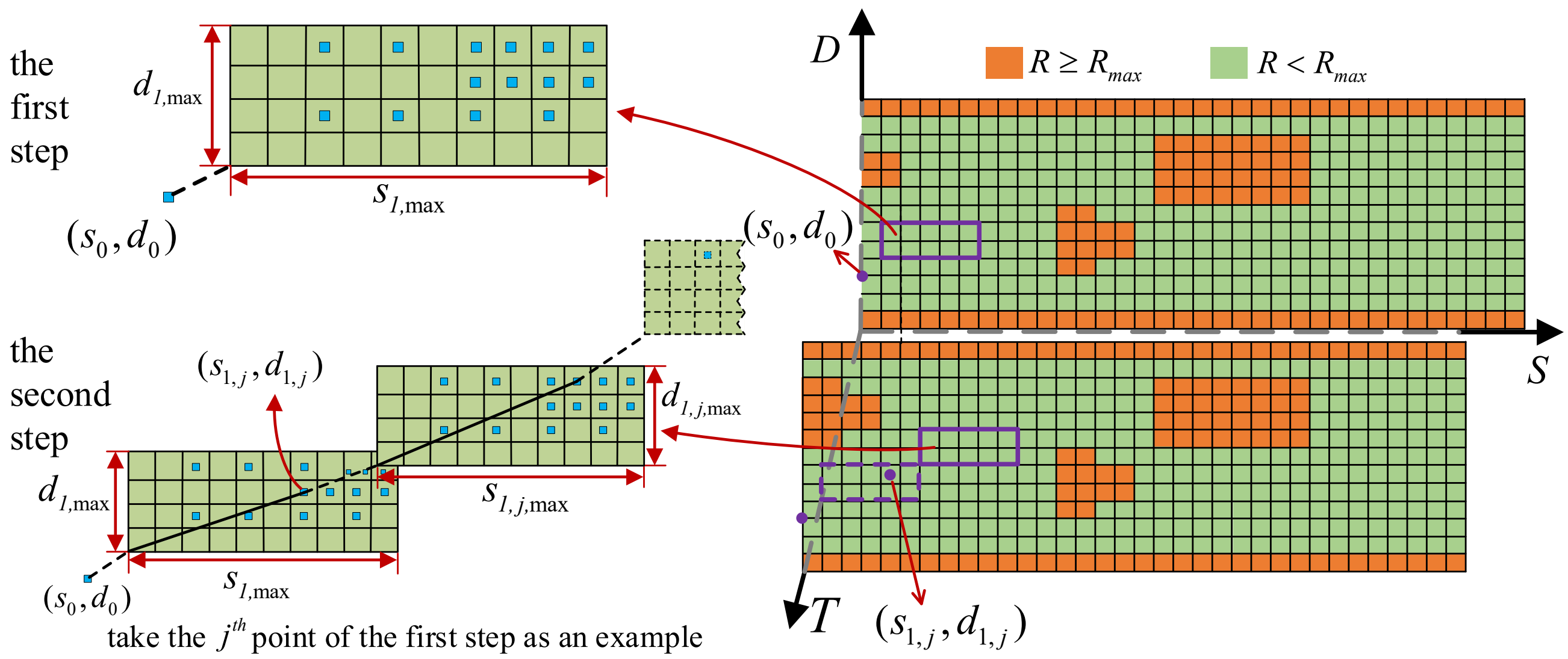

**Fig. 20** Sampling process based on STROM.

**Step 3: Dynamic iterative sampling.** The union of the STROM usage ranges of the $n_t$ sampling points in the time domain $t$ is used as the global sampling range. Taking the global sampling range as the reference, we apply a non-uniform sampling within this region to accelerate lane-change completion. Sampling is restricted to grid cells permitted by the risk constraints, ensuring that the generated samples comply with the risk constraints at the time of selection. When global sampling of the $t+1^{th}$ planning time domain is completed, the next set of sample points associated with each sample point is claimed using the hierarchical attribution. This is explained in **Section 5.2**.

**Step 4: Verify risk compliance.** For each newly generated sampling point, check whether the trajectory curve connecting it to the previous corresponding sampling point intersects with a high-risk area. If an intersection is detected, immediately discard the sampling point.

**Step 5: Termination condition.** Repeat **Steps 3** and **4** until the constrain (31) is satisfied, at which point the iteration terminates.

$$\delta \geq \frac{W}{2} - d_{t,k} \geq -\delta \tag{31}$$

where, $\delta$ is a positive number greater than 0, take 0.5 m. $W$ is the width of the road.

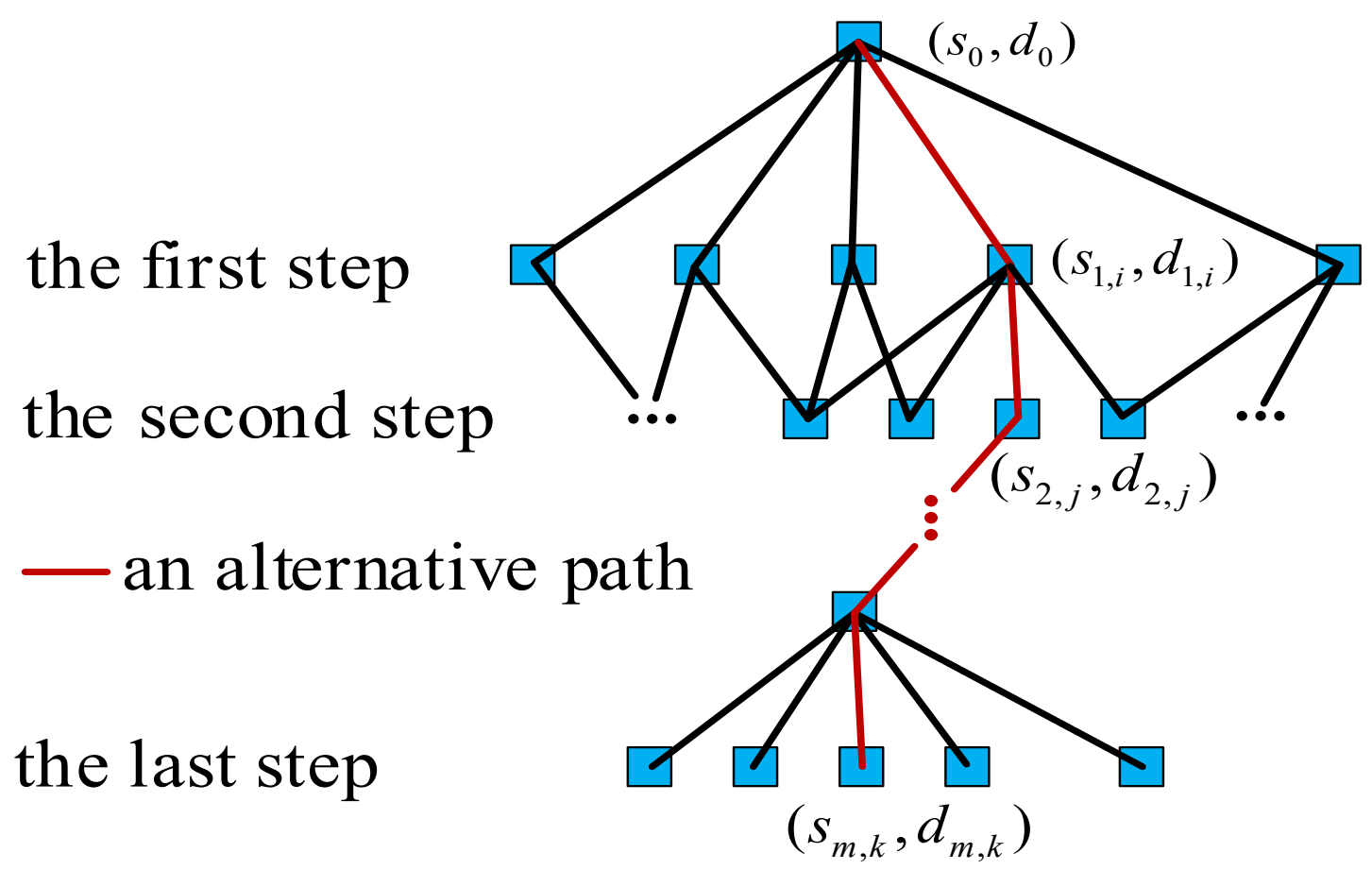

**Fig. 21** Set of candidate trajectories.

After the above five steps of sampling, connecting the sampling points associated with each index, as shown in **Fig. 21**. A total of paths is formed. $n_t$ is the number of samples at step $t$, and $m$ is the number of steps.

*5.3.2 Path evaluation based on dynamic programming*

We obtain the optimal rough path through cost function evaluation, which takes into account efficiency, vehicle dynamics constraints, and driving comfort. To achieve a balanced evaluation, the cost function for the $p^{th}$ path is divided into three components. They are provided in Eqs. (32) to (36).

$$J_p = \omega_{eff} \cdot J_{eff,p} + \omega_{dyn} \cdot J_{dyn,p} + \omega_{smo} \cdot J_{smo,p} \tag{32}$$

$$J_{eff,p} = m \tag{33}$$

$$J_{dyn,p} = \begin{cases} \sum_{t=2}^{m-1} \kappa(p_i), \kappa(p_i) \le \kappa_{max} \\ \inf, \kappa(p_i) > \kappa_{max}, \exists p_i \end{cases} \tag{34}$$

$$\kappa = \frac{1}{R} = \frac{a \cdot b \cdot c}{S} \tag{35}$$

$$J_{smo,p} = \sum_{t=1}^{m} \left[ (d_t - d_{t-1}) - (d_{t-1} - d_{t-2}) \right]^2 \tag{36}$$

where $J_p$ is the cost function evaluation index, $J_{eff,p}$ and $\omega_{eff}$ are the efficiency index of the $p^{th}$ path and its weight, respectively, $J_{dyn,p}$ and $\omega_{dyn}$ are the curvature cost of the $p^{th}$ path and its weight, respectively. Only when the curvature on the path meets the maximum curvature limit, it will enter the cost evaluation session, otherwise, it will be eliminated directly. The Eq. (34) shows the process in this way. $J_{smo,p}$ and $\omega_{smo2}$ are the comfort indexes and their weights of the $p^{th}$ path, and $m$ is the total step length of the $p^{th}$ path, the smaller the value means the shorter the time of lane-changing, which represents the efficiency. $\kappa(p_i)$ is the curvature of the $i^{th}$ scattering point on the $p^{th}$ path, the local curvature is calculated using the outer circle formed by this point and the scattering points before and after it, $R$ is the radius of the outer circle, $a, b, c$ are the distances between the two of the $i-1^{th}, i^{th}, i+1^{th}$ scatters, respectively, $S$ is the area of the triangle formed by these three scatters, and $J_{smo2,p}$ denotes the lateral mean velocity fluctuations between neighboring steps.

Algorithms for solving such problems are well established, using DP algorithms to solve for optimal rough paths.

5.4 Path smoothing base on quadratic programming and STFOM

We employ a segmental jerk approach to fit curves, assuming that the third-order derivative of a curve connecting two scattering points between them is always constant. In other words, second-order derivatives are continuous, while higher-order derivatives (third and above) are discontinuous. We propose a quadratic path planning method that leverages dynamic time domains and STROM, as illustrated in **Fig. 22**. The key idea is to use STROM as constraints while allowing horizontal positions of sampling points to vary, keeping the vertical positions fixed. Since horizontal variation is significantly greater than vertical variation, this ensures that each time domain has a unique and distinct STROM. This secondary planning process optimizes both curve comfort and risk minimization. The proposed approach offers several advantages:

(1) **Enhanced safety:** Traditional sampling methods treat vehicles as point masses, neglecting their

geometric contours. The proposed QP explicitly accounts for vehicle geometry, effectively eliminating this risk.

(2) **Reduced computational time:** By keeping the vertical distance between points constant while allowing horizontal distance variations, the speed optimization problem becomes completely known. Unlike traditional serial computation, where speed optimization depends on path planning, our method enables parallel computation of path and speed optimization, reducing computational overhead.

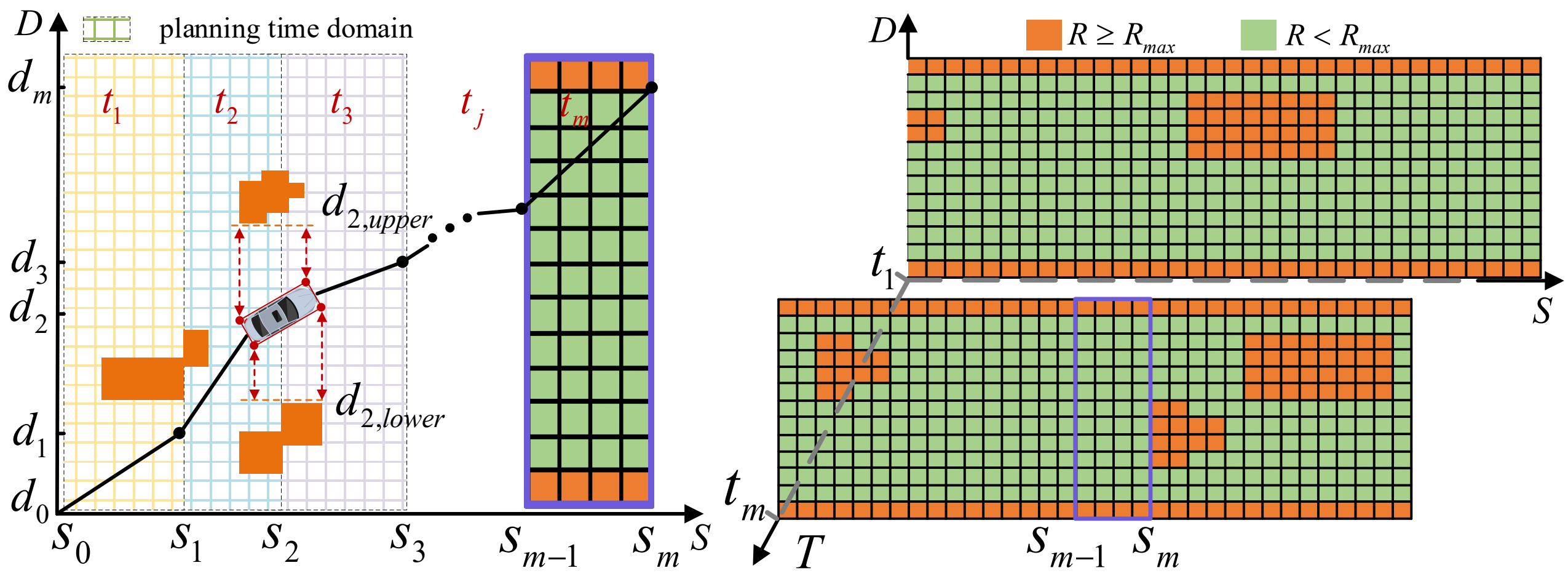

**Fig. 22** Schematic of secondary planning of paths based on STROM.

Combined with the uniqueness of this problem, the corresponding cost functions, continuity constraints, and inequality constraints have the following form:

(1) **General form:** the general form of the secondary plan is shown in formula (37) ~ (40).

$$\arg\max_{X} \frac{1}{2} X^T H X + f^T X \tag{37}$$

$$A_{eq} X = 0 \tag{38}$$

$$a \leq A_{ineq} \leq b \tag{39}$$

$$X = \left( d_1 \quad d_1^{'} \quad d_1^{''} \quad \cdots \quad d_m \quad d_m^{'} \quad d_m^{''} \right)^T \tag{40}$$

where, $A_{eq}$ denotes the coefficients of the equation constraints, $A_{ineq}$ denotes the coefficients of the inequality constraints, $a$ and $b$ are the upper and lower bounds of the inequality constraints, respectively, $X$ is the variable of optimization, i.e., the position of each scattering point, the first-order derivatives, and the second-order derivatives are optimized, respectively. When these variables are known, combined with the basic assumption that the third-order derivative is constant, the shape of the curve is known.

(2) **Equation constraints:** 1) The first-order Taylor expansion of $d_{i+1}$ and $d_{i+1}^{'}$, at $d_i$ and $d_i^{'}$, respectively, yields Eqs. (41) and (42), which are deformed to give the equation constraints (43) and (44). 2) To make the vehicle at the end point of the trajectory planning move parallel to the direction of the reference line, the point at the end point should satisfy the equation constraints as in Eq. (45).

$$d_{i+1} = d_i + d_i^{'} \cdot (s_{i+1} - s_i) + \frac{1}{2} d_i^{"} \cdot (s_{i+1} - s_i)^2 + \frac{1}{6}(d_{i+1}^{"} - d_i^{"}) \cdot (s_{i+1} - s_i)^2 \tag{41}$$

$$d_{i+1}^{'} = d_i^{'} + d_i^{"} \cdot (s_{i+1} - s_i) + \frac{1}{2}(d_{i+1}^{"} - d_i^{"}) \cdot (s_{i+1} - s_i)^2 \tag{42}$$

$$d_i + d_i^{'} \cdot (s_{i+1} - s_i) + \frac{1}{3} d_i^{"} \cdot (s_{i+1} - s_i)^2 - d_{i+1} + \frac{1}{6} d_{i+1}^{"} \cdot (s_{i+1} - s_i)^2 = 0 \tag{43}$$

$$d_i^{'} + \frac{1}{2} d_i^{"} \cdot (s_{i+1} - s_i) - d_{i+1}^{'} + \frac{1}{2} d_{i+1}^{"} \cdot (s_{i+1} - s_i)^2 = 0 \tag{44}$$

$$d_m^{'} = 0 \tag{45}$$

(3) **Inequality constraints:** 1) In the path evaluation stage, we introduce an evaluation function to limit the curvature at discrete sampling points. This constraint facilitated the optimization process, making it easier to obtain an optimal solution in quadratic programming. Unlike the previous approach, where the curvature was approximated using an outer circle, the new method directly computes the real curvature for more precise trajectory planning. The curvature-related inequality constraints are formulated in constrain (46). 2) In the sampling process, we consider the STROM constraints. However, the vehicle was simplified as a mass point. While this ensured that the geometric center satisfied the risk constraints, it did not account for the actual vehicle contour, potentially exposing the vehicle to higher risks. To address this gap, we must also make the four endpoints of the vehicle (representing its physical boundary) satisfy the risk constraints simultaneously. **Fig. 22** illustrates this approach, and the corresponding inequality constraints are defined in constrains (47) and (48).

$$\kappa(d_i) < \kappa_{max} \tag{46}$$

$$min(d_{i,p1}, d_{i,p2}, d_{i,p3}, d_{i,p4}) \geq d_{i,lower} \tag{47}$$

$$max(d_{i,p1}, d_{i,p2}, d_{i,p3}, d_{i,p4}) \leq d_{i,upeer} \tag{48}$$

(4) **Cost function:** Since we have considered the efficiency in the sampling stage, the completion time of trajectory planning is consistent with the completion time given in the sampling stage. The efficiency is not optimized separately here. The cost function mainly considers the comfort of trajectory planning, as shown in Eq. (49).

$$cost\ function = \omega_{dl} \cdot \sum_{i=1}^{m} d_i^{'2} + \omega_{ddl} \cdot \sum_{i=1}^{m} d_i^{"2} + \omega_{dddl} \cdot \sum_{i=1}^{m} (d_{i+1}^{"} - d_{i+1}^{"})^2 \tag{49}$$

where $\omega_{dl}$, $\omega_{ddl}$ and $\omega_{dddl}$ are the weights of the first-, second-, and third-order derivatives of $d$, respectively.

The above constraints and cost functions form the standard form of convex QP, and methods for solving convex quadratic programming problems are well established. Using the MATLAB Cplex toolbox, it is possible to implement the optimization of $X$.

5.5 Speed smoothing: parallel computation with path smoothing

In **Section 5.3**, the time intervals between adjacent points of the obtained optimal path remain constant, all equal to the sampling time window, but the distances between them are different. Therefore, the average speed between each adjacent step is different, as can be seen in the red dashed line in **Fig. 23**.

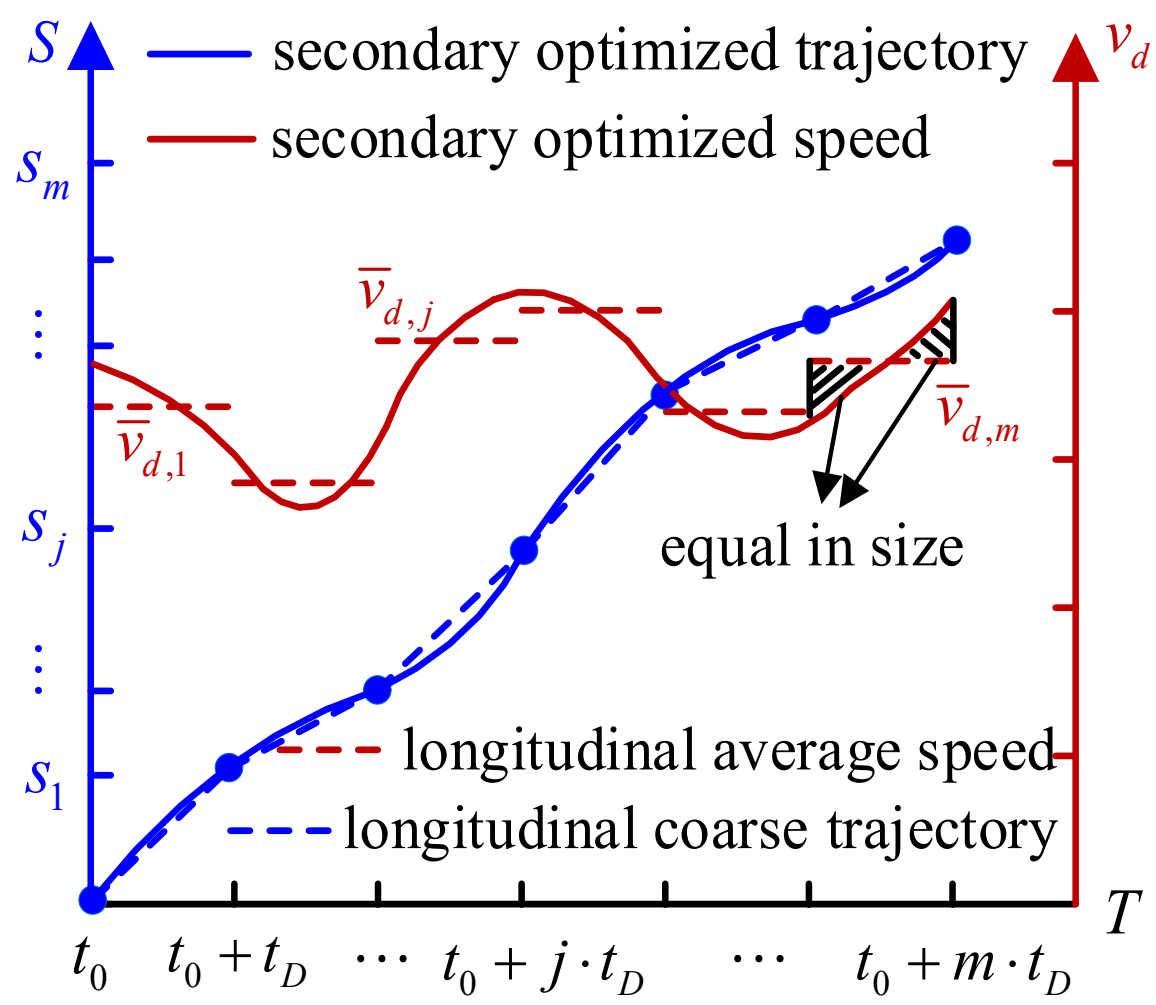

**Fig. 23** Schematic of quadratic programming of velocity with segmental mean velocity.

From **Fig. 23**, it's evident that this method leads to speed discontinuities between neighboring points, which are difficult to control in real-world driving. To address this gap, we introduce a speed smoothing method based on QP, similar to the technique used for path smoothing. Since lane-change trajectory planning typically involves minor variations in transverse speed, the longitudinal speed is the primary focus of optimization. Once the longitudinal speed is optimized, it is used—along with the known path—to derive the transverse speed through inverse computation.

Without going into too much detail due to the similarity with the path planning approach, the variables optimized by this process are $Y=\left(s_1 \quad s_1^{'} \quad s_1^{''} \quad \cdots \quad s_m \quad s_m^{'} \quad s_m^{''}\right)^T$. The form of the equational constraints, inequality constraints, and cost function are given directly as follows:

(1) **Equation constraints:** 1) Similar to the Taylor expansion for path constraints, the derivation is not shown again, and constraints (50) and (51) express the process. 2) To satisfy the risk constraint, the curve needs to pass strictly through the sampling points. The additional equational constraint is obtained as Eq. (52).

$$s_i + s_i^{'}\cdot(t_{i+1}-t_i)+\frac{1}{3}s_i^{''}\cdot(t_{i+1}-t_i)^2 - s_{i+1}+\frac{1}{6}s_{i+1}^{''}\cdot(t_{i+1}-t_i)^2=0 \tag{50}$$

$$s_i^{'}+\frac{1}{2}s_i^{''}\cdot(t_{i+1}-t_i)-s_{i+1}^{'}+\frac{1}{2}s_{i+1}^{''}\cdot(t_{i+1}-t_i)^2=0 \tag{51}$$

$$s(t_i)-s_{i_profile}=0 \tag{52}$$

where $(t_i, s_{i_profile})$ is the position of the $i^{th}$ point sampled.

(2) **Inequality constraints:** The main constraints are placed on the vehicle's dynamics parameters, as shown in the constraints (53) to (55).

$$0 \le \dot{s}_i \le v_{max} \tag{53}$$

$$-Dec_{max} \le \ddot{s}_i \le Acc_{max} \tag{54}$$

$$s_{i+1} - s_i \ge 0 \tag{55}$$

(3) **Cost function:** Like the cost function of path QP, it mainly measures the comfort, as shown in Eq. (56).

$$cost\ \ function = \omega_{ds} \cdot \sum_{i=1}^{m} s_i^{'2} + \omega_{dds} \cdot \sum_{i=1}^{m} s_i^{''2} + \omega_{ddds} \cdot \sum_{i=1}^{m} \left( s_{i+1}^{''} - s_{i+1}^{''} \right)^2 \tag{56}$$

Using the MATLAB Cplex toolbox, the optimization of the $Y$.

**6 Case study and discussion: trajectory planning performance analysis**

To assess the safety, real-time performance, efficiency, and comfort of the proposed trajectory planning method, we conducted a series of comparative experiments. **Section 6.1** outlines the sources of all trajectory planning cases, the comparative methods, and the parameter settings. Building on this foundation, **Section 6.2** provides a comprehensive analysis of the performance of the three methods across four scenarios, emphasizing the proposed method's advantages in efficiency and comfort under strict safety constraints. **Section 6.3** presents the running time of the algorithm across all cases, demonstrating that the proposed method meets the stringent real-time requirements of autonomous driving trajectory planning modules. Finally, to validate the accuracy of model calibration, **Section 6.4** details a sensitivity analysis on the maximum allowable risk $R_{max}$, exploring how variations in $R_{max}$ influence driving styles.

6.1 Simulation setup

Using the YOLO tool, we extracted 30 cases from aerial trajectory videos. Each case includes an ego vehicle (tasked with executing an on-ramp or off-ramp maneuver) and several surrounding interfering vehicles. Based on the number of interfering vehicles and the behavior of the Ego vehicle, we categorized these cases into four main types: on-ramp, simple non-congested scenarios (10 cases), off-ramp, simple non-congested scenarios (10 cases), on-ramp, complex congestion scenarios (5 cases), and off-ramp, complex congestion scenarios (5 cases). The parameter values are listed in **Table 6**. Where the bounds of acceleration and deceleration are referenced from the study of Gu et al. (Gu et al., 2024) and Wang et al. (Wang et al., 2024c). In addition, since the scenario we study is an urban expressway, in order to be consistent with the real speed limit of the case scenario, we set the speed limit to 80 $km/h$, i.e., 22 $m/s$.

**Table 6** Parameter values.

| Symbol | Description | Value |
|---|---|---|
| $v_{max}$ | Maximum speed | 22 $m/s$ |
| $Acc_{max}$ | Maximum acceleration | 4 $m/s^2$ |
| $\text{Dec}_{max}$ | Maximum deceleration | 6 $m/s^2$ |
| $\kappa_{max}$ | Maximum allowable curvature | 2 $m^{-1}$ |
| $R_{max}$ | Maximum allowable risk | 4 |

To comprehensively and objectively evaluate the superiority of the proposed algorithm, we introduced multiple comparison schemes in the aforementioned extracted real-world cases. The specific algorithm settings and characteristics are as follows:

(1) The proposed method in this study. This refers to the trajectory planning method driven by the STRF, which deeply integrates the dynamic risk assessment mechanism into the underlying trajectory sampling and optimization process, achieving high efficiency and safe obstacle avoidance.
(2) Real human driving trajectories. These are the actual driving trajectories of human drivers extracted directly from the aerial dataset. These trajectories serve as the absolute baseline reference under natural driving conditions, utilized to evaluate the practical performance of the planning algorithms in terms of human-like behavior and driving efficiency.
(3) Baseline 1 (DP and QP algorithm). This represents the traditional hierarchical method combining dynamic programming and quadratic programming. As a classic baseline widely applied in the current autonomous driving field, its primary characteristic is to first generate a rough path based on static sampling through dynamic programming, and subsequently employ quadratic programming for smoothing.
(4) Baseline 2 (SAC algorithm). This refers to the weaving segment lane changing strategy based on the soft actor critic algorithm. To reflect the current cutting edge level of deep reinforcement learning in the continuous action space trajectory planning field, we introduced this algorithm as the second competitive baseline. To ensure the fairness of the comparison, we utilized the gym library to reproduce a 1:1 high fidelity simulation environment of the weaving segment, and constructed the SAC algorithm relying on the stable-baselines 3 library. Targeting the highly dynamic characteristics of the weaving segment, this advanced baseline adopted a multi task random selection mechanism to avoid policy overfitting, which involves taking all extracted vehicle interaction groups as training cases and randomly selecting them for training in each episode. This algorithm possesses a 0.25 second fine grained planning capability to ensure the dynamic continuity of the environment, and introduced surrounding vehicle observation noise during training to enhance the generalization capability. In terms of model architecture, its observation space is designed to be 33 dimensional, comprehensively covering the 6 dimensional longitudinal and lateral kinematic state of the ego vehicle, the 24 dimensional relative kinematic characteristics of the nearest 6 surrounding interfering vehicles, and the 3 dimensional lane changing task state including target deviation and remaining distance. Its action space directly outputs two dimensional continuous control commands for longitudinal and lateral acceleration, achieving underlying physical control. Regarding the reward function design, this baseline is configured with a mixed mechanism combining sparse rewards and dense guidance, simultaneously considering economic cruising speed maintenance, lane changing progress advancement, driving comfort, and safe collision avoidance termination penalties. By inputting various cases into the trained and converged model for testing, this baseline provides a powerful learning based reference for the performance comparison of the proposed method in this paper.

It is important to note that the planning method proposed in this paper and the simulation execution of Baseline 1 were all tested on a computer equipped with an Intel i5 13500 processor relying on the MATLAB 2023a software platform. Conversely, Baseline 2, which demands substantial computing power, has its environment setup, model training, and inference testing specifically deployed on a separate device. This device is equipped with an Intel Core 12900KF processor with a base frequency of 3.2

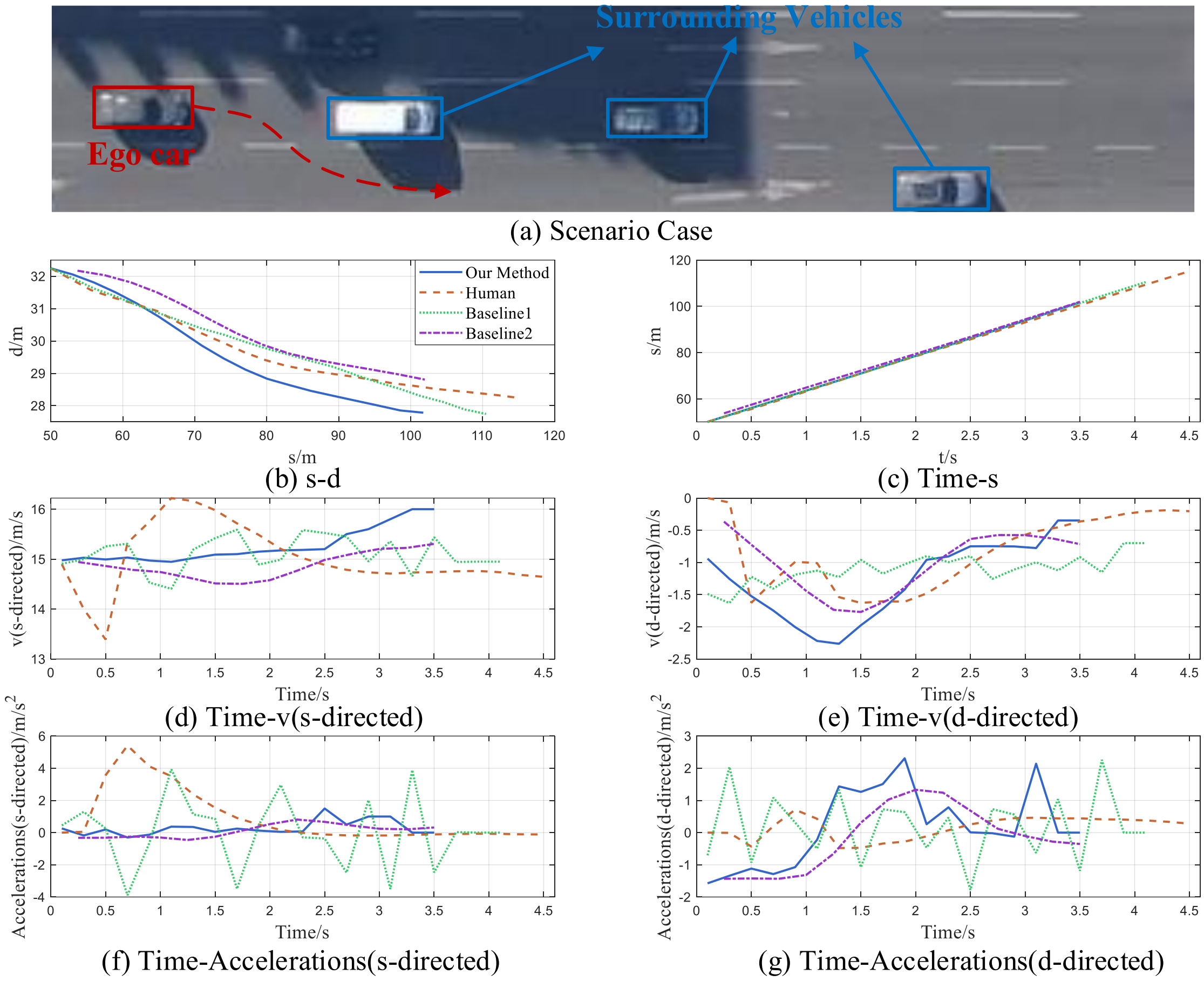

**Fig. 24** Case 1: off-ramp, simple non-congested scenario.

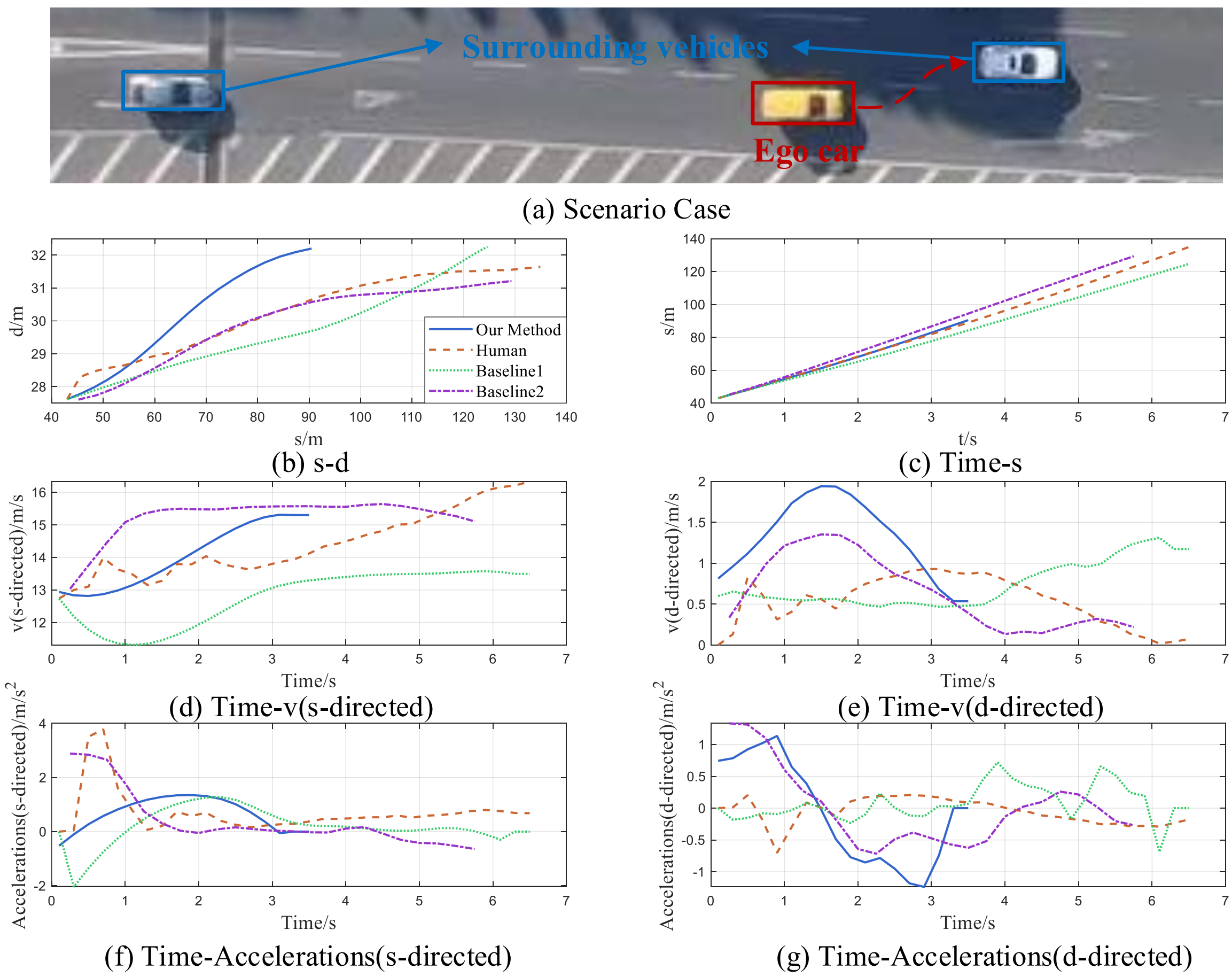

**Fig. 25** Case 2: on-ramp, simple non-congested scenario.

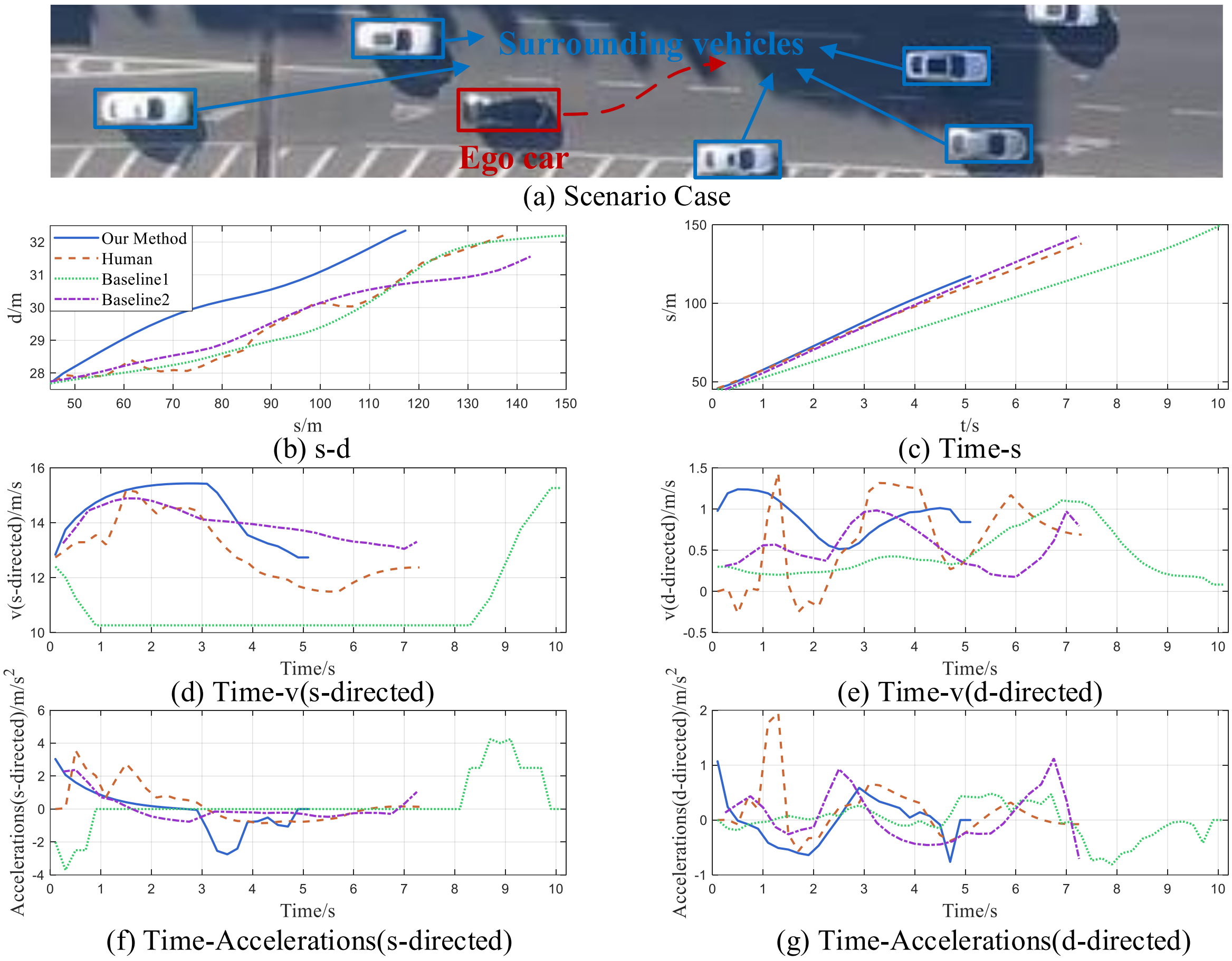

**Fig. 26** Case 3: on-ramp, complex congestion scenario.

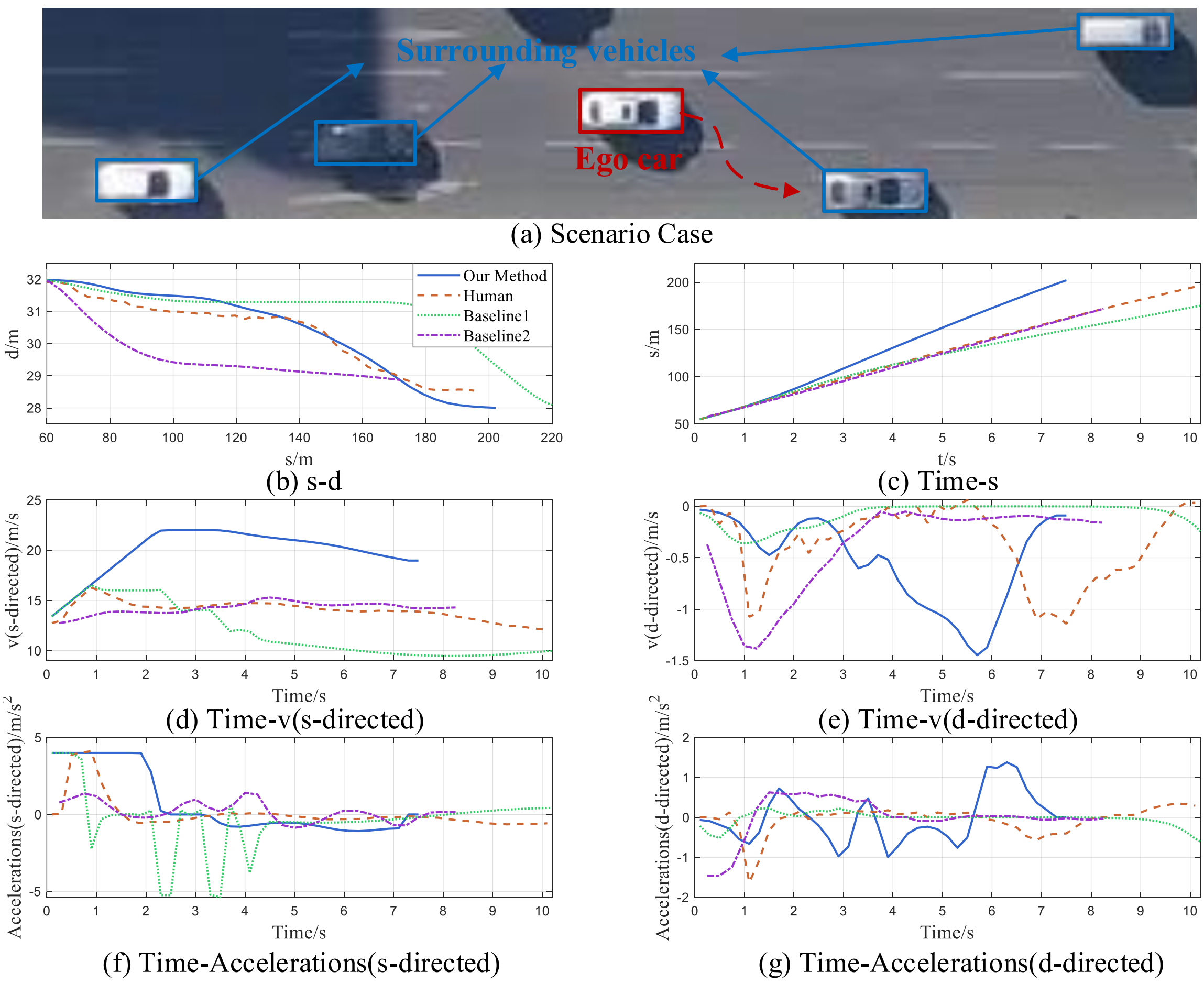

**Fig. 27** Case 4: off-ramp, complex congestion scenario.

GHz and an NVIDIA RTX 3090 graphics processing unit with 24 GB of video memory, and Python 3.8 is used to conduct the corresponding simulations.

6.2 Performance validation in multiple scenarios

We chose a typical case in each scenario to show and analyze, four cases with all the trajectory planning results contained in the three methods are shown in **Fig. 24** ~ **Fig. 27**. Furthermore, to display the effect of trajectory planning more clearly and intuitively, we have generated trajectory planning videos of the relevant simulations, and the link to these videos can be found in **Appendix A**. In **Fig. 24** ~ **Fig. 27**, the red dashed box highlights the ego vehicle, while the blue box represents. surrounding vehicles. The trajectories generated by our proposed method are represented by blue solid lines, the DP+QP method (Baseline 1) by green dashed lines, the real human driving trajectories by orange dashed lines, and the SAC method (Baseline 2) by purple dashed lines. A comparative analysis of these trajectories intuitively demonstrates the advantages of the proposed method in real-world driving environments.

*6.2.1 Performance validation: efficiency and comfort*

According to **Fig. 24** ~ **Fig. 27**, we can find: The proposed method demonstrates superior efficiency across all scenarios. Specifically, in Cases 1, 2, 3, and 4, our method achieves lane-changing completion times of 3.5 s, 3.5 s, 5.1 s, and 7.5 s, with corresponding average longitudinal velocities of approximately 15.1 m/s, 12.8 m/s, 14.3 m/s, and 19.3 m/s, respectively. These metrics rank our approach first among the four compared methods. Regarding the average lane-changing completion time, our algorithm reduces the duration by 30.98% compared to human drivers, 44.00% compared to the DP+QP method, and 20.48% compared to the SAC method. Furthermore, the average longitudinal velocity is 12.41% higher than that of human driving, 25.20% higher than the DP+QP method, and 8.20% higher than the SAC method. This efficiency enhancement is attributed to the STRF-based candidate trajectory generation strategy, which determines the lane-changing timing during the optimal coarse trajectory selection phase—a robust and controllable process. In contrast, the lane-changing time in the DP+QP method is heavily influenced by the S-T graph and the coarsely sampled path. When numerous infeasible regions emerge in the S-T graph due to dynamic obstacles, the obstacle avoidance maneuvers inevitably prolong the lane-changing duration. This phenomenon is particularly evident in Case 4, where the DP+QP method performs significantly poorly. The baseline determines its path solely depending on the risk level at the initiation of the lane change; once it selects a path that is optimal initially but becomes sub-optimal seconds later, it can only mitigate risks by adjusting its speed. Consequently, the velocity (see **Fig. 27** (d)) and acceleration (see **Fig. 27** (f)) of the baseline are highly undesirable, as it is forced to decelerate to counteract potential risks. Notably, although the SAC algorithm also demonstrates commendable efficiency—performing on par with our method in specific scenarios (e.g., both achieved nearly identical lane-changing completion times and average longitudinal velocities in Case 1, and the SAC method even surpassed our average longitudinal velocity in Case 2)—our method still maintains a higher overall efficiency than this data-driven approach, highlighting the immense application potential of our proposed framework.

The proposed method exhibits satisfactory comfort performance. Regarding longitudinal acceleration fluctuations, it maintains a controllable average acceleration and displays a smooth acceleration profile in most scenarios, indicating highly acceptable longitudinal comfort. A temporary degradation in comfort is observed only during the initial phase of the fourth case when the maximum acceleration is reached—a necessary trade-off between efficiency and comfort. The DP+QP method exhibits the most drastic velocity fluctuations, followed by human driving. This can be attributed to the baseline's tendency to prioritize safety over comfort during spatial-temporal planning, leading to aggressive speed adjustments to evade dynamic obstacles. The trajectories of human drivers are extracted from the YOLO videos and processed through a one-time Kalman filtering; slight errors might exist between consecutive time frames. Therefore, we will not extensively discuss the underlying reasons for their fluctuations. The SAC method also demonstrates satisfactory comfort, generally performing comparably to our method. Regarding lateral acceleration, although the proposed method is not optimal, it generally satisfies the comfort requirements: lateral acceleration is maintained between -1 $m/s^2$ and 1 $m/s^2$ for most of the trajectory, with a maximum absolute acceleration of less than 2 $m/s^2$.

*6.2.2 Performance validation: safety*

To verify the safety performance of the four algorithms, we analyzed case 1 and case 4 (introduced in **Section 6.2.1**) to characterize the risk levels endured by the ego vehicle across the three methods during the planned lane-changing cycle, thereby comparing their safety differences.

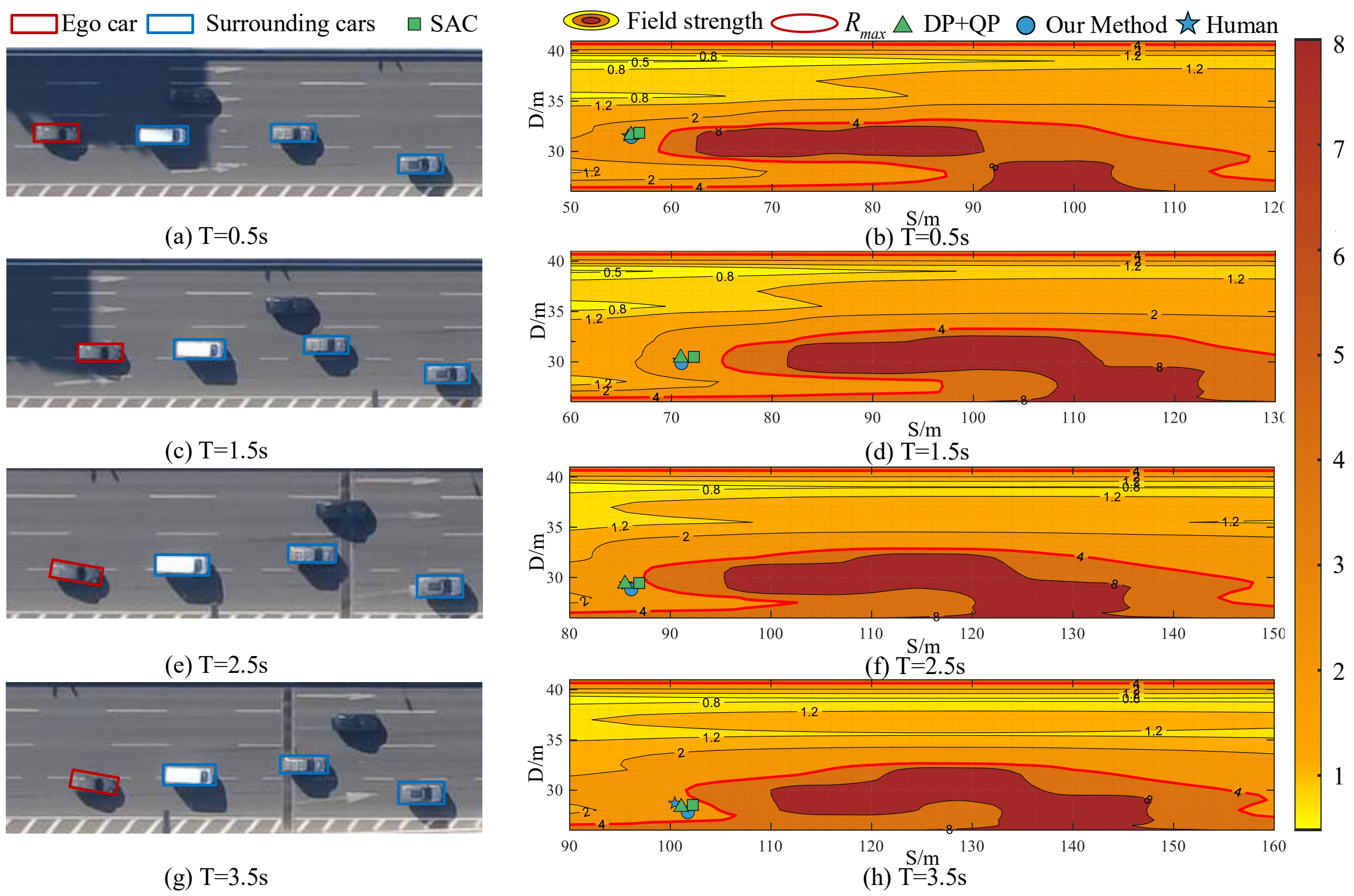

**Fig. 28** Trajectory planning safety analysis (Case 1).

In a simple scenario (case 1, as shown in **Fig. 28**), all surrounding vehicles were positioned ahead of the primary vehicle, with no rear-interfering vehicles and minimal multi-directional risk interference. As a result, all three algorithms maintained a relatively safe level at all times, indicating that human drivers, the baseline method, and the proposed method can all complete lane changes under safe conditions in simple scenarios. In complex multi-vehicle interaction environments (Case 4, as shown in **Fig. 29**), however, significant differences emerged. When the ego vehicle faced risk interference from multiple directions, human drivers struggled to attend to all risk sources simultaneously, leading to trajectories with measurable risk. The baselines also performed poorly: the DP+QP method lack robust risk quantification models, and their coarse sampling relies solely on initial planning-phase risk distributions, often resulting in suboptimal sampled trajectories. While velocity planning can still avoid obstacles, its effectiveness is limited in complex scenarios when relied on alone. Due to the absence of a more accurate risk quantification mechanism embedded within the state space construction and reward function formulation, the SAC method may yield a trained model that, while capable of avoiding collisions, remains susceptible to high-risk situations characterized by close inter-vehicle proximity. For example, as depicted in **Fig. 29**, the risk level of the ego vehicle persistently exceeds the threshold throughout the entire process. In contrast, the proposed method effectively leverages dynamic iterative sampling and speed planning to avoid high-risk regions across both spatial and temporal dimensions, regardless of scenario complexity. Additionally, as evidenced by **Fig. 27** and **Fig. 29**, our method achieves optimal efficiency under safe conditions.

In summary, the results indicate that the planned trajectories by our method not only effectively navigate risk at any given moment but also maximize efficiency and comfort. This is achieved through the STRF, which accounts for the influence of future obstacle trajectories on current risk, reducing constraints on the ego vehicle and minimizing its sensitivity to surrounding vehicle trajectory changes. The expanded field of view enables trajectory planning to adopt a global perspective, optimizing the overall strategy while maintaining smoother speed transitions and improving efficiency.

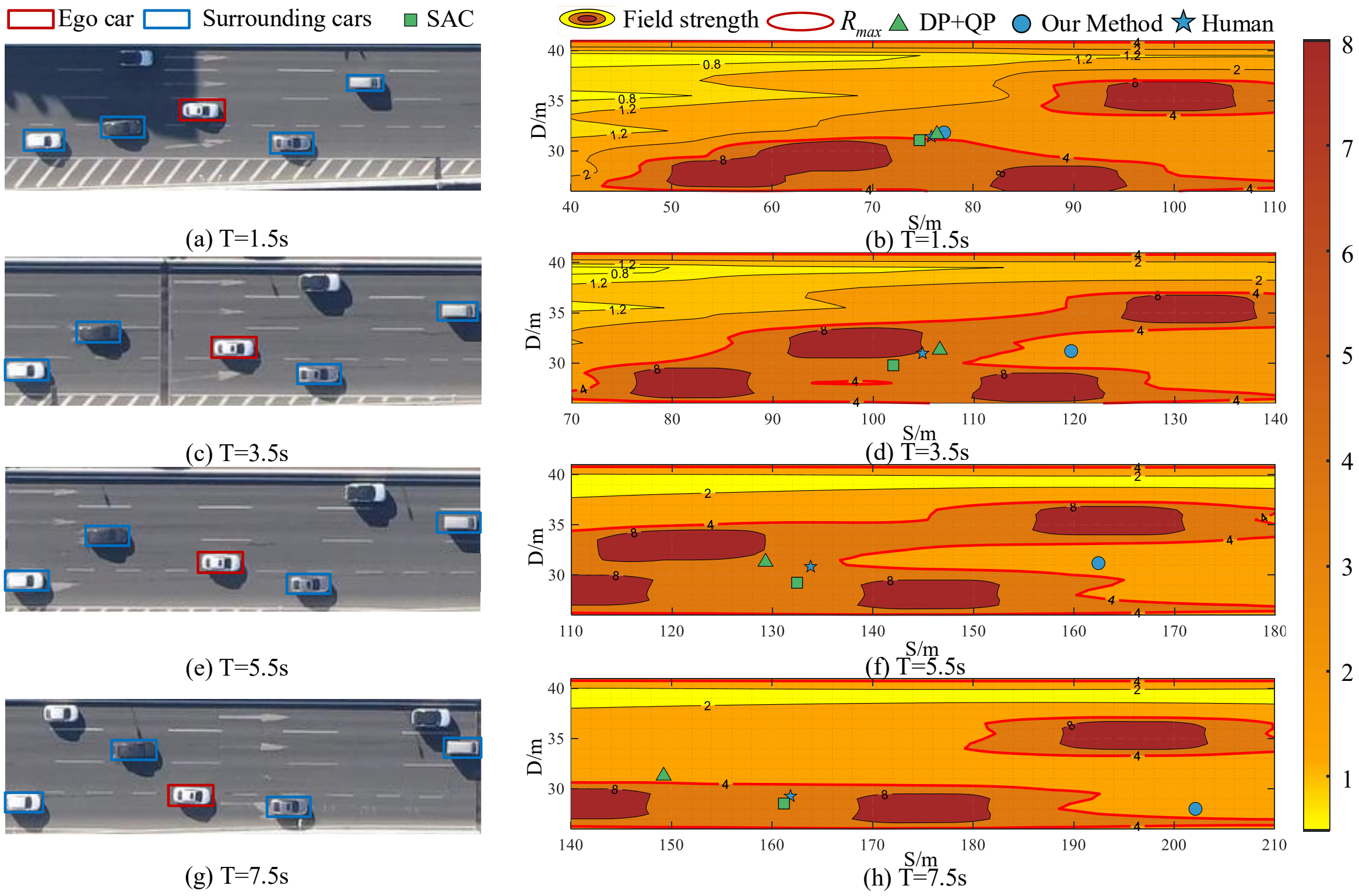

**Fig. 29** Trajectory planning safety analysis (Case 4).

6.3 Analysis of running time

We further evaluated the runtime of each module within the proposed algorithm across 30 cases. **Table 7** presents the average runtime for key trajectory planning modules, including the construction of the spatial-temporal risk occupancy map, dynamic iterative sampling, path evaluation, path smoothing, and speed smoothing. The results indicate the following:

(1) Despite being implemented in MATLAB without systematic running time optimization, the proposed algorithm achieves an average runtime of under 100 ms across all scenarios, within the acceptable range for practical application. Further reductions in running time are anticipated when deploying the code to in-vehicle systems.

(2) The proposed path-speed parallel optimization approach reduces running time, with savings varying by case complexity: larger or more complex scenarios yield greater time savings. This efficiency gain is expected to become more pronounced for highly complex, large-scale tasks.

(3) The construction of the spatial-temporal risk occupancy map constitutes a significant portion of computational time, with duration dependent on the risk assessment scope, grid density, and number of surrounding vehicles. Future work will explore leveraging edge computing within a "vehicle-road-cloud integration" framework to reduce computational resource consumption.
(4) The running time of dynamic iterative sampling is closely linked to the spatial-temporal risk occupancy map: as high-risk regions expand and reduce the number of viable sampling points, fewer valid sampling paths can be generated, leading to a corresponding decrease in sampling time.

**Table 7** Running time analysis for trajectory planning.

| Scenario / Module | On-ramp, simple and uncongested scenario (ms) | Off-ramp, simple and uncongested scenario (ms) | On-ramp, complex and congested scenario (ms) | Off-ramp, complex and congested scenario (ms) |
|---|---|---|---|---|
| Construction of STROM | 25 | 32 | 45 | 53 |
| Dynamic iterative sampling | 24 | 20 | 17 | 16 |
| Path evaluation | 3 | 3 | 2 | 2 |
| Path smoothing | 9 | 8 | 11 | 11 |
| Speed smoothing | 9 | 9 | 12 | 12 |
| Path-speed parallel smoothing | 15 | 15 | 18 | 18 |
| Total running time | 67 | 70 | 82 | 91 |

6.4 Sensitivity analysis of $R_{max}$

The maximum allowable risk $R_{max}$ is a key parameter in calibrating the spatial-temporal risk field model, theoretically influencing a vehicle's driving style. Intuitively, a higher $R_{max}$ corresponds to a more aggressive driving style—prioritizing efficiency but increasing risk exposure—while a lower $R_{max}$ aligns with a conservative style, ensuring safety at the potential cost of efficiency. Identifying an appropriate value is therefore critical for accurate risk quantification and effective trajectory planning. To validate the reasonableness of the calibrated parameters and explore how variations in this value affect driving style, we conducted a sensitivity analysis. Case 2 (from **Section 6.2.1**) was selected for this analysis, with **Fig. 30** illustrating the trajectories generated by our method under three distinct $R_{max}$ (3, 4, and 6). Simulation videos related to sensitivity analysis are available online, and the link can be found in **Appendix A**.

From these results, the following observations can be made:

(1) The calibrated parameter value of $R_{\max}$ is reasonable and valid. On the one hand, our calibrated parameter values avoid excessive transverse and longitudinal accelerations ($R_{\max} = 6$). On the other hand, they do not exhibit overly conservative behavior ($R_{\max} = 3$), which suffers from reduced efficiency as a result of such conservatism. This balance ensures compliance with comfort and safety constraints. This feasibility is rooted in extensive real-world data: during parameter calibration, we collected a large volume of real driving data and employed a rigorous calibration framework, enabling the parameters to effectively quantify driving risks in a manner that aligns with practical scenarios.

(2) Notably, as $R_{\max}$ increases, efficiency improves significantly. In specific scenarios, we can select different values of $R_{\max}$ to enable trajectory planning with varied driving styles, provided feasibility is fully demonstrated. For instance, a larger $R_{\max}$ may be chosen for cargo transportation where comfort is irrelevant or no occupants are present, while a smaller $R_{\max}$ is more suitable for scenarios involving special occupants. This aligns with the diverse driving styles observed in real traffic environments, offering a foundation for diversified, human-like autonomous driving.

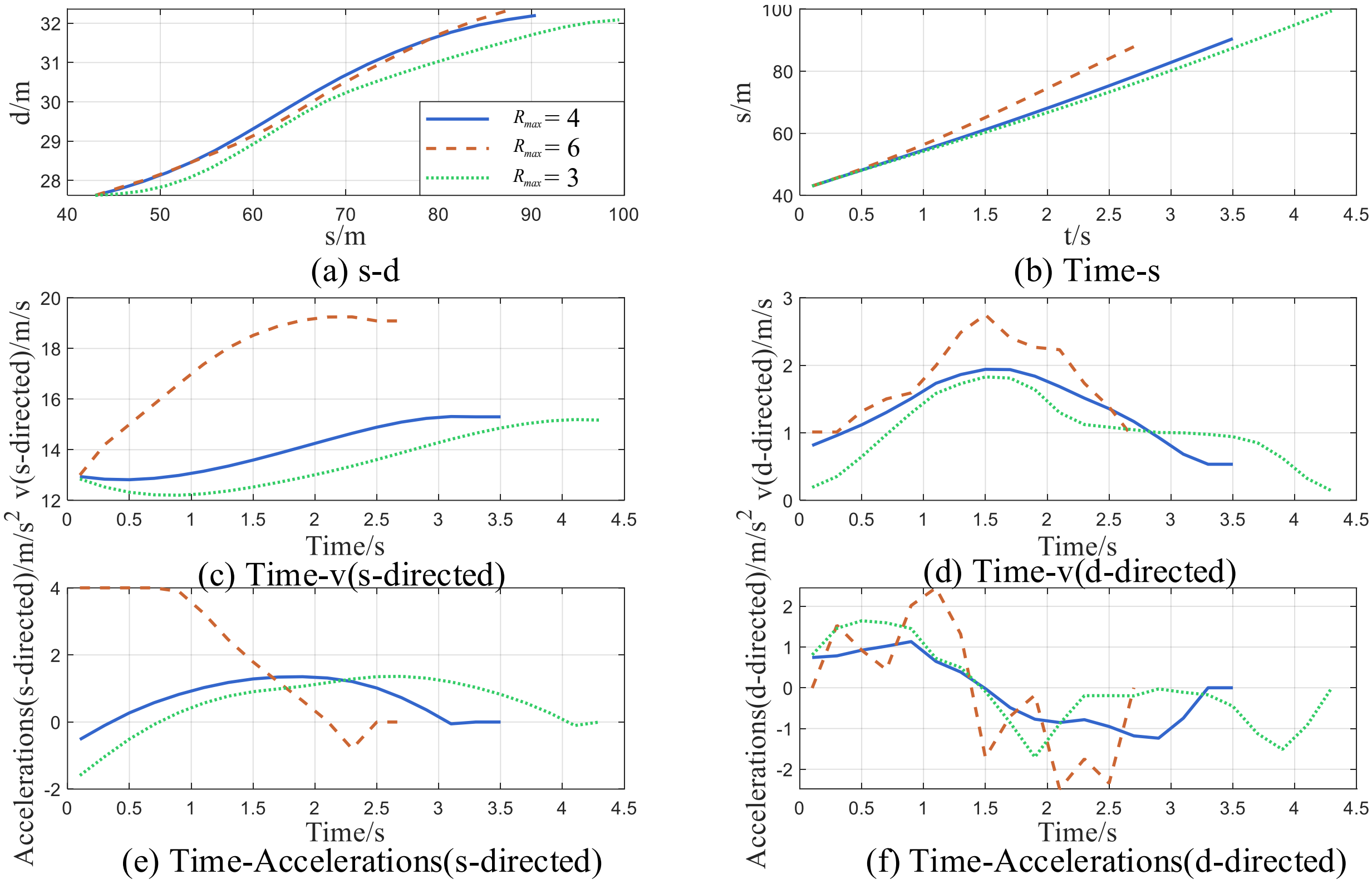

**Fig. 30** Sensitivity analysis of $R_{\max}$ (Case 2).

## 7 Conclusion

This paper addresses two key challenges: dynamic risk assessment in mixed traffic scenarios and trajectory planning for CAVs in weaving segments. First, we introduced a new paradigm of spatial-temporal distance, integrated the time dimension and future trajectories of dynamic obstacles, and constructed a three-dimensional STRF model, which overcomes the limitations of traditional two-dimensional static fields in capturing dynamic risks. In addition, we also imposed constraints on specific vehicles based on the geometric characteristics of the on-ramps and off-ramps at the weaving segments, significantly improving the accuracy of risk assessment. We calibrate STRF's parameters using real-world aerial video data, leveraging YOLO and dynamic risk balance theory. Experimental results demonstrate that the calibrated risk field effectively distinguishes between high- and low-risk areas while dynamically reflecting vehicle motion states. Compared with the traditional risk field, the STRF demonstrates excellent risk foreseeability. Second, we explore the application of STRF in CAV trajectory planning in weaving segments. The proposed trajectory planning framework integrates STROM, dynamic iterative sampling, and QP, effectively balancing safety, efficiency, and comfort. Meanwhile, this framework also supports a two-stage planning process that optimizes path and speed in parallel, reducing the algo-rithm's completion time to a certain extent. Experimental comparisons show that our method outperforms both the baseline and human-driven scheme in terms of efficiency, smoothness, and speed fluctuation.

Despite its strong performance in simulation experiments, this study still has several objective limitations, which also point toward directions for future research. First, regarding the limitations of data calibration, the parameter calibration of the STRF model in this study primarily relies on trajectory data from HDVs. Considering the objective physical differences between CAVs and human drivers in terms of perception precision, response latency, and underlying vehicle dynamics control, using human driving data to calibrate CAV model parameters may introduce certain errors. This could lead to a subtle deviation between the risk assessment results and the actual physical boundaries of CAVs. Second, concerning the limitations of system functionality, although this study is framed around CAVs, the current methodological framework primarily reflects the obstacle avoidance and planning functions of a single CAV. Since active communication and cooperative game mechanisms between vehicles have not been explicitly introduced into the model, the core advantages of CAV systems in system-level traffic flow optimization and cooperative conflict resolution have not yet been fully utilized. Therefore, future research will aim to extend STRF as a global optimization objective encompassing ramp control and multi-vehicle cooperative control. Finally, the conclusion that the average execution time of a single optimization cycle is less than 100 ms serves as a preliminary evaluation of the algorithm's efficiency under specific environmental conditions. It must be explicitly stated that, at the current stage, this study has not yet deployed the proposed models directly onto physical onboard hardware for real-world autonomous driving tests. Considering the objective physical differences between the laboratory MATLAB simulation environment and actual vehicle platforms, the practical response time in real-world engineering implementations and complex onboard application scenarios may exhibit slight deviations from the values reported herein. Therefore, future research aims to systematically optimize the relevant code and deploy it to onboard computing units for real-world testing, thereby further verifying its real-time performance.

### CRediT authorship contribution statement

**Guodong Ma:** Writing – original draft, Writing – review & editing, Software, Methodology, Conceptualization, Formal analysis. **Baofeng Sun:** Writing – review & editing, Supervision, Funding acquisition. **Hongchao Liang:** Methodology, Formal analysis. **Wenyu Yang:** Formal analysis. **Huxing Zhou:** Supervision.

### Declaration of competing interest

The authors declare that they have no known competing financial interests or personal relationships that could have appeared to influence the work reported in this paper.

### Acknowledgments

This research was supported by National Natural Science Foundation of China (grant number 52472313). The authors would like to acknowledge the editor and the anonymous reviewers, whose insightful comments and suggestions helped to improve the paper.

### Data availability

Data will be made available on request.

**Appendix A. Links to trajectory planning simulation video**

The links to the simulation-generated trajectory planning videos and sensitivity analysis videos are as follows: https://github.com/GuodongMa11/Trajectory-Planning-Simulation-Video/issues/1.

**Appendix B. Data used in this study**

The relevant data can be obtained by contacting the first author at: magd22@mails.jlu.edu.cn.